\documentclass[11pt]{article}
\usepackage{amsmath,amssymb,amsthm,array,cite,version,chngcntr}
\usepackage{graphicx,graphics,epsfig,subfigure,xcolor,color,float}
\usepackage{pifont} 
\usepackage{cleveref}
\crefformat{section}{\S} %
\allowdisplaybreaks

\theoremstyle{definition}

\newtheorem{alg}{Algorithm}[section]

\theoremstyle{plain}
\newtheorem{thm}{Theorem}[section]
\newtheorem{lem}{Lemma}[section]

\theoremstyle{remark}
\newtheorem{rem}{Remark}[section]
\newtheoremstyle{notes}%
{3pt}
{3pt}
{}
{}
{\itshape\color{red}}
{.}
{.5em}
{}
\theoremstyle{notes}

\numberwithin{equation}{section}
\numberwithin{figure}{section}
\numberwithin{table}{section}
\numberwithin{footnote}{section}
\counterwithout{alg}{section}
\counterwithout{nts}{section}
\counterwithout{sch}{section}

\newcommand{\bena}{\begin{eqnarray}\begin{array}{l}}
\newcommand{\eena}{\end{array}\end{eqnarray}}
\newcommand{\ben}{\begin{eqnarray}}
\newcommand{\een}{\end{eqnarray}}
\newcommand{\bea}{\begin{array}}
\newcommand{\eea}{\end{array}}
\newcommand{\bes}{\begin{subequations}}
\newcommand{\ees}{\end{subequations}}
\newcommand{\bec}{\begin{cases}}
\newcommand{\eec}{\end{cases}}
\newcommand{\bef}{\begin{figure}[H]}
\newcommand{\eef}{\end{figure}}
\newcommand{\bet}{\begin{tikzpicture}}
\newcommand{\eet}{\end{tikzpicture}}
\newcommand{\beq}{\begin{equation}}
\newcommand{\eeq}{\end{equation}}
\newcommand{\bep}{\begin{proof}}
\newcommand{\eep}{\end{proof}}
\def\beg#1\eeg{\begin{align}#1\end{align}}
\def\besl#1\eesl{\begin{subequations}\begin{align}#1\end{align}\end{subequations}}

\def\inc(#1){\includegraphics[height=3 cm]{pics/#1}}

\def\bn{\ensuremath{{\bf n}}}

\def\bv{\ensuremath{{\bf v}}}

\def\bI{\ensuremath{{\bf I}}}

\begin{document}
\title{A Second Order Fully-discrete Linear  Energy Stable Scheme for a Binary Compressible Viscous Fluid Model}
\author{Xueping Zhao\footnote{Department of Mathematics,
  University of South Carolina,
   Columbia, SC 29208, USA} and Qi Wang\footnote{Department of Mathematics,
  University of South Carolina,
   Columbia, SC 29208, USA and Beijing Computational Science Research Center, Beijing  100193, China}}

\date{\today}
\maketitle

\begin{abstract}

We present a linear, second order fully discrete numerical scheme  on a staggered grid for a thermodynamically consistent hydrodynamic phase field model of binary compressible fluid flow mixtures derived from the generalized Onsager Principle. The hydrodynamic model not only possesses the variational structure, but also warrants the mass, linear momentum conservation as well as  energy dissipation. We first reformulate the model in an equivalent form using the energy quadratization method and then discretize the reformulated model to obtain a semi-discrete partial differential equation system using the Crank-Nicolson method in time. The numerical scheme so derived preserves the mass conservation and energy dissipation law at the semi-discrete level. Then, we discretize the semi-discrete PDE system on a staggered grid in space to arrive at a fully discrete scheme using the 2nd order finite difference method, which respects a discrete energy dissipation law.  We prove the unique solvability of the linear system resulting from the fully discrete scheme. Mesh refinements and two numerical examples on phase separation due to the spinodal decomposition in two polymeric fluids and interface evolution in the gas-liquid mixture are presented to show the convergence property and the usefulness of the new scheme in applications.

\end{abstract}

\maketitle

\section{Introduction}
\noindent \indent Material systems comprising of multi-components, some of which are  compressible while others are incompressible,  are ubiquitous in nature and industrial applications. For example, in growing tissues, cell proliferation makes the material volume changes so that it can not be described as incompressible \cite{Hannaezo&Jonanny2016}.  Another example of the mixture of compressible fluids is the binary fluid flows of non-hydrocarbon (e.g. $CO_2$) and hydrocarbons encountered in the enhanced oil recovery(EOR) process. Since gas (e.g. $CO_2$) injection offers considerable potential benefits to oil recovery and is attracting the most new market interest since 1972, properties (viscosity, density et al.) of multi-component compressible mixtures of nonhydrocarbon and hydrocarbons have been studied by a number of investigators \cite{Cullick_Mathis1984, Hong_Duan03, hydrocarbons}.

Phase field methods have been used successfully to formulate models for fluid mixtures in many applications ranging from life sciences \cite{ShaoD1, ShaoD2, WiseJTB2008, AransonPLOSONE2013} (cell biology \cite{KapustinaPLOSCB2015, Najem&GrantPRE2016, ShaoD1, ZhaoWang-cell-divisionII, AransonJRSI2012, JulicherNP2016}, biofilms \cite{Zhao&WangMB2016, ZhaoWang-biofilm-datafit3D, Zhao&WangBMB2016}, cell adhesion and motility \cite{Camley&Zhao&Li&Levine&RappelPRE2017, AransonSoftMatter2014, Najem&GrantSM2015, Najem&GrantPRE2016,  NonomuraPLOSONE2012, ShaoD1}, cell membrane \cite{Aland&Lowengrub&VoigtJCP2014, AlandS1, GavishN1, WangX1, WitkowskiT1}, tumor growth \cite{WiseJTB2008}), to  materials science \cite{BordenCMAME2012, Chen&Y1994, DuLiRyWa052}, fluid dynamics \cite{LowengrubPRE2009, LowengrubRSA1998, TorabiS1}, image processing \cite{BertozziIEEE, Li&KimCMA2011, BertozziSINUMA2000} et al. The most widely studied phase field model for binary fluid mixtures is the one for fluid mixtures of  two incompressible fluids of identical densities \cite{Hohenberg&Halperin1977, Liu&Shen2003, Abels2009}. While modeling binary fluid mixtures using phase field models, one commonly uses a labeling or a phase variable (a volume fraction or a mass fraction) $\phi$ to distinguish between distinct material phases. For instance $\phi = 1$ indicates one fluid phase while $\phi = 0$ denotes another fluid phase in the immiscible, binary fluid mixture. For immiscible mixtures, the interfacial region is described by $0 < \phi < 1$. A transport equation for the volume fraction $\phi$ along with conservation equations of mass and momentum constitute the governing system of  equations for the binary incompressible fluid mixture.

In the compressible fluid flow, we use the mass density $\rho_i$ or molar density $n_i$ in place of volume fraction $\phi_i$ ($i =1,2$), to represent the distribution of each compressible component in the fluid mixture. The material compressibility comes from two sources. One is the material compressibility itself and the other is the mass-generating source. In general, the transport equation for the mass density of each component is given by
\ben\bea{l}
\frac{\partial \rho_i}{\partial t} + \nabla \cdot ({\rho_i \bf v_i}) = j_i, \qquad i = 1, \cdots, N,
\eea\een
or
\ben\bea{l}
\frac{\partial n_i}{\partial t} + \nabla \cdot ({n_i \bf v_i}) = j_i, \qquad i = 1, \cdots, N,
\eea\een
where ${\bf v}_i$ is the velocity of the ith component, $j_i$ is the mass source or molar source of the ith component. The transport equations for the mass or molar densities along with the conservation laws of mass and momentum constitute the governing equations of the hydrodynamic phase field models of the compressible fluid mixtures.

Distinguishing properties of the compressible hydrodynamic phase field models include that the density of each compressible material component is a variable, the mass average velocity of the fluid flow is most likely not solenoidal, and the pressure is determined by the equation of state or the free energy of the mixture system (at least in the isothermal case). In \cite{LowengrubRSA1998}, Truskinovsky and Lowengrub derived the Navier$-$Stokes$-$Cahn$-$Hilliard (NSCH) system for a binary mixture of two incompressible fluid flows with unmatched densities in the fluid components, in which the mass concentration of one fluid component in the binary fluid flow is used as the phase variable. They termed the hydrodynamic phase field model quasi-incompressible.  In \cite{SunsyMulti, SunsyMulti2}, Sun et al. propose a general diffuse interface model with a given equation of state (e.g. Peng-Robinson equation of state) to describe the multi-component fluid flow based on the principles of the NVT-based framework. In \cite{Zhao&W2018}, we systematically derived a thermodynamically consistent hydrodynamic phase field model for multi-component compressible fluid mixtures through a variational approach coupled with the generalized Onsager Principle \cite{Yang&Li&Forest&Wang2016} and discussed various means to arrive at the quasi-incompressible limit and the fully incompressible limit. In this paper, we develop an unconditionally  energy stable numerical algorithm to solve the thermodynamically consistent, hydrodynamic phase field model.

The hydrodynamic phase field model is nonlinear, exemplified  in its free energy, mobility coefficients and in the advection in the transport equations. Higher order approximation, unconditional energy stability as well as computational efficiency are desired properties to attain  in developing its  numerical approximation. To preserve the energy dissipation property, several time-marching approaches have been developed in the past: convex splitting method \cite{Ell.S93, Eyre1998, WQF2018_1, WQF2018_2}, stabilization method \cite{ZhuJ1, Shen10_1}, and energy quadratization (EQ, including SAV) approach \cite{ZYGW_LCP_2017, ZYLW_SISC_2016, Gong&Z&Y&W2018_SIAM_JSC, Gong&Z&Y&W2018_SIAM_JSC}.The convex splitting method has been used to obtain a series of first order energy stable schemes for various PDE models exhibiting energy dissipation properties. However, the convex-splitting scheme is usually nonlinear and therefore can be expensive to solve from time to time. On the other hand, even though it is possible to construct a second order convex splitting scheme in some cases, it was usually done on a case by case basis and a general formulation is not yet available. The stabilization method is another method for obtaining energy stable numerical approximations, which is equivalent to a   convex splitting method in some cases. By adding a linear, stablizing operator in the order of the truncation error, one can obtain an energy stable algorithm. In general, a second order stabilizing scheme  can be derived, it preserves the discrete energy decay but not the dissipation rate. The energy quadratization(EQ), also known as the invariant energy quadratization(IEQ),  method was proposed recently \cite{YANG16} and well developed in various gradient flows and  hydrodynamic phase field models \cite{ZYGW_LCP_2017, ZYLW_SISC_2016, Gong&Z&Y&W2018_SIAM_JSC}. By introducing intermediate variables, one can rewrite the nonlinear free energy functional into a quadratic from, from which a linear second order or even higher order numerical scheme can be constructed \cite{Zhao&Yang&Gong&WangCMAME2017, Gong&Zhao&WangCPC, Jia_XF_YZ_XP_XG_Jun_Q2018}.

Recently, Sun et al. \cite{SunsyMulti, SunsyMulti2} used the convex splitting approach and the scalar auxiliary variable method \cite{SAV_SHEN2018407}, which is developed based on the EQ strategy, to solve binary compressible hydrodynamic phase field models, respectively. They obtained some first order semi-discrete schemes.
In this paper, we develop a linear, second order, fully discrete  numerical scheme for the hydrodynamic phase field model for binary fluid mixtures based on the energy quadratization strategy. We will show that this scheme is unconditionally energy stable and the linear system resulting from the second order numerical scheme is uniquely solvable.  At each time step, the linear algebraic system is solved within 3 iterations, with a linear pre-conditioner. Two examples on phase separation dynamics in viscous polymeric blends and interface evolution in gas-liquid mixtures are presented to show the usefulness of the new scheme in some practical applications.

The paper is organized as follows. In \S 2, we  briefly recall the derivation of the compressible hydrodynamic phase field model. Its non-dimensionalization is given in \S 3. In \S 4, we reformulate the model using the energy quadratization method. The fully discrete numerical scheme, where we use second order finite difference in space and "linearized" Crank-Nicolson method in time, is given in \S 5 where the unique solvability of the scheme and the property of energy dissipation are proved as well. In \S 6, we show several numerical experiments that validate the accuracy, stability and efficiency of the numerical scheme. We give concluding remarks in \S 7.


\section{Thermodynamically Consistent Hydrodynamic Phase Field  Models for Binary Compressible Viscous Fluid Flows}

\noindent \indent A general thermodynamically consistent hydrodynamic phase field model for fluid mixture of n viscous fluid components has been derived in \cite{Zhao&W2018}. Here, we brief recall the basic ingredients  in the binary fluid model and discuss its energy dissipation property. We consider a fluid mixture flow of two compressible viscous fluids with densities $\rho_1$ and $\rho_2$, respectively. The mass conservation equation for each fluid component is respectively given by
\ben
\frac{\partial \rho_i}{\partial t} + \nabla \cdot(\rho_i {\bf v}_i)  = 0, \quad i=1,2,
\een
where ${\bf v}_i$ is the velocity of the ith fluid component, $i=1,2$.
We define the total mass of the fluid mixture as $\rho=\rho_1+\rho_2$ and the mass average velocity as
$
{\bf v}=\frac{1}{\rho}(\rho_1{\bf v}_1+\rho_2{\bf v}_2).
$
Then, the mass conservation equation for the total mass density $\rho$ is given by
\ben\bea{l}
\frac{\partial \rho}{\partial t} + \nabla \cdot(\rho {\bf v})  = 0.
\eea \label{eq:mass_conservation}
\een
Using the mass average velocity, we rewrite the mass transport equation as follows
\ben\bea{l}
\frac{\partial \rho_i}{\partial t} + \nabla \cdot(\rho_i {\bf v})  = j_i = \nabla \cdot  {\bf J}_i, \quad i = 1,2,
\eea \label{eq:transport_eqn}
\een
where  ${\bf J}_i=\rho_i({\bf v} - {\bf v}_i)$ is the excessive mass flux of fluid $i=1,2$, and $j_1 + j_2 = 0$ according to the total mass conservation law.
The linear momentum conservation law of the fluid mixture is given by
\ben\bea{l}
\frac{\partial (\rho {\bf v})}{\partial t}  +  \nabla \cdot (\rho {\bf v} {\bf v})  = \nabla \cdot \sigma+{\bf b}
\eea \label{eq:momentum_conservation}
\een
from the momentum conservation for each fluid component,
where $\bf b$ is the body force, $\sigma$ is the total stress tensor, $\sigma = \sigma^s + \sigma^e$,  $\sigma^s$ is the symmetric viscous stress tensor, and $\sigma^e$ is the Ericksen stress tensor. Both ${\bf J}_i$, i = 1,2 and $\sigma^s$ would be determined  by constitutive relations later.

For the compressible fluid mixture, we assume the free energy of the system is given by
\ben\bea{l}
F = \int_V f(\rho_1, \rho_2, \nabla \rho_1, \nabla \rho_2) d{\bf x},
\eea\een
where $f$  is the free energy density function and $V$ is the domain in which the fluid mixture occupies.
The total energy of the fluid system is given by the sum of the kinetic energy and the free energy:
\ben\bea{l}
E_{total} = \int_V [\frac{1}{2}\rho ||{\bf v}||^2 + f ]d{\bf x}.
\eea\een

Considering the conservation laws of mass and linear momentum,  we calculate the  energy dissipation rate as follows
\ben\bea{l}
\frac{dE_{total}}{dt}
=  \int_{V}[- \sigma^s : {\bf D}   +  ({\bf b}+\nabla \cdot \sigma^e + \rho_1 \nabla \mu_1 + \rho_2 \nabla \mu_2) \cdot {\bf v}  +   \mu_1  j_1 + \mu_2 j_2 ] d{\bf x} \\
+ \int_{\partial V} [(\sigma^s \cdot {\bf v}) \cdot {\bf n}   - \frac{1}{2} (\rho {\bf v} \|{\bf v}\|^2) \cdot {\bf n}
+   ( -  \mu_1 \rho_1 {\bf v} -  \mu_2 \rho_2 {\bf v}  +\frac{\partial f}{\partial (\nabla \rho_1)}\frac{\partial \rho_1}{\partial t}    +   \frac{\partial f}{\partial (\nabla \rho_2)}\frac{\partial \rho_2}{\partial t} ) \cdot {\bf n}] dS.
\eea\een
where   ${\bf D} = \frac{1}{2}(\nabla {\bf v} + \nabla {\bf v}^T)$ is the rate of strain tensor, $\bf n$ is the unit external normal of the domain boundary $\partial V$,
$\mu_1 = \frac{\partial f}{\partial \rho_1} - \nabla \cdot \frac{\partial f}{\partial \nabla \rho_1}$, $\mu_2= \frac{\partial f}{\partial \rho_2} - \nabla \cdot \frac{\partial f}{\partial \nabla \rho_2}$ are the chemical potentials with respect to $\rho_1$ and $\rho_2$, respectively.
 We identify the Erickson stress by the equation
\ben\bea{l}
\nabla \cdot \sigma^e = -\rho_1 \nabla \mu_1 - \rho_2  \nabla \mu_2.
\eea\een
The energy dissipation rate reduces to
\ben\bea{l}
\frac{dE_{total}}{dt}
=  \int_{V}[{\bf b}\cdot {\bf v}- \sigma^s :{\bf D}  +  \mu_1 j_1 + \mu_2 j_2]  d{\bf x}+\int_{\partial V} [(\sigma^s \cdot {\bf v}) \cdot {\bf n}   - \frac{1}{2} (\rho {\bf v} \|{\bf v}\|^2) \cdot {\bf n}\\
+   ( -  \mu_1 \rho_1 {\bf v} -  \mu_2 \rho_2 {\bf v}  +\frac{\partial f}{\partial (\nabla \rho_1)}\frac{\partial \rho_1}{\partial t}    +   \frac{\partial f}{\partial (\nabla \rho_2)}\frac{\partial \rho_2}{\partial t} ) \cdot {\bf n}] dS.
\eea\een
In the bulk integral, we propose the following constitutive relations following the generalized Onsager principle
\ben\bea{l}
\sigma^s =  2 \eta {\bf D} + \overline{\eta} tr({\bf D}){\bf I},\\
j_i = -\sum_{k=1}^2 \nabla \cdot M_{ik} \cdot \nabla \mu_{k},
\eea\een
where $\eta$, $\bar{\eta}$ are the shear and volumetric viscosity respectively, and $\mathcal{M} = (M_{ik})_{2\times2} \geq 0$ is the symmetric mobility matrix. Since  $\sum_{i=1}^2 j_i = 0$ according to the mass conservation law, this imposes a constraint $\mathcal{M}\cdot {\bf 1}={\bf 0},$ where ${\bf 1}^T=(1,1)$.
Examining the surface integral, we notice that  if we assume the following conditions
\ben\bea{l}
{\bf v}|_{\partial V}=0, \qquad {\bf n}  \cdot \nabla \mu_i|_{\partial V} = 0, \qquad {\bf n} \cdot \frac{\partial f}{ \partial (\nabla \rho_i)}|_{\partial V} = 0, \quad i = 1, 2.
\eea \label{pdes-bc}
\een
on the boundary, the surface integral vanishes in the energy dissipation function. So, at the absence of the body force ${\bf b}={\bf 0}$, the total energy dissipation rate reduces to
\ben\bea{l}
\frac{dE_{total}}{dt}
=  -\int_{V}[ 2\eta {\bf D} :{\bf D} +\overline{\eta} tr({\bf D})^2  + (\nabla \mu_1, \nabla \mu_2) \cdot \mathcal{M} \cdot (\nabla \mu_1, \nabla \mu_2)^T ]  d{\bf x}\leq 0,
\eea\een
provided $\eta, \overline{\eta}\geq 0, \mathcal{M} \geq  0$.

\begin{rem}
If we choose the boundary conditions as follows
\ben\bea{l}
\bv\cdot \bn=0, \qquad  \sigma_s\cdot \bn=-\beta (\bI-\bn\bn)\cdot \bv, \qquad  \bn\cdot \frac{\partial f}{\partial (\nabla \rho_1)}=-\gamma_1\frac{\partial \rho_1}{\partial t}, \qquad  \bn\cdot \frac{\partial f}{\partial (\nabla \rho_2)}=-\gamma_2\frac{\partial \rho_2}{\partial t},
\eea\een
where $\beta, \gamma_1, \gamma_2 \geq 0$, the energy dissipation rate is given by
\ben\bea{l}
\frac{dE_{total}}{dt}
=  -\int_{V}[ 2\eta {\bf D} :{\bf D} +\overline{\eta} tr({\bf D})^2  + (\nabla \mu_1, \nabla \mu_2) \cdot \mathcal{M} \cdot (\nabla \mu_1, \nabla \mu_2)^T ]  d{\bf x}\\

-\int_{\partial V}[   { {\beta (\bI-\bn\bn)  \| \bv\|^2 }} +\gamma_1 (\frac{\partial \rho_1}{\partial t})^2+\gamma_2  (\frac{\partial \rho_2}{\partial t})^2 ] ds.
\eea\een
These boundary conditions allow fluid flows slip at the boundary and mass fluxes to move through the boundary, which leads to additional  energy dissipation due to energy dissipation at the surface. We will not pursue these boundary conditions in this study, which worthy of a complete study of its own.
\end{rem}

We summarize the governing equations of the compressible binary fluid system in the hydrodynamic phase field model as follows:
\ben\bea{l}
\begin{cases}
\frac{\partial \rho_1}{\partial t} + \nabla \cdot (\rho_1 {\bf v}) =   \nabla \cdot M_{11} \cdot \nabla \mu_1 +  \nabla \cdot M_{12} \cdot \nabla \mu_2,\\

\frac{\partial \rho_2}{\partial t} + \nabla \cdot (\rho_2 {\bf v}) =   \nabla \cdot M_{21} \cdot \nabla \mu_1 +  \nabla \cdot M_{22} \cdot \nabla \mu_2,\\

\frac{\partial (\rho {\bf v})}{\partial t}  +  \nabla \cdot (\rho {\bf v} {\bf v})  = 2 \nabla \cdot ( \eta{\bf D}) + \nabla ( \overline{\eta} \nabla \cdot {\bf v})-  \rho_1 \nabla \mu_1 - \rho_2 \nabla \mu_2,
\end{cases}
\eea
\een
where
 $\sum_{i,k=1}^2 \nabla \cdot M_{ik} \cdot \nabla \mu_{k} = 0$. One particular mobility matrix satisfying the constraint is consisted of the entries $M_{1} = M_{11} = -M_{12} = -M_{21} = M_{22}$. The governing equations reduce to
\ben\bea{l}
\begin{cases}
\frac{\partial \rho_1}{\partial t} + \nabla \cdot (\rho_1 {\bf v}) =   \nabla \cdot M_{1} \cdot \nabla (\mu_1 - \mu_2),\\

\frac{\partial \rho_2}{\partial t} + \nabla \cdot (\rho_2 {\bf v}) =   -\nabla \cdot M_{1} \cdot \nabla (\mu_1 - \mu_2),\\

\frac{\partial (\rho {\bf v})}{\partial t}  +  \nabla \cdot (\rho {\bf v} {\bf v})  = 2 \nabla \cdot ( \eta{\bf D}) + \nabla ( \overline{\eta} \nabla \cdot {\bf v})-  \rho_1 \nabla \mu_1 - \rho_2 \nabla \mu_2.
\end{cases}
\eea \label{pdes}
\een
For the viscosity coefficients, we denote $\eta_1,\eta_2$ as the shear viscosities of the fluid component 1 and 2 respectively, and $\overline{\eta}_1, \overline{\eta}_2$ as the volumetric viscosities of the two components. $\eta, \overline{\eta}$ are chosen as the mass average  viscosities of the two components:
\ben
\eta=\frac{1}{\rho}[\rho_1 \eta_1+\rho_2 \eta_2], \qquad \overline{\eta}=\frac{1}{\rho}[\rho_1 \overline{\eta}_1+\rho_2 \overline{\eta}_2].
\een

In this study, we focus on the free energy density function $f$ in the following form
\ben
f(\rho_1, \rho_2, \nabla \rho_1, \nabla \rho_2, T) = h(\rho_1, \rho_2, T) + \frac{1}{2}[\kappa_{\rho_1 \rho_1} (\nabla \rho_1)^2 + 2\kappa_{\rho_1 \rho_2} (\nabla \rho_1, \nabla \rho_2) + \kappa_{\rho_2 \rho_2} (\nabla \rho_2)^2].
\een
where $h(\rho_1, \rho_2, T)$ is the homogeneous or the bulk free energy density function, T is the  absolute temperature, assumed a constant in this study, and $\kappa_{\rho_i \rho_j}, i, j = 1, 2$ are model parameters measuring the strength of the conformational entropy (which are assumed constant in this study).

Sometimes, we have to use molar densities  $n_i$  as the fundamental variables in the model $i=1,2$, system (\ref{pdes}) can be rewritten as follows
\ben\bea{l}
\begin{cases}
m_1 (\frac{\partial n_1}{\partial t} + \nabla \cdot (n_1 {\bf v}) )=    \nabla \cdot   M_{1} \cdot \nabla ( \frac{1}{m_1}  \mu_{n1} -  \frac{1}{m_2}  \mu_{n2}),\\

m_2 (\frac{\partial n_2}{\partial t} + \nabla \cdot (n_2 {\bf v}) )=   - \nabla \cdot  M_{1} \cdot \nabla ( \frac{1}{m_1} \mu_{n1} -  \frac{1}{m_2} \mu_{n2}),\\

\frac{\partial (\rho {\bf v})}{\partial t}  +  \nabla \cdot (\rho {\bf v} {\bf v})  = 2 \nabla \cdot ( \eta{\bf D}) + \nabla ( \overline{\eta} \nabla \cdot {\bf v})-   n_1 \nabla \mu_{n1} -  n_2 \nabla \mu_{n2},
\end{cases}
\eea \label{pdes_mole}
\een
where $n_i = \frac{\rho_i}{m_i}$, $m_i$ is the molar mass of the ith component and $\mu_{ni} = \frac{\delta f}{\delta n_i} = \frac{\delta f}{\delta \rho_i} m_i$, i = 1, 2. Correspondingly, The shear and volumetric viscosities are given respectively by $\eta =\sum_{i=1}^2 \frac{n_i m_i}{n_1 m_1 + n_2 m_2} \eta_i $ and $\overline{\eta} = \sum_{i=1}^2 \frac{n_i m_i}{n_1 m_1 + n_2 m_2} \overline{\eta_i}$.

With  molar densities $n_i, i=1,2$ as the primitive variables, we rewrite free energy density $f$ as follows
\ben\bea{l}
f(n_1 m_1, n_2 m_2,  m_1 \nabla n_1, m_2 \nabla n_2, T) = h(m_1 n_1, m_2  n_2, T) \\
+ \frac{1}{2}[\kappa_{n_1 n_1} (\nabla  n_1)^2 + 2\kappa_{n_1 n_2} (\nabla n_1, \nabla n_2) + \kappa_{n_2 n_2} (\nabla n_2)^2].
\eea\een
Where $\kappa_{n_i n_i} = m_i^2\kappa_{\rho_i \rho_i}$, i = 1, 2 and $\kappa_{n_1 n_2} = m_1 m_2 \kappa_{\rho_1 \rho_2}$.

The free energy density function is specific to the fluid system studied.
\begin{itemize}
\item  For polymeric binary fluid mixtures while approximated as a viscous fluid, the Flory-Huggins type free energy density function can be used to describe fluid mixing \cite{Doi&E1986, LowengrubRSA1998}
\ben\bea{l}
h(\rho_1, \rho_2, T) = \frac{k_B T}{m}\rho( \frac{1}{N_1}  \frac{\rho_1}{\rho} ln \frac{\rho_1}{\rho}+\frac{1}{N_2}   \frac{\rho_2}{ \rho} ln \frac{\rho_2}{\rho}+\chi \frac{\rho_1\rho_2}{\rho^2}),
\eea\label{eq:flory-huggins}
\een
Where $k_B$ is the Boltzmann constant, T is the absolute temperature and m the average mass of a molecule.

\item For compressible gas-liquid mixtures, the semi-empirical Peng-Robinson free energy density is often used \cite{SunsyMulti},
\ben\bea{l}
h(n_1, n_2, \cdots, n_N, n , T) = f^{ideal} + f^{repulsion} + f^{attraction},
\eea \label{eq:P_R}
\een
where
\ben\bea{l}
f^{ideal} = R T \sum_{i=1}n_i (ln (n_i) - 1), \\
f^{repulsion} = - n R T ln(1-bn),\\
f^{attraction} = \frac{a(T)n}{2 \sqrt{2} b} ln \big (  \frac{1 + (1- \sqrt{2} )b n}{ 1 + (1+ \sqrt{2} )b n}  \big).
\eea\een
Here $n = \sum_{i=1}^N n_i$ is the total molar density. The corresponding chemical potential of the ith component is given by
\ben\bea{l}
\mu_{ni} = \frac{\partial h}{\partial n_i} - \nabla \cdot \frac{\partial h}{\partial \nabla n_i} = R T \big (ln(n_i) + \frac{b_i n}{1- bn}  - ln(1- bn) \big ) + \frac{a b_i n}{b ((\sqrt{2}-1)bn - 1)  (1 + (1 + \sqrt{2})bn)}\\
+ \frac{1}{2\sqrt{2}} (\frac{2 \sum_{j=1}^M n_j (a_i a_j)^{1/2} (1 - k_{i j})  }{b n}  -   \frac{a b_i}{b^2} ) ln (\frac{1 + (1 - \sqrt{2})bn}{1 + (1 + \sqrt{2})bn}) - \kappa_{n_i n_i}\Delta n_i - \kappa_{n_i n_j} \Delta n_j, \qquad j \neq i,
\eea\een
where  $b(n_1, n_2)$ is the volume parameter and $a(n_1, n_2, T)$ is the interaction parameter. This free energy was proposed to improve that of the Van der Waals' to mitigate the deviation away from the ideal gas model.
\end{itemize}

\section{Non-dimensionalization}
\noindent \indent For system (\ref{pdes}), using characteristic time $t_0$, characteristic length $l_0$, and characteristic density $\rho_0$, we nondimensionalize the physical variables and parameters as follows
\ben\bea{l}
\tilde{t} = \frac{t}{t_0}, \quad \tilde{x} = \frac{x}{l_0}, \quad \tilde{\rho}_i = \frac{\rho_i}{\rho_0}, \quad i = 1,2, \quad \tilde{{\bf v}} = \frac{{\bf v}t_0}{l_0},  \quad \tilde{M}_1  = \frac{M_1 }{t_0 \rho_0},
\quad \frac{1}{{Re}_s} = \tilde{\eta} = \frac{t_0}{\rho_0 l_0^2} \eta, \\
\frac{1}{{Re}_v} = \tilde{\overline{\eta}} = \frac{t_0}{\rho_0 l_0^2} \overline{\eta}, \quad \tilde{\mu}_i = \frac{t_0^2}{ l_0^2} \mu_i, \quad i = 1,2, \quad \tilde \kappa_{\rho_i \rho_j} = \kappa_{\rho_i \rho_j} \frac{\rho_0 t_0^2}{l_0^4}, \quad i, j = 1,2,
\eea\een
where $Re_s$, $Re_v$ are the Reynolds numbers.  We rewrite the dimensionless governing equations, after dropping the $\tilde{}$s for simplicity, as follows
\ben\bea{l}
\begin{cases}
\frac{\partial \rho_1}{\partial t} + \nabla \cdot (\rho_1 {\bf v}) =   \nabla \cdot M_{1} \cdot \nabla ( \mu_1 -  \mu_2 ),\\

\frac{\partial \rho_2}{\partial t} + \nabla \cdot (\rho_2 {\bf v}) =   -\nabla \cdot M_{1} \cdot \nabla ( \mu_1 -  \mu_2 ),\\

\frac{\partial (\rho {\bf v})}{\partial t}  +  \nabla \cdot (\rho {\bf v} {\bf v})  = 2 \nabla \cdot ( \frac{1}{Re_s} {\bf D}) + \nabla ( \frac{1}{Re_v} \nabla \cdot {\bf v})-  \rho_1 \nabla \mu_1 - \rho_2 \nabla \mu_2.
\end{cases}
\eea\label{eq:nondim}
\een
where
\ben
\mu_1 =\frac{\partial h}{\partial \rho_1} - \kappa_{\rho_1 \rho_1} \Delta \rho_1 - \kappa_{\rho_1 \rho_2} \Delta \rho_2, \qquad
\mu_2 = \frac{\partial h}{\partial \rho_2} - \kappa_{\rho_1 \rho_2} \Delta \rho_1 - \kappa_{\rho_2 \rho_2} \Delta \rho_2.
\een

\noindent \indent Similarly, for system (\ref{pdes_mole}) with molar density as fundamental variables, using  characteristic molar density $n_0$ ($mol \cdot m^{-d}$), characteristic mass density $\rho_0 = n_0 m_2$($kg \cdot m^{-d}, d=3$) and characteristic temperature $T_0$ (Kelvin), we nondimensionalize the physical variables and parameters as follows
\ben\bea{l}
\tilde{t} = \frac{t}{t_0}, \quad \tilde{x} = \frac{x}{l_0}, \quad \tilde{\rho}= \frac{\rho}{\rho_0}, \quad \tilde{n} = \frac{n}{n_0}, \quad \tilde{T} = \frac{T}{T_0},  \quad
 \frac{1}{{Re}_s} = \tilde{\eta} = \frac{t_0}{\rho_0 l_0^2} \eta, \quad \tilde{m_1} = \frac{m_1 n_0}{\rho_0}, \quad \tilde{m_2} = \frac{m_2 n_0}{\rho_0}, \\
\frac{1}{{Re}_v} = \tilde{\overline{\eta}} = \frac{t_0}{\rho_0 l_0^2} \overline{\eta}, \quad  \tilde{\mu}_{ni} = \frac{n_0 t_0^2}{\rho_0 l_0^2} \mu_{ni},\quad i = 1,2, \quad
\tilde{M}_1  =  \frac{M_1 }{t_0 \rho_0} ,\quad \tilde{\kappa_{n_i n_j}} = \kappa_{n_i n_j} \frac{n_0^2 t_0^2}{\rho_0  l_0^4}, \quad i, j = 1,2.
\eea\een
Dropping $\tilde{}$s for simplicity, we rewrite the dimensionless governing equations as follows
\ben\bea{l}
\begin{cases}
m_1 (\frac{\partial n_1}{\partial t} + \nabla \cdot (n_1 {\bf v}) )=    \nabla \cdot   M_{1} \cdot \nabla ( \frac{1}{m_1}  \mu_{n1} -   \mu_{n2}),\\

 (\frac{\partial n_2}{\partial t} + \nabla \cdot (n_2 {\bf v}) )=   - \nabla \cdot  M_{1} \cdot \nabla ( \frac{1}{m_1} \mu_{n1} -   \mu_{n2}),\\

\frac{\partial (\rho {\bf v})}{\partial t}  +  \nabla \cdot (\rho {\bf v} {\bf v})  = 2 \nabla \cdot ( \eta{\bf D}) + \nabla ( \overline{\eta} \nabla \cdot {\bf v})-   n_1 \nabla \mu_{n1} -  n_2 \nabla \mu_{n2}.
\end{cases}
\eea\label{eq:nondim_mole}
\een
where we set $\tilde{m_2} = \frac{m_2 n_0}{\rho_0} = 1$, i.e. $m_1$ is the ratio of the specific masses, a dimensionless model parameter. The dimensionless chemical potentials are given by
\ben
\mu_{n1} =\frac{\partial h}{\partial n_1} - \kappa_{n_1 n_1} \Delta n_1 - \kappa_{n_1 n_2} \Delta n_2,\qquad
\mu_{n2} = \frac{\partial h}{\partial n_2} - \kappa_{n_1 n_2} \Delta n_1 - \kappa_{n_2 n_2} \Delta n_2.
\een

In the following, we focus on  developing an energy stable numerical scheme for system \eqref{eq:nondim} on staggered grids. An energy stable numerical scheme for system \eqref{eq:nondim_mole} can be obtained analogously. First, we reformulate the equation system using the energy quadratization strategy.

\section{Reformulation of the  Model using Energy Quadratization}

\noindent \indent  In order to use the Energy Quadratization (EQ) method to design numerical schemes, we need to reformulate the model equations.
We first transform the energy of the system into  a quadratic form
\ben\bea{l}
E_{total} = \int_V [\frac{1}{2}\rho {\bf v}^T {\bf v} + f] d{\bf x} =  \int_V [\frac{1}{2} {\bf u}^T {\bf u} +  q_1^2 +  \frac{1}{2} {\bf p}^T \cdot {\bf K} \cdot {\bf p} - A ]  d{\bf x}.
\eea\een
where ${\bf u} = \sqrt{\rho}{\bf v}$, $q_1 = \sqrt{h(\rho_1, \rho_2, T) + A}$  and A is a constant such that $h(\rho_1, \rho_2, T) + A > 0$. We note that we can always find a constant A if the bulk free energy density function is bounded below. In addition, ${\bf p} = (\nabla \rho_1, \nabla \rho_2)^T$ and $\bf K$ is the coefficient matrix of the conformational entropy
\ben
{\bf K} = \left(
\bea{cc}
\kappa_{\rho_1 \rho_1 }  & \kappa_{\rho_1 \rho_2}  \\
\kappa_{\rho_1 \rho_2}  & \kappa_{\rho_2 \rho_2 }
\eea
\right) > 0.
\een
Using identity
\ben\bea{l}
\frac{\partial (\sqrt{\rho}{\bf u})}{\partial t}  = \frac{1}{2\sqrt{\rho}} \frac{\partial \rho}{\partial t} {\bf u} + \sqrt{\rho} \frac{\partial {\bf u}}{\partial t}= -\frac{1}{2\sqrt{\rho}}\nabla \cdot (\sqrt{\rho} {\bf u}) {\bf u} + \sqrt{\rho} \frac{\partial {\bf u}}{\partial t},
\eea\een
we rewrite the governing equations  into
\ben\bea{l}
\begin{cases}
\frac{\partial \rho_1}{\partial t} + \nabla \cdot (\frac{\rho_1}{\sqrt{\rho}} {\bf u}) =   \nabla \cdot M_1 \cdot \nabla  (\mu_1 -  \mu_2)  ,\\

\frac{\partial \rho_2}{\partial t} + \nabla \cdot (\frac{\rho_2}{\sqrt{\rho}} {\bf u}) =   -\nabla \cdot M_1 \cdot \nabla  (\mu_1 -  \mu_2)  ,\\

\frac{\partial {\bf u}}{\partial t}   + \frac{1}{2}( \frac{1}{\sqrt{\rho}}\nabla \cdot ({\bf u}{\bf u}) + {\bf u} \cdot \nabla \frac{\bf u}{\sqrt{\rho}})  = \frac{1}{\sqrt{\rho}} \nabla \cdot \sigma ,\\

\frac{\partial q_1}{\partial t}= {\frac{\partial q_1}{\partial \rho_1}} \frac{\partial \rho_1}{\partial t} + {\frac{\partial q_1}{\partial \rho_2}}  \frac{\partial \rho_2}{\partial t},
\end{cases}
\eea\label{Model2}
\een
where
\ben\bea{l}
\sigma = \sigma^s + \sigma^e, \quad
\sigma^s =  2\frac{1}{Re_s}  {\bf D} + \frac{1}{Re_v}  (\nabla \cdot \frac{\bf u}{\sqrt{\rho}}){\bf I},\\

\sigma^e =  (f - \rho_1 \mu_1 - \rho_2 \mu_2) {\bf I} - \frac{\partial f}{\partial \nabla \rho_1} \nabla \rho_1 - \frac{\partial f}{\partial \nabla \rho_2} \nabla \rho_2,\\

\nabla \cdot \sigma = \nabla \cdot (\sigma^s + \sigma^e) =  2 \nabla \cdot ( \frac{1}{Re_s} {\bf D}) + \nabla ( \frac{1}{Re_v}  \nabla \cdot \frac{\bf u}{\sqrt{\rho}}) -  \rho_1 \nabla \mu_1 - \rho_2 \nabla \mu_2,\\

\mu_1 = \frac{\delta f}{\delta \rho_1 } = \frac{\partial f}{\partial \rho_1} - \nabla \cdot \frac{\partial f}{\partial \nabla \rho_1} = 2q_1 {\frac{\partial q_1}{\partial \rho_1}} - \kappa_{\rho_1 \rho_1} \Delta \rho_1 - \kappa_{\rho_1 \rho_2} \Delta \rho_2,\\

\mu_2 = \frac{\delta f}{\delta \rho_2 } = \frac{\partial f}{\partial \rho_2} - \nabla \cdot \frac{\partial f}{\partial \nabla \rho_2} = 2q_1 {\frac{\partial q_1}{\partial \rho_2}} - \kappa_{\rho_2 \rho_2} \Delta \rho_2 - \kappa_{\rho_1 \rho_2} \Delta \rho_1,\\

{\bf D} = \frac{1}{2} (\nabla \frac{\bf u}{\sqrt{\rho}} + (\nabla \frac{\bf u}{\sqrt{\rho}})^T), \quad

 \frac{1}{Re_s}  = \frac{\rho_1}{\rho}  \frac{1}{Re_{s1}}  + \frac{\rho_2}{\rho}  \frac{1}{Re_{s2}} , \quad

 \frac{1}{Re_v}  = \frac{\rho_1}{\rho} \frac{1}{Re_{v1}}  + \frac{\rho_2}{\rho}  \frac{1}{Re_{v2}} .
\eea\een
\rem{We define the inner product of two functions $f$ and $g$ as follows:
\ben
(f,g)=\int_V fg d{\bf x}.
\een
}
{\thm{System (\ref{Model2}) is dissipative, and the corresponding energy dissipation rate is given by }}
\ben\bea{l}
\frac{\partial E}{\partial t}  = - 2 (\frac{1}{Re_s} , {\bf D}: {\bf D} ) - (\frac{1}{Re_v}  \nabla \cdot \frac{{\bf u}}{\sqrt{{\rho}}}, \nabla \cdot \frac{{\bf u}}{\sqrt{{\rho}}})  -  ( \nabla \mu_1, \nabla \mu_2) \cdot \mathcal{M} \cdot ( \nabla \mu_1, \nabla \mu_2)^T\leq  0,
\eea\een
where $Re_s, Re_v \geq 0, \mathcal{M} =\left(
\bea{cc}
M_1  & - M_1  \\
- M_1  & M_1
\eea
\right) \geq 0$.

{\bf Proof:} By the definition of E, we have
\ben\bea{l}
\frac{\partial E}{\partial t} =  \int_V \big [{\bf u}^T \frac{\partial {\bf u}}{\partial t} + 2q_1 \frac{\partial q_1}{\partial t} + (\nabla \rho_1, \nabla \rho_2) \cdot {\bf K} \cdot (\nabla \frac{\partial \rho_1}{\partial t},  \nabla \frac{\partial \rho_2}{\partial t})^T \big ] d{\bf x}.
\eea\een
Taking the inner product of (\ref{Model2}-3) with ${\bf u}$ and using integration by parts, we obtain
\ben\bea{l}
({\bf u}, \frac{\partial {\bf u}}{\partial t} )
= - 2 (\frac{1}{Re_s} , {\bf D} : {\bf D}) - (\frac{1}{Re_v}   \nabla \cdot \frac{{\bf u}}{{\sqrt{\rho}}}, \nabla \cdot \frac{{\bf u}}{{\sqrt{\rho}}}) - ({\bf u},    \rho_1\frac{1}{{\sqrt{\rho}}} \nabla \mu_1 + \rho_2\frac{1}{{\sqrt{\rho}}} \nabla \mu_2).
\eea \label{u_relation-continuous}
\een
Taking the inner product of (\ref{Model2}-4) with $2 q_1$, using the identities of $\mu_i$, i= 1,2, and performing integration by parts, we obtain
\ben\bea{l}
(2q_1, \frac{\partial q_1}{\partial t})

=  -   ( \nabla \mu_1, \nabla \mu_2) \cdot \mathcal{M} \cdot ( \nabla \mu_1, \nabla \mu_2)^T + (\frac{ {\rho}_1}{{\sqrt{\rho}}} {\bf u} , \nabla {\mu_1}) + (\frac{ {\rho}_2}{{\sqrt{\rho}}} {\bf u} , \nabla  {\mu}_2)\\

-  (\nabla \rho_1, \nabla \rho_2) \cdot {\bf K} \cdot (\nabla \frac{\partial \rho_1}{\partial t},  \nabla \frac{\partial \rho_2}{\partial t})^T.
\eea \label{q_relation-continuous}
\een
Combining (\ref{u_relation-continuous}) and (\ref{q_relation-continuous}), we obtain
\ben\bea{l}
\frac{\partial E}{\partial t}   = - 2 (\frac{1}{Re_s} , {\bf D} : {\bf D}) - (\frac{1}{Re_v}   \nabla \cdot \frac{{\bf u}}{{\sqrt{\rho}}}, \nabla \cdot \frac{{\bf u}}{{\sqrt{\rho}}})  -  ( \nabla \mu_1, \nabla \mu_2) \cdot \mathcal{M} \cdot ( \nabla \mu_1, \nabla \mu_2)^T \leq 0
\eea\een
provided $ \mathcal{M}  \geq 0$.

We next design a second order energy stable numerical scheme based on the reformulated governing system of equations.


\section{Linear,  Second Order Energy Stable  Numerical Scheme}
\subsection{Notations and Useful Lemmas}
\begin{figure}
\centering
{\includegraphics[width=0.5\textwidth]{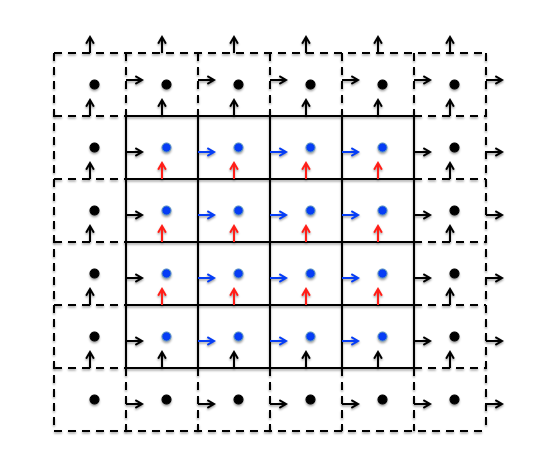}}
\caption{Staggered grid in 2D space.}
\label{fig:stagger}
\end{figure}
\noindent \indent We first introduce some notations,  finite difference operators and useful lemmas. Here, we follow the notations  in \cite{Chen&ShenJCP2016, SWWW2011, Wise&Kim&LowengrubJCP2007}.   Let $\Omega = [0, L_x] \times [0, L_y]$ be the computational domain with $L_x = h_x  \times N_x$, $L_y = h_y \times N_y$, where $N_x, N_y$ are positive integers, and $h_x, h_y$ are spatial step sizes in the x and y direction, respectively. We define three sets for the grid points as follows
 \ben\bea{l}
 E_x := \{x_{i+1/2}  = i \cdot h \quad | \quad i = 0, 1, \cdots, N_x \}, \\
 C_x := \{x_{i}  = (i-\frac{1}{2}) \cdot h \quad | \quad i = 1, \cdots, N_x \}, \\
 C_{\overline{x}} := \{x_{i}  = (i - \frac{1}{2}) \cdot h \quad  | \quad i = 0, 1, \cdots, N_x, N_x + 1 \},
 \eea\een
 where $E_x$ is a uniform partition of $[0, L_x]$ of size $N_x$ in the x-direction and its elements are called edge-centered points. The elements of $C_x$ and $C_{\overline{x}}$ are called cell-centered points. The two points belonging to $C_{\overline{x}} \string\ C_{x}$ are called ghost points. Analogously, we define $E_y$ as the uniform partition of $[0, L_y]$ of size $N_y$, called edge-centered points in the y-direction, and $C_y, C_{\overline{y}}$ the cell-centered points of the interval $[0, L_y]$. In Figure \ref{fig:stagger}, we show the staggered grid in 2D space. In this paper, we discretize the variables with the Neumann boundary conditions at the cell-center and the ones with the  Dirichlet boundary conditions at the edge-center. We define the corresponding discrete function space on this staggered grid as follows
 \ben\bea{l}
 \mathcal{C}_{x \times y} := \{ \phi : C_x \times C_y \rightarrow \mathcal{R}   \} , \quad  \mathcal{C}_{\overline{x} \times y} := \{ \phi : C_{\overline{x}} \times C_y \rightarrow \mathcal{R}   \} , \quad  \mathcal{C}_{x \times \overline{y}} := \{ \phi : C_x \times C_{\overline{y}} \rightarrow \mathcal{R}   \} , \\
   \mathcal{C}_{\overline{x} \times \overline{y}} := \{ \phi : C_{\overline{x}} \times C_{\overline{y}} \rightarrow \mathcal{R}   \} , \quad \mathcal{E}_{x \times y}^{ew} := \{ \phi : E_x \times C_y \rightarrow \mathcal{R}   \} , \quad \mathcal{E}_{x \times \overline{y}}^{ew} := \{ \phi : E_x \times C_{\overline{y}} \rightarrow \mathcal{R}   \} , \\
 \mathcal{E}_{x \times y}^{ns} := \{ \phi : C_x \times E_y \rightarrow \mathcal{R}   \} , \quad \mathcal{E}_{\overline{x} \times y}^{ns} := \{ \phi : C_{\overline{x}} \times E_y \rightarrow \mathcal{R}   \} , \quad \mathcal{V}_{x \times y} := \{ \phi : E_x \times E_y \rightarrow \mathcal{R}   \} .
 \eea\een
 $\mathcal{C}_{x \times y}, \mathcal{C}_{\overline{x} \times y}, \mathcal{C}_{x \times \overline{y}}$ and $\mathcal{C}_{\overline{x} \times \overline{y}}$ are the sets for discrete cell-centered functions, and $\mathcal{E}_{x \times y}^{ew}, \mathcal{E}_{x \times y}^{ns}$ east-west and north-south edge-centered functions, respectively.

\subsubsection{Average and Difference Operators}

\noindent \indent Assume $u, r \in \mathcal{E}_{x \times y}^{ew} \cup \mathcal{E}_{x \times \overline{y}}^{ew} $, $v, w \in \mathcal{E}_{x \times y}^{ns} \cup  \mathcal{E}_{\overline{x} \times y}^{ns} $, $\phi, \psi \in \mathcal{C}_{x \times y} \cup \mathcal{C}_{\overline{x} \times y} \cup \mathcal{C}_{x \times \overline{y}} \cup \mathcal{C}_{\overline{x} \times \overline{y}}$ and $f \in \mathcal{V}_{x \times y}$, we define the east-west-edge-to-center average and difference operator as $a_x, d_x : \mathcal{E}_{x \times \overline{y}}^{ew} \cup \mathcal{V}_{x \times y} \rightarrow \mathcal{C}_{x \times \overline{y}} \cup \mathcal{E}_{x \times y}^{ns}$ component-wise as follows
\ben\bea{l}
a_x u_{i,j} := \frac{1}{2} (u_{i+\frac{1}{2},j} + u_{i-\frac{1}{2},j} ), \quad d_x u_{i,j} := \frac{1}{h_x} (u_{i+\frac{1}{2},j} - u_{i-\frac{1}{2},j} ), \\
a_x f_{i,j +\frac{1}{2} } := \frac{1}{2} (f_{i+\frac{1}{2},j+\frac{1}{2} } + f_{i-\frac{1}{2},j+\frac{1}{2}} ), \quad d_x f_{i,j+\frac{1}{2}} := \frac{1}{h_x} (f_{i+\frac{1}{2},j+\frac{1}{2}} - f_{i-\frac{1}{2},j+\frac{1}{2}} ).
\eea\een
The north-south-edge-to-center average and difference operators are defined as $a_y, d_y : \mathcal{E}_{\overline{x} \times {y}}^{ns} \cup \mathcal{V}_{x \times y} \rightarrow \mathcal{C}_{\overline{x} \times {y}} \cup \mathcal{E}_{x \times y}^{ew}$ component-wise as follows
\ben\bea{l}
a_y v_{i,j} := \frac{1}{2} (v_{i,j+\frac{1}{2}} + v_{i,j-\frac{1}{2}} ), \quad d_y v_{i,j} := \frac{1}{h_y} (v_{i,j+\frac{1}{2}} - v_{i,j-\frac{1}{2}} ), \\
a_y f_{i+\frac{1}{2}, j} := \frac{1}{2} (f_{i+\frac{1}{2},j+\frac{1}{2} } + f_{i+\frac{1}{2},j-\frac{1}{2}} ), \quad d_y f_{i+\frac{1}{2},j} := \frac{1}{h_y} (f_{i+\frac{1}{2},j+\frac{1}{2}} - f_{i+\frac{1}{2},j-\frac{1}{2}} ).
\eea\een
We denote the center-to-east-west-edge average and difference operators as $A_x, D_x : \mathcal{C}_{\overline{x} \times \overline{y}} \cup \mathcal{E}_{\overline{x} \times y}^{ns} \rightarrow \mathcal{E}_{x \times \overline{y}}^{ew} \cup \mathcal{V}_{x \times y}$ in component-wise forms:
 \ben\bea{l}
A_x \phi_{i + \frac{1}{2},j} := \frac{1}{2} (\phi_{i+1,j} + \phi_{i,j} ), \quad D_x \phi_{i + \frac{1}{2},j} := \frac{1}{h_x} (\phi_{i+1,j} - \phi_{i,j} ), \\
A_x v_{i +\frac{1}{2} ,j +\frac{1}{2} } := \frac{1}{2} (v_{i+1,j+\frac{1}{2} } + v_{i,j+\frac{1}{2}} ), \quad D_x v_{i +\frac{1}{2} ,j+\frac{1}{2}} := \frac{1}{h_x} (v_{i+1,j+\frac{1}{2}} - v_{i,j+\frac{1}{2}} ).
\eea\een
Analogously, the center-to-north-south-edge average and difference operator are defined as $A_y, D_y : \mathcal{C}_{\overline{x} \times \overline{y}} \cup \mathcal{E}_{{x} \times \overline{y}}^{ew} \rightarrow \mathcal{E}_{\overline{x} \times {y}}^{ns} \cup \mathcal{V}_{x \times y}$ in component-wise forms:
 \ben\bea{l}
A_y \phi_{i,j + \frac{1}{2}} := \frac{1}{2} (\phi_{i,j + 1} + \phi_{i,j} ), \quad D_y \phi_{i ,j+ \frac{1}{2}} := \frac{1}{h_y} (\phi_{i, j+1} - \phi_{i,j} ), \\
A_y u_{i +\frac{1}{2} ,j +\frac{1}{2} } := \frac{1}{2} (u_{i+\frac{1}{2},j+1 } + u_{i+\frac{1}{2}, j} ), \quad D_y u_{i +\frac{1}{2} ,j+\frac{1}{2}} := \frac{1}{h_y} (u_{i+\frac{1}{2}, j + 1} - u_{i+\frac{1}{2},j} ).
\eea\een
The standard 2D discrete Laplacian operator is defined as $\Delta_h : \mathcal{E}_{x \times \overline{y}}^{ew} \cup \mathcal{E}_{\overline{x} \times y}^{ns} \cup \mathcal{C}_{\overline{x} \times \overline{y}} \rightarrow  \mathcal{E}_{x \times {y}}^{ew} \cup \mathcal{E}_{{x} \times y}^{ns} \cup \mathcal{C}_{{x} \times {y}}$:
\ben\bea{l}
\Delta_h u := D_x (d_x u) + d_y(D_y u), \quad \Delta_h v := d_x(D_x v) + D_y (d_yv), \quad \Delta_h \phi := d_x(D_x \phi) + d_y (D_y \phi).
\eea\een

\subsubsection{Boundary Conditions}
\noindent \indent The  homogenous Neumann boundary conditions are discretized as follows
\ben\bea{l}
\phi_{0,j} = \phi_{1,j}, \quad \phi_{N_x, j} = \phi_{N_x +1, j}, \quad j = 0, 1, 2, \cdots, N_y + 1,\\
\phi_{i,0} = \phi_{i,1}, \quad \phi_{i,N_y} = \phi_{i, N_y + 1}, \quad i = 0, 1, 2, \cdots, N_x +1.
\eea \label{eq:neumann_discrete}
\een
We denote it as ${\bf n} \cdot \nabla_h \phi  |_{\partial V}= 0$.

The homogeneously Dirichlet boundary conditions are discretized as follows
\ben\bea{l}
u_{\frac{1}{2}, j} = u_{N_x + \frac{1}{2}, j} = 0, \quad j = 1, 2, \cdots, N_y,\\
A_y u_{i + \frac{1}{2}, \frac{1}{2}} = A_y u_{i + \frac{1}{2}, N_y + \frac{1}{2}} = 0, \quad i = 0, 1, 2, \cdots, N_x,\\
v_{i, \frac{1}{2}} = v_{i, N_y + \frac{1}{2}} = 0, \quad j = 1, 2, \cdots, N_x,\\
A_x v_{\frac{1}{2}, j+\frac{1}{2}} = A_x v_{N_x + \frac{1}{2}, j+\frac{1}{2}} = 0, \quad j =0, 1, 2, \cdots, N_y,\\
\eea \label{eq:dirichlet_discrete}
\een
We denote it as $u_h |_{\partial V} = 0$ and $v_h |_{\partial V} = 0$.

If $f \in \mathcal{V}_{x \times y}$ satisfies homogenous Dirichelet boundary condition, we have
\ben\bea{l}
f_{\frac{1}{2}, j + \frac{1}{2} } = f_{N_x + \frac{1}{2}, j + \frac{1}{2} }  = f_{i + \frac{1}{2},  \frac{1}{2} } = f_{i + \frac{1}{2}, N_y + \frac{1}{2} }  = 0.
\eea\een
where $i = 0, 1, 2, \cdots, N_x,  j =0, 1, 2, \cdots, N_y$. We denote is as $f_h |_{\partial V} = 0$.

\subsubsection{Inner products and norms}

\noindent \indent We defined the following inner products for discrete functions
\ben\bea{l}
(\phi, \psi)_2 := h_x h_y \sum_{i=1}^{N_x} \sum_{j=1}^{N_y} \phi_{i,j} \psi_{i,j},\quad [u, r]_{ew} := (a_x (u r), 1)_2, \quad
[v, w]_{ns} := (a_y (uw), 1)_2, \\ (f, g)_{vc} := (a_x (a_y (fg)), 1)_2,\quad (\nabla \phi, \nabla \psi)_h := [D_x(\phi), D_x(\psi)]_{ew} + [D_y(\phi), D_y(\psi)]_{ns},
\eea\een
The corresponding norms are defined as follows
\ben\bea{l}
\| \phi  \|_2 := (\phi, \phi)_2^{\frac{1}{2}}, \quad \| u \|_{ew} := [u, u]_{ew}^{\frac{1}{2}}, \quad \| v\|_{ns} := [v, v]_{ns}^{\frac{1}{2}}, \quad \| f\|_{vc} := (f, f)_{vc}^{\frac{1}{2}}.
\eea\een
For $\phi = \mathcal{C}_{\overline{x} \times \overline{y}}$, we define $\| \nabla \phi \|_2$ as
\ben\bea{l}
\| \nabla \phi \|_2  := \sqrt{\| D_x \phi  \|^2_{ew}  +  \|  D_y \phi  \|^2_{ns}},
\eea\een
For the edge-centered velocity vector ${\bf v} = (u, v)$, $u \in \mathcal{E}_{x \times \overline{y}}^{ew}, v \in \mathcal{E}_{\overline{x} \times y}^{ns}$, we define $\| {\bf v} \|_2, \|  \nabla {\bf v} \|_2$ as
\ben\bea{l}
\| {\bf v} \|_2  :=  \sqrt{\| u  \|^2_{ew}  +  \|  v \|^2_{ns}}, \quad  \| \nabla {\bf v} \|_2  :=  \sqrt{\| d_x u  \|^2_2  +  \|  D_y u \|^2_{vc}  + \|  D_x v \|^2_{vc}  + \| d_y v  \|^2_2},\\
\| {\bf D} \|_2   := \sqrt{\| d_x u  \|^2_2  +  \frac{1}{2} \|  D_y u \|^2_{vc}  + \frac{1}{2}\|  D_x v \|^2_{vc}  + (D_y u, D_x v)_{vc} + \| d_y v  \|^2_2},\\
(\phi, {\bf D}: {\bf D})_2 :=\\
 \sqrt{( \phi,  (d_x u)^2 )_2  +  \frac{1}{2} (  A_x (A_y\phi) , (D_y u)^2 )_{vc} + \frac{1}{2}( A_x (A_y\phi),  (D_x v)^2)_{vc}  + (A_x (A_y\phi) D_y u, D_x v)_{vc} + ( \phi, (d_y v) ^2)_2}.
\eea\een
Where ${\bf D} = \frac{1}{2} (\nabla {\bf v} + \nabla {\bf v}^T)$.
From these definitions, we obtain the following lemmas \cite{Gong&Zhao&Wang2018_SIAM_JSC2}:
\begin{lem} \label{lem1}
(Summation by parts): If $\phi \in \mathcal{C}_{\overline{x} \times \overline{y}}$, $u \in \mathcal{E}_{{x} \times {y}}^{ew}$,  $v \in \mathcal{E}_{{x} \times {y}}^{ns}$, and $u_h |_{\partial V} = 0$ and $v_h |_{\partial V} = 0$, then
\ben\bea{l}
\quad [A_x \phi, u]_{ew} = (\phi, a_x u)_2, \quad [A_y \phi, v]_{ns} = (\phi, a_y v)_2,\\
\quad [D_x \phi, u]_{ew} + (\phi, d_x u)_2 = 0,  \quad [D_y \phi, v]_{ns} + (\phi, d_y v)_2 = 0,
\eea
\een
\end{lem}
\begin{lem}\label{lem2}
If $f \in \mathcal{V}_{{x} \times {y}}$, and $f_h |_{\partial V} = 0$, $u \in \mathcal{E}_{{x} \times \overline{y}}^{ew}$,  $v \in \mathcal{E}_{\overline{x} \times {y}}^{ns}$, then
\ben\bea{l}
[a_y f, u]_{ew} = (f, A_y u)_{vc}, \quad [a_x f, v]_{ns} = (f, A_x v)_{vc}.
\eea
\een
\end{lem}
\begin{lem} \label{lem3}
If $f \in \mathcal{V}_{{x} \times {y}}$, $u \in \mathcal{E}_{{x} \times \overline{y}}^{ew}$,  $v \in \mathcal{E}_{\overline{x} \times {y}}^{ns}$, and $u_h |_{\partial V} = 0$ and $v_h |_{\partial V} = 0$, then
\ben\bea{l}
[d_yf, u]_{ew} + (f, D_y u)_{vc} = 0, \quad [d_x f, v]_{ns} + (f, D_x v)_{vc} = 0.
\eea
\een
\end{lem}
\noindent With these notations and lemmas, we are ready to introduce the fully-discrete numerical scheme in the following section.


\subsection{Semi-discrete scheme in time}

\noindent \indent First, we discretize the governing equations using Crank-Nicolson method in time.
We denote
\ben\bea{l}
\delta_t (\cdot)^{n+1/2} =\frac{1}{\Delta t}( (\cdot)^{n+1} - (\cdot)^n), \qquad \overline{(\cdot)}^{n+1/2} = \frac{1}{2} (3(\cdot)^{n} - (\cdot)^{n-1}).
\eea\een
The second order algorithm is given below.

\begin{alg}
\ben\bea{l}
\begin{cases}
\delta_t \rho_1^{n+1/2} +  \nabla \cdot ( \overline{\rho}_1^{n + 1/2}    \overline{\frac{1}{\sqrt{\rho}}}^{n+1/2}      {\bf u}^{n + 1/2}) =   \nabla \cdot {M}_1\cdot \nabla {\mu}_1^{n+1/2} -  \nabla \cdot {M}_1\cdot \nabla {\mu}_2^{n+1/2} ,\\

\delta_t \rho_2^{n+1/2} +  \nabla \cdot ( \overline{\rho}_2^{n + 1/2}    \overline{\frac{1}{\sqrt{\rho}}}^{n+1/2}      {\bf u}^{n + 1/2}) =  - \nabla \cdot {M}_1\cdot \nabla {\mu}_1^{n+1/2} +  \nabla \cdot {M}_1\cdot \nabla {\mu}_2^{n+1/2} ,\\

\delta_t {\bf u}^{n+1/2}   + \frac{1}{2}( \overline{\frac{1}{\sqrt{\rho}}}^{n+1/2}     \nabla \cdot ( \overline{\bf u}^{n+1/2}{\bf u}^{n+1/2}) + \overline{\bf u}^{n+1/2} \cdot \nabla (\overline{\frac{1}{\sqrt{\rho}}}^{n+1/2}  {\bf u}^{n+1/2}))=\\

 \overline{\frac{1}{\sqrt{\rho}}}^{n+1/2}
(2 \nabla \cdot ( \frac{1}{Re_s}  {\bf D}^{n+1/2}) + \nabla (  \frac{1}{Re_v}  \nabla \cdot (\overline{\frac{1}{\sqrt{\rho}}}^{n+1/2}  {\bf u}^{n+1/2}) )

-  \overline{\rho_1}^{n+1/2} \nabla {\mu}_1^{n+1/2} - \overline{ \rho_2}^{n+1/2} \nabla {\mu_2}^{n+1/2}),\\

\delta_t q_1^{n+1/2} = \overline{\frac{\partial q_1}{\partial \rho_1}}^{n+1/2} \delta_t \rho_1^{n+1/2} + \overline{\frac{\partial q_1}{\partial \rho_2}}^{n+1/2} \delta_t \rho_2^{n+1/2},
\end{cases}
\eea
\label{eq:CN1}
\een
where
\ben\bea{l}
 \mu_1^{n+1/2} = 2q_1^{n+1/2} \overline{\frac{\partial q_1}{\partial \rho_1}}^{n+1/2} - \kappa_{\rho_1 \rho_1} \Delta \rho_1^{n+1/2} - \kappa_{\rho_1 \rho_2} \Delta \rho_2^{n+1/2},\\

 \mu_2^{n+1/2} = 2q_1^{n+1/2} \overline{\frac{\partial q_1}{\partial \rho_2}}^{n+1/2} - \kappa_{\rho_1 \rho_2} \Delta \rho_1^{n+1/2} - \kappa_{\rho_2 \rho_2} \Delta \rho_2^{n+1/2},\\

{\bf D}^{n+1/2} = \frac{1}{2} (\nabla (\overline{\frac{1}{\sqrt{\rho}}}^{n+1/2} {\bf u}^{n+1/2}  )+  \nabla (\overline{\frac{1}{\sqrt{\rho}}}^{n+1/2} {\bf u}^{n+1/2} )^T),\\

\frac{1}{Re_s}  =  \overline{\frac{\rho_1}{\rho}}^{n+1/2}  \frac{1}{Re_{s1}}  +  \overline{\frac{ \rho_2}{\rho}}^{n+1/2}  \frac{1}{Re_{s2}} , \qquad

\frac{1}{Re_v}  =  \overline{\frac{\rho_1}{\rho}}^{n+1/2}  \frac{1}{Re_{v1}}  +  \overline{\frac{ \rho_2}{\rho}}^{n+1/2}  \frac{1}{Re_{v2}}.
\eea\een
\end{alg}

For the scheme, we have the following theorem.
{\thm{Scheme \eqref{eq:CN1} is unconditional energy stable, and satisfies the following discrete energy identity }
\ben\bea{l}
\frac{E^{n+1} - E^n}{\Delta t}  = - 2 (\frac{1}{Re_s} , {\bf D}^{n+1/2} : {\bf D}^{n+1/2}) - (\frac{1}{Re_v}  \nabla \cdot ( \overline{\frac{1}{\sqrt{\rho}}}^{n+1/2} {\bf u}^{n+1/2})  , \nabla \cdot ( \overline{\frac{1}{\sqrt{\rho}}}^{n+1/2}{\bf u}^{n+1/2} )  ) \\

- ( \nabla \mu_1^{n+1/2}, \nabla \mu_2^{n+1/2}) \cdot \mathcal{M}  \cdot ( \nabla \mu_1^{n+1/2}, \nabla \mu_2^{n+1/2})^T< 0,
\eea\een
Where
\ben\bea{l}
E^n = \int_V [ \frac{1}{2}|| {\bf u}^n||^2 + (q_1^n)^2 + \frac{1}{2} ({\bf p}^n)^T \cdot {\bf K}^n \cdot {\bf p}^n - A ] d{\bf x}.
\eea\een}
and ${\bf p}^n = (\nabla \rho_1^n, \nabla \rho_2^n)$.
\rem{
We note that a useful identity in the proof of the theorem.
\ben\bea{l}
({\bf u}^{n+1/2} ,  \frac{1}{2}( \overline{\frac{1}{\sqrt{\rho}}}^{n+1/2}   \nabla \cdot ( \overline{\bf u}^{n+1/2}{\bf u}^{n+1/2}) + \overline{\bf u}^{n+1/2} \cdot \nabla (\overline{\frac{1}{\sqrt{\rho}}}^{n+1/2}  {\bf u}^{n+1/2})  )  ) = 0.
\eea \label{Identity1}
\een
}
{\bf Proof:} By the definition of $E^n$, we have
\ben\bea{l}
\frac{E^{n+1} - E^n}{\Delta t} =  \int_V {\bf u}^{n+1/2} \delta_t {\bf u}^{n+1/2} +  2q_1^{n+1/2} \delta_t q_1^{n+1/2} +  \kappa_{\rho_1 \rho_1} \nabla \rho^{n+1/2} \delta_t \nabla \rho_1^{n+1/2} \\

+ \kappa_{\rho_2 \rho_2} \nabla \rho_2^{n+1/2} \delta_t \nabla \rho_2^{n+1/2} + \kappa_{\rho_1 \rho_2} [ \nabla \rho_1^{n+1/2} \delta_t \nabla  \rho_2^{n+1/2}  +  \nabla \rho_2^{n+1/2} \delta_t \nabla \rho_1^{n+1/2} ]

d{\bf x},
\eea \label{eq:energy_differnce}
\een
Taking the inner product of \eqref{eq:CN1}-3 with ${\bf u}^{n+1/2}$, using identity (\ref{Identity1}), and performing integration by parts, we obtain
\ben\bea{l}
({\bf u}^{n+1/2}, \delta_t {\bf u}^{n+1/2})
= - 2 (\frac{1}{Re_s} {\bf D}^{n+1/2} : {\bf D}^{n+1/2}) - ( \frac{1}{Re_v}   \nabla \cdot ( \overline{\frac{1}{\sqrt{\rho}}}^{n+1/2} {\bf u}^{n+1/2})  , \nabla \cdot ( \overline{\frac{1}{\sqrt{\rho}}}^{n+1/2}{\bf u}^{n+1/2} )  ) \\

-  ({\bf u}^{n+1/2} ,  \overline{\frac{1}{\sqrt{\rho}}}^{n+1/2}  \overline{\rho_1}^{n+1/2} \nabla {\mu}_1^{n+1/2} + \overline{\frac{1}{\sqrt{\rho}}}^{n+1/2}  \overline{ \rho_2}^{n+1/2} \nabla {\mu_2}^{n+1/2}) ).
\eea \label{u_relation}
\een
Taking the inner product of (\ref{eq:CN1}-4) with $2q_1^{n + \frac{1}{2}}$, using (\ref{eq:CN1}-1,2), and performing integration by parts, we obtain
\ben\bea{l}
2 (q_1^{n+1/2} , \delta_t q_1^{n+1/2})

= (  \overline{\frac{1}{\sqrt{\rho}}}^{n+1/2} \overline{\rho}_1^{n+1/2}  , \nabla {\mu_1}^{n+1/2})

+ (\overline{{\rho_2}}^{n+1/2} \overline{\frac{1}{\sqrt{\rho}}}^{n+1/2} , \nabla  {\mu_2}^{n+1/2})  \\

-  \kappa_{\rho_1 \rho_1} \nabla \rho_1^{n+1/2} \delta_t  \nabla \rho_1^{n+1/2}

- \kappa_{\rho_2 \rho_2} \nabla \rho_2^{n+1/2} \delta_t \nabla  \rho_2^{n+1/2} - \kappa_{\rho_1 \rho_2} [ \nabla \rho_1^{n+1/2} \delta_t \nabla \rho_2^{n+1/2}  \\
+  \nabla \rho_2^{n+1/2} \delta_t \nabla \rho_1^{n+1/2} ]

-   ( \nabla \mu_1^{n+1/2}, \nabla \mu_2^{n+1/2}) \cdot \mathcal{M}  \cdot ( \nabla \mu_1^{n+1/2}, \nabla \mu_2^{n+1/2})^T.
\eea \label{q_relation}
\een
Utilizing (\ref{eq:energy_differnce}), (\ref{u_relation}) and (\ref{q_relation}), we arrive at the conclusion
\ben\bea{l}
\frac{E^{n+1} - E^n}{\Delta t}  = - 2 (\frac{1}{Re_s} {\bf D}^{n+1/2} : {\bf D}^{n+1/2}) - ( \frac{1}{Re_v}   \nabla \cdot ( \overline{\frac{1}{\sqrt{\rho}}}^{n+1/2} {\bf u}^{n+1/2})  , \nabla \cdot ( \overline{\frac{1}{\sqrt{\rho}}}^{n+1/2}{\bf u}^{n+1/2} )  ) \\

-  ( \nabla \mu_1^{n+1/2}, \nabla \mu_2^{n+1/2}) \cdot  \mathcal{M}  \cdot ( \nabla \mu_1^{n+1/2}, \nabla \mu_2^{n+1/2})^T \leq 0
\eea\een
provided $ \mathcal{M} \geq 0$.

\subsection{Fully Discrete Numerical Scheme}

\noindent \indent
We discretize the semidiscrete equations in \eqref{eq:CN1} using the second order finite difference discretization  on staggered grids in space to obtain a fully discrete scheme as follows
\begin{alg}
\ben\bea{l}
\begin{cases}
\big \{ \delta_t \rho_1^{n+1/2}  + d_x (A_x(   \overline{\rho}_1^{n + 1/2}    \overline{\frac{1}{\sqrt{\rho}}}^{n+1/2}    ) u^{n+1/2}) +  d_y (A_y(\overline{\rho}_1^{n + 1/2}    \overline{\frac{1}{\sqrt{\rho}}}^{n+1/2}    ) v^{n+1/2}) =   \\
M_1 \Delta_h \mu_1^{n+1/2} -   M_1 \Delta_h  \mu_2^{n+1/2} \big \}|_{i, j}, i = 1, \cdots, N_x, j = 1, \cdots, N_y, \\
\\
\big \{ \delta_t \rho_2^{n+1/2}  + d_x (A_x(\overline{\rho}_2^{n + 1/2}    \overline{\frac{1}{\sqrt{\rho}}}^{n+1/2}    ) u^{n+1/2}) +  d_y (A_y(\overline{\rho}_2^{n + 1/2}    \overline{\frac{1}{\sqrt{\rho}}}^{n+1/2}    ) v^{n+1/2})=  \\
 - M_1 \Delta_h \mu_1^{n+1/2} +  M_1 \Delta_h \mu_2^{n+1/2} \big \}|_{i, j},
 i = 1, \cdots, N_x, j = 1, \cdots, N_y, \\
\\
\big \{ \delta_t u^{n+1/2}    + \frac{1}{2} (\overline{u}^{n+1/2} D_x( \overline{\frac{1}{\sqrt{\rho}}}^{n+1/2} a_x u^{n+1/2}) + A_x( \overline{\frac{1}{\sqrt{\rho}}}^{n+1/2} d_x (\overline{ u}^{n+1/2} u^{n+1/2}))) \\

+ \frac{1}{2}(a_x (A_x \overline{v}^{n+1/2} D_y(A_x( \overline{\frac{1}{\sqrt{\rho}}}^{n+1/2} )u^{n+1/2})) + A_x( \overline{\frac{1}{\sqrt{\rho}}}^{n+1/2} )d_y(A_yu^{n+1/2} A_x(\overline{v}^{n+1/2})) \\
= g_{v1} \big \}|_{i+\frac{1}{2}, j},
 i = 1, \cdots, N_x - 1, j = 1, \cdots, N_y ,\\
\\

\big \{ \delta_t v^{n+1/2}   + \frac{1}{2} (a_x(A_y \overline{u}^{n+1/2} D_x(A_y( \overline{\frac{1}{\sqrt{\rho}}}^{n+1/2} )v^{n+1/2})) + A_y( \overline{\frac{1}{\sqrt{\rho}}}^{n+1/2} )d_x(A_y\overline{u}^{n+1/2} A_xv^{n+1/2})) \\
+ \frac{1}{2} (\overline{v}^{n+1/2} D_y( \overline{\frac{1}{\sqrt{\rho}}}^{n+1/2} a_yv^{n+1/2}) + A_y( \overline{\frac{1}{\sqrt{\rho}}}^{n+1/2} d_y(\overline{v}^{n+1/2}v^{n+1/2}))) \\
=  g_{v2} \big \}|_{i, j+\frac{1}{2}}, i = 1, \cdots, N_x, j = 1, \cdots, N_y - 1,\\
\\

\big \{\delta_t q_1^{n+1/2} = \overline{{\frac{\partial q_1}{\partial \rho_1}}}^{n+1/2}  \delta_t \rho_1^{n+1/2} + \overline{\frac{\partial q_1}{\partial \rho_2}}^{n+1/2}  \delta_t \rho_2^{n+1/2} \big \} |_{i,j}, i = 1, \cdots, N_x , j=1, \cdots, N_y,
\end{cases}
\eea \label{eq:full_discrete}
\een
where
\ben\bea{l}
g_{v1} = A_x( \overline{\frac{1}{\sqrt{\rho}}}^{n+1/2} )  (2 D_x( \frac{1}{Re_s^{n+1/2}} d_x  (A_x( \overline{\frac{1}{\sqrt{\rho}}}^{n+1/2} ) u^{n+1/2})) + d_y( A_x(A_y \frac{1}{Re_s^{n+1/2}})  D_y  (A_x( \overline{\frac{1}{\sqrt{\rho}}}^{n+1/2} ) u^{n+1/2}))) \\

+ A_x( \overline{\frac{1}{\sqrt{\rho}}}^{n+1/2} )  d_y( A_x(A_y \frac{1}{Re_s^{n+1/2}})  D_x( A_y( \overline{\frac{1}{\sqrt{\rho}}}^{n+1/2} ) v^{n+1/2})) \\
 + A_x( \overline{\frac{1}{\sqrt{\rho}}}^{n+1/2} ) D_x( \frac{1}{Re_v^{n+1/2}} d_x( A_x( \overline{\frac{1}{\sqrt{\rho}}}^{n+1/2} ) u^{n+1/2}))
 + A_x( \overline{\frac{1}{\sqrt{\rho}}}^{n+1/2} )  D_x( \frac{1}{Re_v^{n+1/2}} d_y( A_y( \overline{\frac{1}{\sqrt{\rho}}}^{n+1/2} ) v^{n+1/2}))  \\

- A_x(\overline{\rho_1}^{n+1/2}  \overline{\frac{1}{\sqrt{\rho}}}^{n+1/2} ) D_x(\mu_1^{n+1/2})
-  A_x( \overline{\rho_2}^{n+1/2}  \overline{\frac{1}{\sqrt{\rho}}}^{n+1/2} ) D_x(\mu_2^{n+1/2}),\\
\eea\een

\ben\bea{l}
g_{v2} = A_y( \overline{\frac{1}{\sqrt{\rho}}}^{n+1/2} )  (d_x( A_x(A_y \frac{1}{Re_s^{n+1/2}})  D_x(  A_y( \overline{\frac{1}{\sqrt{\rho}}}^{n+1/2} ) v^{n+1/2})) + 2 D_y( \frac{1}{Re_s^{n+1/2}} d_y ( A_y( \overline{\frac{1}{\sqrt{\rho}}}^{n+1/2} ) v^{n+1/2}))) \\

+  A_y( \overline{\frac{1}{\sqrt{\rho}}}^{n+1/2}  ) d_x( A_x(A_y \frac{1}{Re_s^{n+1/2}})  D_y( A_x( \overline{\frac{1}{\sqrt{\rho}}}^{n+1/2} ) u^{n+1/2}))  \\

+  A_y( \overline{\frac{1}{\sqrt{\rho}}}^{n+1/2}  )  D_y( \frac{1}{Re_v^{n+1/2}} d_x( A_x(\frac{1}{\sqrt{\rho}}) u^{n+1/2}))
+  A_y( \overline{\frac{1}{\sqrt{\rho}}}^{n+1/2}   ) D_y( \frac{1}{Re_v^{n+1/2}} d_y( A_y( \overline{\frac{1}{\sqrt{\rho}}}^{n+1/2} ) v^{n+1/2})) \\

 - A_y( \overline{\rho_1}^{n+1/2} \overline{\frac{1}{\sqrt{\rho}}}^{n+1/2}  ) D_y(\mu_1^{n+1/2}) -  A_y( \overline{\rho_2}^{n+1/2}  \overline{\frac{1}{\sqrt{\rho}}}^{n+1/2} ) D_y(\mu_2^{n+1/2}).
\eea\een

For any time step $t_n$, $\rho_i^n, \mu_i^n, i=1, 2$ and $q_1^n$ satisfy discrete homogeneous Neumann boundary conditions \eqref{eq:neumann_discrete}, $u^n, v^n$ satisfy the discrete homogeneous Dirichlet boundary conditions \eqref{eq:dirichlet_discrete}. The discrete Reynolds numbers are defined as follows
\ben\bea{l}
 \{ \frac{1}{Re_s^{n+1/2}}  = \overline{(\frac{\rho_1}{\rho})}^{n+1/2}  \frac{1}{Re_{s1}}  + \overline{(\frac{\rho_2}{\rho})}^{n+1/2}   \frac{1}{Re_{s2}} \} |_{i,j}, i = 1, \cdots, N_x , j=1, \cdots, N_y,\\
 \{ \frac{1}{Re_v^{n+1/2}}  = \overline{(\frac{\rho_1}{\rho})}^{n+1/2}  \frac{1}{Re_{v1}}  + \overline{(\frac{\rho_2}{\rho})}^{n+1/2}   \frac{1}{Re_{v2}}\} |_{i,j}, i = 1, \cdots, N_x , j=1, \cdots, N_y .
\eea\een
\end{alg}

{\thm{Scheme (\ref{eq:full_discrete}) is unconditionally energy stable, and the discrete total energy satisfies the following identity }
\ben\bea{l}
\frac{E_h^{n+1} - E_h^n}{\Delta t}  = - 2 (\frac{1}{Re_s}, {\bf D}_h^{n+1/2} : {\bf D}_h^{n+1/2})_2 - (\frac{1}{Re_v} tr({\bf D}_h^{n+1/2}), tr({\bf D}_h^{n+1/2}) )_2 \\

- M_1 ( \nabla (\mu_1^{n+1/2} -  \mu_2^{n+1/2}),   \nabla (\mu_1^{n+1/2} -  \mu_2^{n+1/2}))_2 \leq 0,
\eea\een
where
\ben\bea{l}
E_h^n =  \frac{1}{2} [u^n, u^n]_{ew} + \frac{1}{2} [v^n, v^n]_{ns} +  ( q_1^n, q_1^n)_2 - (A,1)_2  \\

 + \frac{1}{2} \kappa_{\rho_1 \rho_1} (\nabla \rho_1^n, \nabla \rho_1^n)_h  +  \frac{1}{2} \kappa_{\rho_2 \rho_2} (\nabla \rho_2^n, \nabla \rho_2^n)_h + \kappa_{\rho_1 \rho_2} (\nabla \rho_1^n, \nabla \rho_2^n)_h.
\eea\een}
and
\ben\bea{l}
{\bf D}_h^{n+1/2} =  \left(
\bea{cc}
d_x  (A_x( \overline{\frac{1}{\sqrt{\rho}}}^{n+1/2} ) u^{n+1/2}))  & \frac{1}{2}S \\
\frac{1}{2} S  & d_y ( A_y( \overline{\frac{1}{\sqrt{\rho}}}^{n+1/2} ) v^{n+1/2}))
\eea
\right)
\eea\een
where $S = D_x( A_y( \overline{\frac{1}{\sqrt{\rho}}}^{n+1/2} ) v^{n+1/2})  +  D_y  (A_x( \overline{\frac{1}{\sqrt{\rho}}}^{n+1/2} ) u^{n+1/2}))$.

\rem{
We note that using lemmas (\ref{lem1})-(\ref{lem3}), we could obtain identities as follows
\ben\bea{l}

(u^{n+1/2} ,  \frac{1}{2} (\overline{u}^{n+1/2} D_x( \overline{\frac{1}{\sqrt{\rho}}}^{n+1/2} a_x u^{n+1/2}) + A_x( \overline{\frac{1}{\sqrt{\rho}}}^{n+1/2} d_x (\overline{ u}^{n+1/2} u^{n+1/2}))) \\

+ \frac{1}{2}(a_x (A_x \overline{v}^{n+1/2} D_y(A_x( \overline{\frac{1}{\sqrt{\rho}}}^{n+1/2} )u^{n+1/2})) + A_x( \overline{\frac{1}{\sqrt{\rho}}}^{n+1/2} )d_y(A_yu^{n+1/2} A_x(\overline{v}^{n+1/2})) = 0, \\
\\

(v^{n+1/2}  , \frac{1}{2} (a_x(A_y \overline{u}^{n+1/2} D_x(A_y( \overline{\frac{1}{\sqrt{\rho}}}^{n+1/2} )v^{n+1/2})) + A_y( \overline{\frac{1}{\sqrt{\rho}}}^{n+1/2} )d_x(A_y\overline{u}^{n+1/2} A_xv^{n+1/2})) \\
+ \frac{1}{2} (\overline{v}^{n+1/2} D_y( \overline{\frac{1}{\sqrt{\rho}}}^{n+1/2} a_yv^{n+1/2}) + A_y( \overline{\frac{1}{\sqrt{\rho}}}^{n+1/2} d_y(\overline{v}^{n+1/2}v^{n+1/2})))) = 0.

\eea  \label{eq:identify_full_discrete}
\een
}

{\bf Proof:} It follows from  the definition of $E_h^n$ that
\ben\bea{l}
\frac{E_h^{n+1} - E_h^n}{\Delta t} = [\frac{u^{n+1} + u^n}{2} , \frac{u^{n+1} - u^n}{\Delta t}]_{ew} + [\frac{v^{n+1} + v^n}{2} , \frac{v^{n+1} - v^n}{\Delta t}]_{ns} + 2(\frac{q_1^{n+1} + q_1^n}{2} , \frac{q_1^{n+1} - q_1^n }{\Delta t})_2 \\

+  \kappa_{\rho_1 \rho_1}( \frac{\nabla \rho_1^{n+1} + \nabla \rho_1^n}{2} , \frac{\nabla \rho_1^{n+1} - \nabla \rho_1^n }{\Delta t})_h +  \kappa_{\rho_2 \rho_2} (\frac{\nabla \rho_2^{n+1} + \nabla \rho_2^n}{2} , \frac{\nabla \rho_2^{n+1} - \nabla \rho_2^n }{\Delta t})_h \\

+ \kappa_{\rho_1 \rho_2}  [(\frac{\nabla \rho_1^{n+1} + \nabla \rho_1^n}{2}  , \frac{\nabla \rho_2^{n+1} - \nabla \rho_2^n }{\Delta t} )_h +  (\frac{\nabla \rho_2^{n+1} + \nabla \rho_2^n}{2} , \frac{\nabla \rho_1^{n+1} - \nabla \rho_1^n }{\Delta t})_h ]
\\
\eea \label{eq:full_energy_dissipation1}
\een
Taking the inner product of (\ref{eq:full_discrete}-3,4) with $u^{n+1/2}, v^{n+1/2}$ respectively and using identify (\ref{eq:identify_full_discrete}), we obtain
\ben\bea{l}
[\frac{u^{n+1} + u^n}{2} , \frac{u^{n+1} - u^n}{\Delta t}]_{ew} + [\frac{v^{n+1} + v^n}{2} , \frac{v^{n+1} - v^n}{\Delta t}]_{ns} \\
 = - 2 (\frac{1}{Re_s}, {\bf D}_h^{n+1/2} : {\bf D}_h^{n+1/2})_2 - (\frac{1}{Re_v} tr({\bf D}_h^{n+1/2}), tr({\bf D}_h^{n+1/2}) )_2 \\

 - [u^{n+1/2}, A_x(\overline{\rho_1}^{n+1/2}  \overline{\frac{1}{\sqrt{\rho}}}^{n+1/2} ) D_x(\mu_1^{n+1/2}) +  A_x( \overline{\rho_2}^{n+1/2}  \overline{\frac{1}{\sqrt{\rho}}}^{n+1/2} ) D_x(\mu_2^{n+1/2})]_{ew}\\

 - [v^{n+1/2}, A_y( \overline{\rho_1}^{n+1/2} \overline{\frac{1}{\sqrt{\rho}}}^{n+1/2}  ) D_y(\mu_1^{n+1/2}) +  A_y( \overline{\rho_2}^{n+1/2}  \overline{\frac{1}{\sqrt{\rho}}}^{n+1/2} ) D_y(\mu_2^{n+1/2}) ]_{ns},
\eea \label{eq:full_energy_dissipation2}
\een Where we used lemmas (\ref{lem1}) and (\ref{lem3}).
Taking the inner product of (\ref{eq:full_discrete}-5) with $2 q_1^{n+1/2}$, and performing integration by parts, we obtain
\ben\bea{l}
2(\frac{q_1^{n+1} + q_1^n}{2} , \frac{q_1^{n+1} - q_1^n }{\Delta t})_2

=   -  M_1 ( \nabla (\mu_1^{n+1/2} - \mu_2^{n+1/2}), \nabla (\mu_1^{n+1/2} - \mu_2^{n+1/2}))_h\\

 + [A_x(   \overline{\rho}_1^{n + 1/2}    \overline{\frac{1}{\sqrt{\rho}}}^{n+1/2}   ) u^{n+1/2}, D_x(\mu_1^{n+1/2})]_{ew} +  [A_y(   \overline{\rho}_1^{n + 1/2}    \overline{\frac{1}{\sqrt{\rho}}}^{n+1/2}   ) v^{n+1/2}, D_y(\mu_1^{n+1/2})]_{ns} \\

 + [A_x(   \overline{\rho}_2^{n + 1/2}    \overline{\frac{1}{\sqrt{\rho}}}^{n+1/2}   ) u^{n+1/2}, D_x(\mu_2^{n+1/2})]_{ew} +  [A_y(   \overline{\rho}_2^{n + 1/2}    \overline{\frac{1}{\sqrt{\rho}}}^{n+1/2}   ) v^{n+1/2}, D_y(\mu_2^{n+1/2})]_{ns} \\

-    \kappa_{\rho_1 \rho_1}( \frac{\nabla \rho_1^{n+1} + \nabla \rho_1^n}{2} , \frac{\nabla \rho_1^{n+1} - \nabla \rho_1^n }{\Delta t})_h -  \kappa_{\rho_2 \rho_2} (\frac{\nabla \rho_2^{n+1} + \nabla \rho_2^n}{2} , \frac{\nabla \rho_2^{n+1} - \nabla \rho_2^n }{\Delta t})_h \\

- \kappa_{\rho_1 \rho_2}  [(\frac{\nabla \rho_1^{n+1} + \nabla \rho_1^n}{2}  , \frac{\nabla \rho_2^{n+1} - \nabla \rho_2^n }{\Delta t} )_h +  (\frac{\nabla \rho_2^{n+1} + \nabla \rho_2^n}{2} , \frac{\nabla \rho_1^{n+1} - \nabla \rho_1^n }{\Delta t})_h ],
\eea \label{eq:full_energy_dissipation3}
\een where we used lemma (\ref{lem1}). Combining (\ref{eq:full_energy_dissipation1}), (\ref{eq:full_energy_dissipation2}) and (\ref{eq:full_energy_dissipation3}), we obtain
\ben\bea{l}
\frac{E_h^{n+1} - E_h^n}{\Delta t}  = - 2 (\frac{1}{Re_s}, {\bf D}_h^{n+1/2} : {\bf D}_h^{n+1/2})_2 - (\frac{1}{Re_v} tr({\bf D}_h^{n+1/2}), tr({\bf D}_h^{n+1/2}) )_2 \\

- M_1 ( \nabla (\mu_1^{n+1/2} -  \mu_2^{n+1/2}),   \nabla (\mu_1^{n+1/2} -  \mu_2^{n+1/2}))_2 \leq 0,
\eea\een
provided $ M_1 \geq 0$.
Having established unconditional energy stability, we now turn to the solvability issue of the linear system of equations.

\subsection{Unique Solvability of the Fully Discrete, Linear Numerical Scheme}
\noindent \indent The linear system resulting from scheme (\ref{eq:full_discrete}) can be written into
\ben\bea{l}
\mathcal{A} \cdot X = \mathcal{G},
\eea\een
where $\mathcal{A}$ is the coefficient matrix of the system  given in Appendix, $X := $ $(\mu_1, \mu_2, u, v , q_1, \rho_1, \rho_2)$ is the solution of the linear system and the right hand term $\mathcal{G} = (g_1, g_2, g_3, g_4, g_5, g_6, g_7)^T$ denotes all the terms at the nth time step.

{\thm{Linear system (\ref{eq:full_discrete}) admits a unique solution. }}\\
{\bf Proof}:
To prove the well-posedness of the system \eqref{eq:full_discrete}, we only need to prove the corresponding homogeneous system admits only the zero solution.
We assume that there is  a solution $X = (\mu_1, \mu_2, u, v, q_1, \rho_1, \rho_2)$ such that $\mathcal{A} \cdot X$ = 0. Using \eqref{eq:full_discrete_0}, we   have
\ben\bea{l}
0 = (\mathcal{A} \cdot X, X)_2
=     M_1 ( \nabla (\mu_1 -  \mu_2),   \nabla (\mu_1 -  \mu_2))_2 + \frac{2}{\Delta t} [u, u]_{ew} +  \frac{2}{\Delta t} [v, v]_{ns} +  \frac{4}{\Delta t} (q_1, q_1)_2 \\
 + 2 (\frac{1}{Re_s}, {\bf D}_h : {\bf D}_h)_2 + (\frac{1}{Re_v} tr({\bf D}_h), tr({\bf D}_h) )_2
+    \frac{2}{\Delta t} [  \kappa_{\rho_1 \rho_1}( \nabla \rho_1 ,\nabla \rho_1)_h  +  \kappa_{\rho_2 \rho_2} (\nabla \rho_2 , \nabla \rho_2 )_h] \\

+  \frac{4}{\Delta t} \kappa_{\rho_1 \rho_2}  (\nabla \rho_1  , \nabla \rho_2)_h

\geq C ((\nabla \rho_1 , \nabla \rho_2)_h + ( \nabla \rho_2 , \nabla \rho_2)_h +  [u, u]_{ew} + [v, v]_{ns} + \|  q_1  \|^2_2 ),
\eea\een
where we used ${\bf K} > 0$, C is a positive constant and ${\bf D}_h$ is defined in \eqref{eq:D_h_definition}.
Thus, we obtain
\ben\bea{l}
D_x \rho_1 = D_y \rho_1 = 0, \quad D_x \rho_1 = D_y \rho_2 = 0, \quad u = v = 0, \quad q_1 = 0,
\eea\een
Based on linear system (\ref{eq:full_discrete_0}), we have
\ben\bea{l}
\mu_1 = \mu_2 = 0, \quad \rho_1 = \rho_2 = 0,
\eea\een
i.e. $X = {\bf 0}$. Thus,  linear system \eqref{eq:full_discrete} admits an unique solution.

\rem{A second order in time, energy stable BDF scheme can be developed as well, which will not be  presented here. }

\section{Numerical results and discussions}

\subsection{Accuracy  Test}

\noindent \indent We conduct a mesh refinement test to verify the convergence rate of the numerical scheme by considering  \eqref{eq:nondim} with a double-well bulk free energy
\ben\bea{l}
h(\rho_1, \rho_2, T) = \rho_1^2(\rho_1-1)^2  + \rho_2^2 (\rho_2 - 1)^2,
\eea\een
in a rectangular domain $\Omega = [0, 1] \times [0, 1]$. We use the following initial conditions
\ben\bea{l}
\rho_1(x, y, t = 0) = 0.5  + 0.01 cos(2 \pi x) ,\quad

\rho_2(x, y, t = 0) = 0.5 - 0.01 cos(2\pi x) , \quad

{\bf v} = (0, 0).
\eea\een
We denote the number of spatial grids as $N_x = N_y = N$, the time step as $\Delta t$. To test the convergence rate in time, we first fix $N = 256$ and vary the time step from $4 \times 10^{-3}$ to $0.125 \times 10^{-3}$ to calculate the $l_2$ norm of the difference between the numerical solutions obtained using consecutive step sizes at $T = 0.1$, i.e. $\| (\cdot)_{\Delta t}(T) - (\cdot)_{2 \Delta t}(T)  \|_{2}$. Then, we fix time step $\Delta t = 10^{-4}$, vary the spatial grid number from $8$ to $256$  and calculate the $l_2$ norm of the difference between the numerical solutions obtained using consecutive grid sizes at $T = 0.1$, i.e. $\| (\cdot)_{h}(T) - (\cdot)_{2 h}(T)  \|_{2}$. In both space and time, we calculate the convergence rate using $p = log_2 \Big ( \frac{  \| (\cdot )_{2h}(T) - (\cdot )_{4 h}(T) \|_2 }  {  \| (\cdot )_{h}(T) - (\cdot )_{2 h}(T)  \|_{2}  }\Big )$, where h is the mesh size in time or space. The refinement results are tabulated in Table \ref{table:error_test_relative_time} and Table \ref{table:error_test_relative_space}, respectively. We observe that the proposed scheme is indeed second-order accurate in both time and space for all variables.

\begin{table}
\begin{center}
\begin{tabular}{|c|c|c|c|c|c|c|c|} \hline
$\Delta$t & $\| {(\rho_1)}_{\Delta t} - {(\rho_1)}_{2\Delta t}  \|_2 $  & order & $\| {(\rho_2)}_{\Delta t} - {(\rho_2)}_{2 \Delta t}  \|_2 $  & order &  $\| {{(\bf u)}}_{\Delta t} - {(\bf u)}_{2 \Delta t}  \|_2$  & order \\ \hline \hline
 4 $\times 10^{-3} $ &    &  & &  & &   \\
2 $\times 10^{-3} $ & 0.5237   $\times 10^{-8}$&    & 0.5240  $\times 10^{-8}$& &0.1498 $\times 10^{-7} $  &    \\
1 $\times 10^{-3} $ & 0.1348   $\times 10^{-8}$& 1.96 & 0.1349   $\times 10^{-8}$& 1.96 &0.3806 $\times 10^{-8} $ & 1.98 \\
 0.5 $\times 10^{-3} $ &  0.3425   $\times 10^{-9}$&  1.98 &  0.3428   $\times 10^{-9}$&  1.98& 0.9594 $\times 10^{-9} $ & 1.99  \\
 0.25 $\times 10^{-3} $ &  0.8644    $\times 10^{-10}$& 1.99  &  0.8651    $\times 10^{-10}$&  1.99  & 0.2435 $\times 10^{-9} $  & 1.98 \\
 0.125 $\times 10^{-3} $ &  0.2129    $\times 10^{-10}$&  2.02  &  0.2130    $\times 10^{-10}$&  2.02 & 0.5779 $\times 10^{-10} $  & 2.08 \\ \hline
 \end{tabular}
\end{center}
\caption{Temporal refinement result for all variables. The model parameter values are chosen as $Re_s = 100 , Re_v = 300, M_{1} = 10^{-7}$, $\kappa_{\rho_1 \rho_1} = \kappa_{\rho_2 \rho_2} = 10^{-4}$, $\kappa_{\rho_1 \rho_2} = \kappa_{\rho_2 \rho_1} = 0$.}

\label{table:error_test_relative_time}
\end{table}

 \begin{table}
\begin{center}
\begin{tabular}{|c|c|c|c|c|c|c|c|} \hline
N & $\| {(\rho_1)}_{h} - {(\rho_1)}_{2 h}  \|_2 $ & order & $\| {(\rho_2)}_{h} - {(\rho_2)}_{2 h}  \|_2 $  & order &  $\| {({\bf u})}_{h} - {(\bf u)}_{2 h}  \|_2 $  & order \\ \hline \hline
8  &                                           &               &                                          &             &                                          &  \\
 16 &  0.2281    $\times 10^{-5}$&    &  0.2282   $\times 10^{-5}$&   & 0.2676 $\times 10^{-7} $  & \\
32 &  0.3417    $\times 10^{-6}$&  1.74  &  0.3421    $\times 10^{-6}$&  1.74  & 0.3487 $\times 10^{-8} $  &1.85 \\
64 &  0.4452    $\times 10^{-7}$&  1.94  &  0.4457    $\times 10^{-7}$&  1.93  & 0.4607 $\times 10^{-9} $  &1.88\\
128 &  0.5623    $\times 10^{-8}$&  1.98  &  0.5631    $\times 10^{-8}$&  1.99  & 0.5898 $\times 10^{-10} $  &1.95\\
256 &  0.7050    $\times 10^{-9}$&  2.00  &  0.7059    $\times 10^{-9}$&  2.00  & 0.7444 $\times 10^{-11} $  &1.98\\  \hline
\end{tabular}
\end{center}
\caption{Spatial refinement result for all variables. The model parameter values are chosen as $Re_s = 1 , Re_v = 3, M_{1} = 10^{-3}$, $\kappa_{\rho_1 \rho_1} = \kappa_{\rho_2 \rho_2} = 10^{-4}$, $\kappa_{\rho_1 \rho_2} = \kappa_{\rho_2 \rho_1} = 0$.}
\label{table:error_test_relative_space}
\end{table}

\subsection{Phase Separation in binary compressible viscous fluids}

\noindent \indent To demonstrate stability and efficiency of the new scheme, we simulate phase separation dynamics using system (\ref{eq:nondim}) with the Flory-Huggins mixing  energy
\ben\bea{l}
h(\rho_1, \rho_2, T) = \frac{k_B T}{m}\rho( \frac{1}{N_1}  \frac{\rho_1}{\rho} ln \frac{\rho_1}{\rho}+\frac{1}{N_2}   \frac{\rho_2}{ \rho} ln \frac{\rho_2}{\rho}+\chi \frac{\rho_1\rho_2}{\rho^2}),
\eea\label{eq:fh}
\een
where we choose the characteristic scales so that  $\frac{k_B T}{m} = 1$ in the  simulation, $N_1, N_2$ are the polymerization indices and $\chi$ is the mixing coefficient, which are given in the simulation  by
\ben\bea{l}
N_1 = N_2 = 1, \qquad \chi = 2.5.
\eea\een

\begin{figure}
\centering
\subfigure[ Flory-Huggins mixing energy density function with respect to $\frac{\rho_1}{\rho}$ ]{\includegraphics[width=0.4\textwidth]{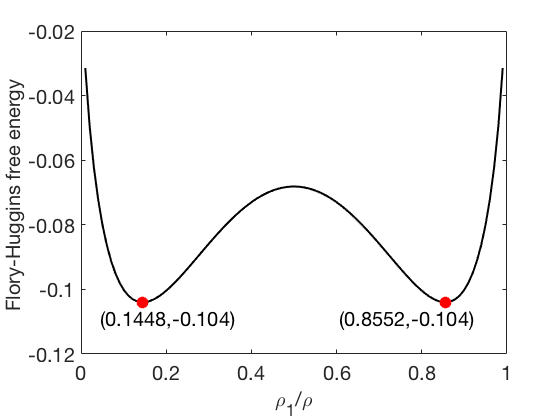}}
\subfigure[ Unstable mode ]{\includegraphics[width=0.4\textwidth]{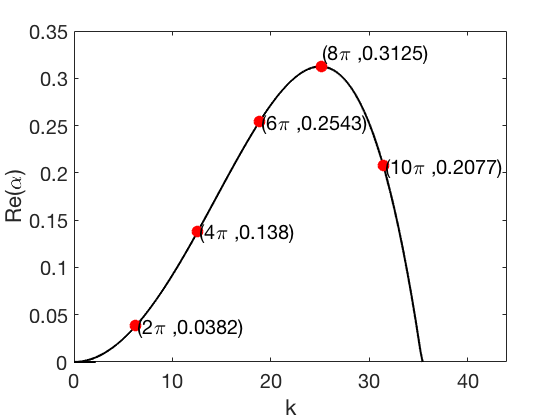}}
\caption{(a) Flory-Huggins mixing energy density function with respect to the mass density fraction $\frac{\rho_1}{\rho}$ at the chosen parameter values. The two minima are labeled by dots in the curve. (b) The unstable mode with parameter values: $N_1 = N_2 = 1, \chi = 2.5$, $M_1 = 10^{-3},  Re_s = 100, Re_v = 300$, $\kappa_{\rho_1 \rho_1} = \kappa_{\rho_2 \rho_2} = 0.0004 $, $\kappa_{\rho_1 \rho_2} = 0 $. }
\label{fig:F_H_1D_unstable_mode}
\end{figure}

The plot of this energy density with the chosen parameter values as a function of $\frac{\rho_1}{\rho}$ is shown in \ref{fig:F_H_1D_unstable_mode}-(a). The other dimensionless model parameters are chosen as follows
\ben\bea{l}
M_1 = 10^{-3}, \quad Re_s = 100, \quad Re_v = 300, \quad \kappa_{\rho_1 \rho_1} = \kappa_{\rho_2 \rho_2} =  4 \times 10^{-4}, \quad \kappa_{\rho_1 \rho_2} = 0.
\eea\een

In order to identify the spinodal decomposition that drives the phase separation in the binary polymer blend, we conduct a simple linear stability analysis on the hydrodynamic phase field model.
We note that this compressible model admits a family of constant solutions:
\ben
{\bf v}={\bf 0},\quad \rho_1=\rho_1^0,\quad \rho_2=\rho_2^0,\label{csol}
\een
where $\rho_1^0, \rho_2^0$ are constants. We perturb the constant solutions with a normal mode as follows:
\ben
{\bf v}= \epsilon e^{\alpha t + i {\bf k} \cdot {\bf x}} {\bf v}^{c}, \quad \rho_1=\rho_1^0 + \epsilon e^{\alpha t + i {\bf k} \cdot {\bf x}} {\rho_1}^{c},\quad  \rho_2=\rho_2^0 +\epsilon e^{\alpha t + i {\bf k} \cdot {\bf x}} {\rho_2}^{c},
\een
where $\epsilon$ is a small parameter, representing the magnitude of the perturbation, and ${\bf v}^c, \rho_1^c, \rho_2^c$ are constants, $\alpha $ is the growth rate, and $\bf k$ is the wave number of the  perturbation. Without loss of generality, we limit our study to 1 dimensional perturbation in $\bf k$ in the $(x,y)$ plane. Substituting these perturbations into the equations in \eqref{eq:nondim} and truncating the equations at order $O(\epsilon)$, we obtain the linearized equations.
The dispersion equation of the linearized equation system of the compressible model \cite{Zhao&W2018} is given by an algebraic equation of $\alpha$:
\ben\bea{l}
(\eta^0 k^2 + \alpha \rho^0)\{\alpha^3 \rho_0 + \alpha^2 k^2[\eta + \rho^0 M_{1} (h_{\rho_1 \rho_1} + \kappa_{\rho_1 \rho_1}k^2) +  \rho^0 M_{1} (h_{\rho_2 \rho_2} + \kappa_{\rho_2 \rho_2}k^2) ] \\

- \alpha^2 k^2[ 2 \rho^0 M_{1}(h_{\rho_1 \rho_2} + \kappa_{\rho_1 \rho_2}k^2)] + \alpha [{\bf p}^T \cdot {\bf C} \cdot {\bf p} + {\bf p}^T\cdot  {\bf K} \cdot {\bf p}k^2 ]k^2 \\

+\alpha \eta M_{1} [   (h_{\rho_1 \rho_1} + \kappa_{\rho_1 \rho_1}k^2)  +  (h_{\rho_2 \rho_2} + \kappa_{\rho_2 \rho_2}k^2) - 2  (h_{\rho_1 \rho_2} + \kappa_{ \rho_1 \rho_2}k^2) ]k^4
\\

+ k^4 M_{1} (\rho_1^0 + \rho_2^0)^2

 [ (h_{\rho_1\rho_1} + \kappa_{\rho_1 \rho_1}k^2) (h_{\rho_2\rho_2} + \kappa_{\rho_2 \rho_2}k^2) -  (h_{\rho_1\rho_2} + \kappa_{\rho_1 \rho_2}k^2)^2]\}
= 0,
\eea \label{eq:Algebraic_com1}
\een
where $\eta = 2\eta^0 + \overline{\eta}^0$, ${\bf p} = (\rho_1^0, \rho_2^0)^T$. In the following, we set  $\rho_1^0 = \rho_2^0 = 0.5$. $\bf K$ is the coefficient matrix of the conformational entropy and $\bf C$ is the Hessian of bulk  energy $h(\rho_1, \rho_2, T)$ with respect to  $\rho_1$ and $\rho_2$,
\ben
{\bf K} = \left(
\bea{cc}
\kappa_{\rho_1 \rho_1 }  & \kappa_{\rho_1 \rho_2}  \\
\kappa_{\rho_1 \rho_2}  & \kappa_{\rho_2 \rho_2 }
\eea
\right),
{\bf C} = \left(
\bea{cc}
h_{\rho_1 \rho_1}  & h_{\rho_1 \rho_2}  \\
h_{\rho_1 \rho_2}  & h_{\rho_2 \rho_2}  \\
\eea
\right).
\een
Obviously, $\alpha = - \frac{\eta^0 }{\rho^0} k^2 < 0$ is a solution of the dispersion equation \eqref{eq:Algebraic_com1}, which contributes a stable mode. To resolve the other modes, we use numerical calculations. Based on the model parameters listed above, we obtain only one unstable mode,  shown in Figure \ref{fig:F_H_1D_unstable_mode}-(b). This unstable mode is dominated by the mixing energy of the model, independent of hydrodynamics of the model. Next, we will numerically simulate phase separation phenomena due to the unstable perturbation on the constant steady state without and with hydrodynamics to show how hydrodynamics can affect the path of phase separation and its outcome.

\subsubsection{Phase separation without hydrodynamics}

\noindent \indent Based on unstable mode shown in Figure (\ref{fig:F_H_1D_unstable_mode}-b), we add a 1D perturbation with wave number $k = 10 \pi$ to the steady state and observe its ensuing nonlinear dynamics. Since the eigenvector corresponding to the unstable mode shown in Figure (\ref{fig:F_H_1D_unstable_mode}-b) is ($\rho_1^c, \rho_2^c$) = (1, -1), we impose the initial conditions specifically as follows
\ben\bea{l}
\rho_1(x, y, t=0) = 0.5  + 0.005 \times cos(10 \pi y), \qquad

\begin{figure}
\centering
\subfigure[ $\rho_1$ at t = 0  ]{\includegraphics[width=0.225\textwidth]{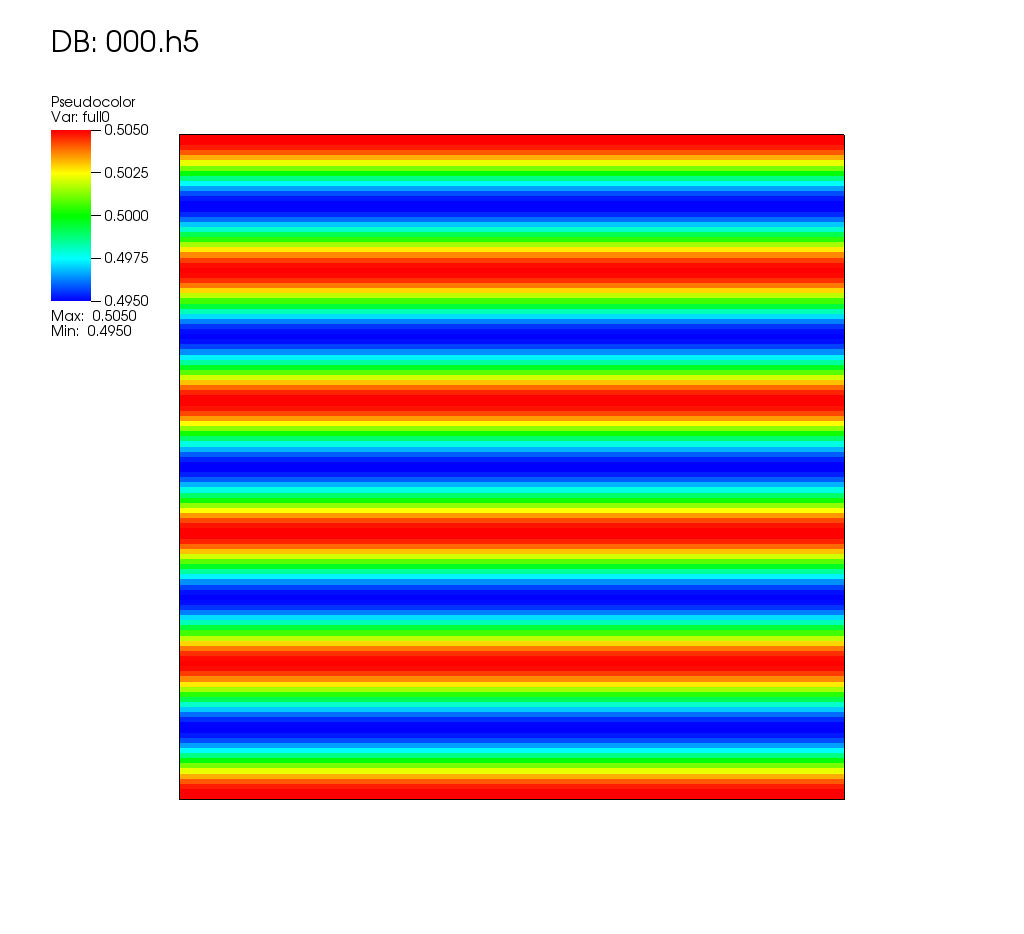}}
\subfigure[ $\rho_1$ at t = 300  ]{\includegraphics[width=0.225\textwidth]{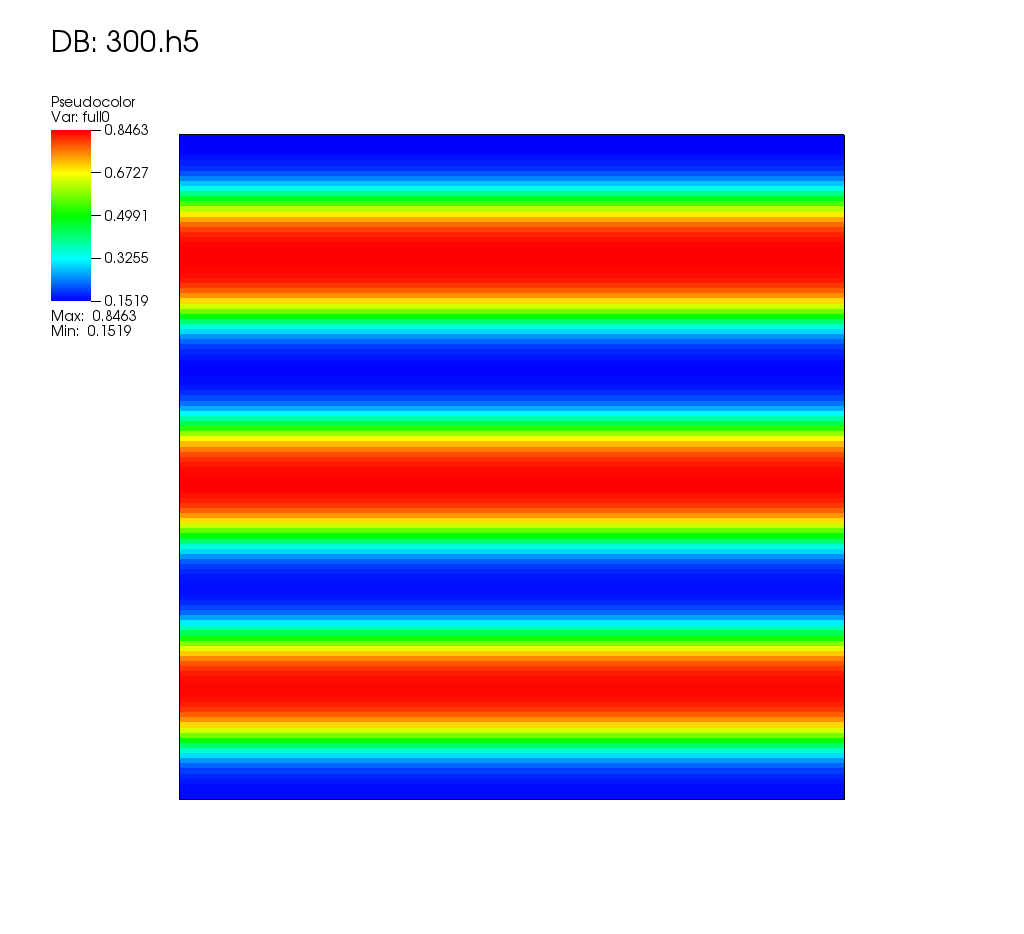}}
\subfigure[ $\rho_1$ at t = 2000  ]{\includegraphics[width=0.225\textwidth]{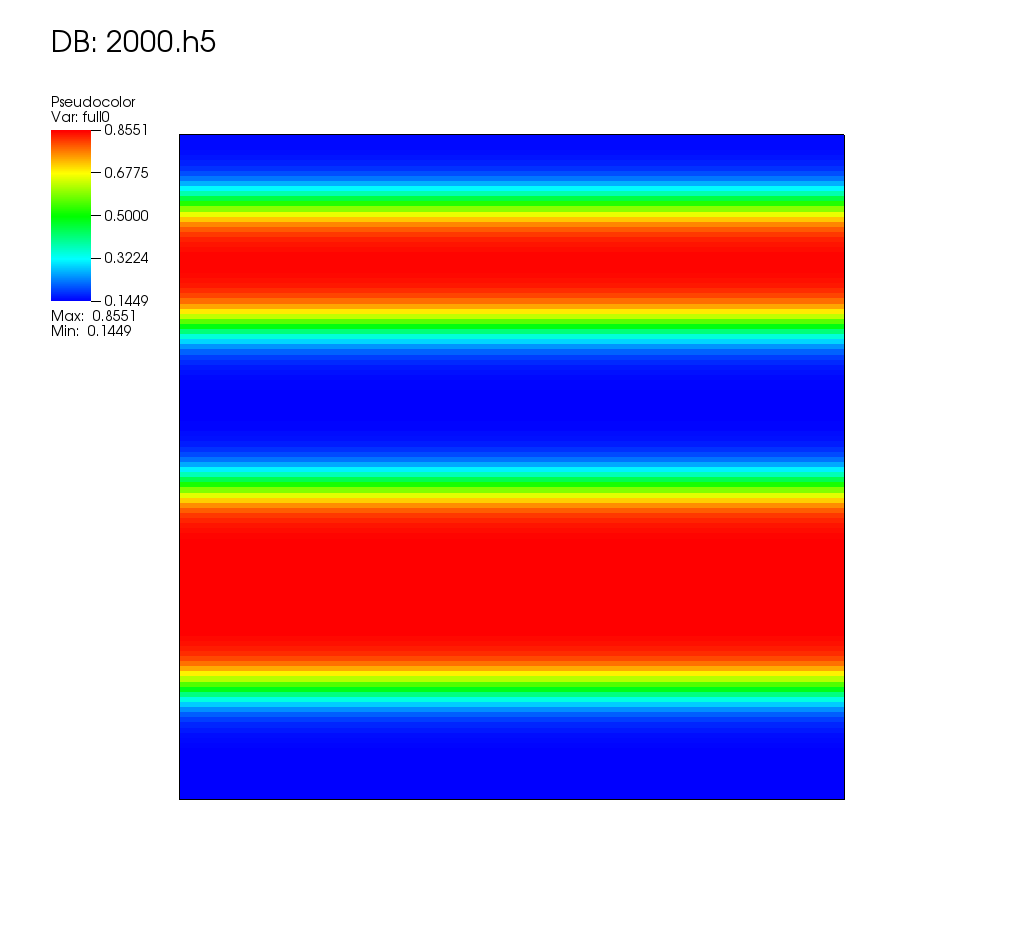}}
\subfigure[ $\rho_1$ at t = 15000  ]{\includegraphics[width=0.225\textwidth]{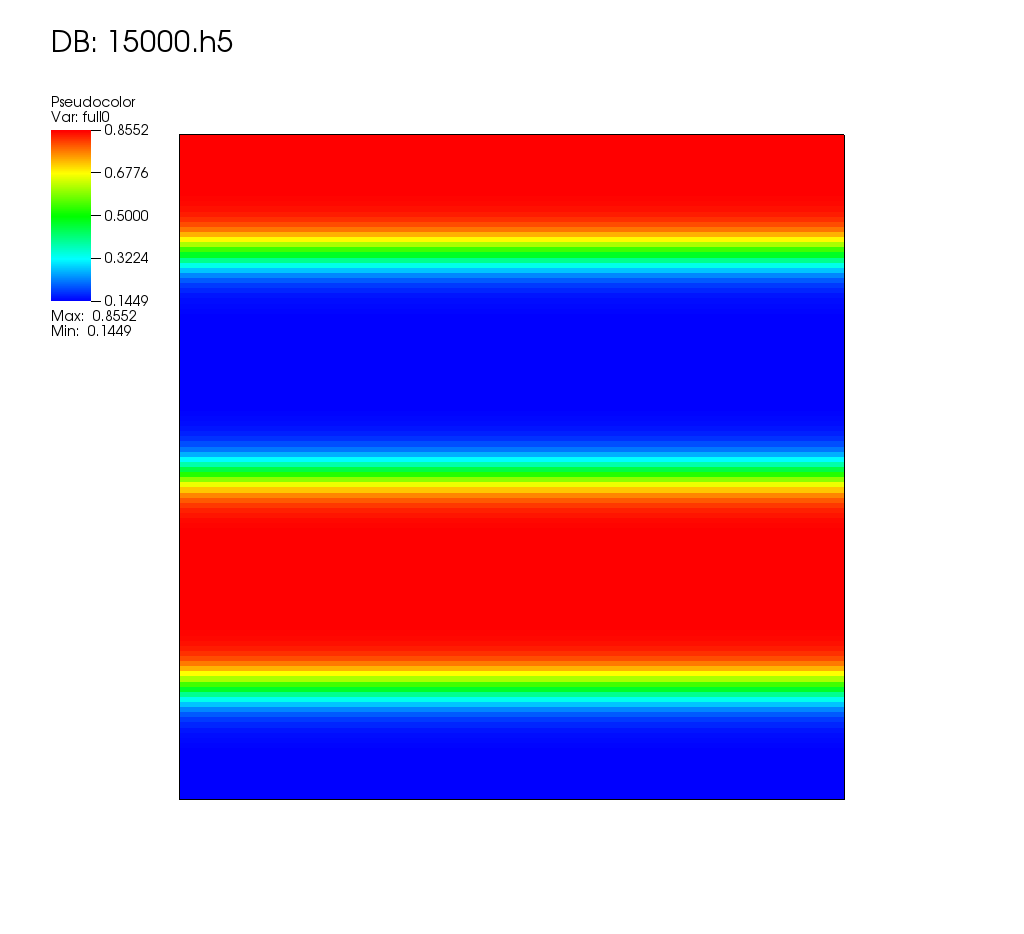}}

\subfigure[ Total Energy ]{\includegraphics[width=0.3\textwidth]{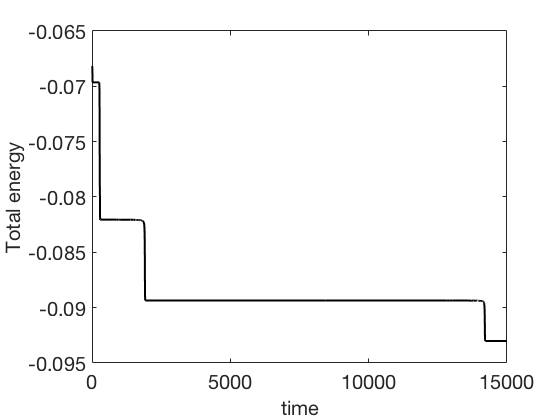}}
\subfigure[ Difference of total mass of the component 1 with its initial total mass ]{\includegraphics[width=0.3\textwidth]{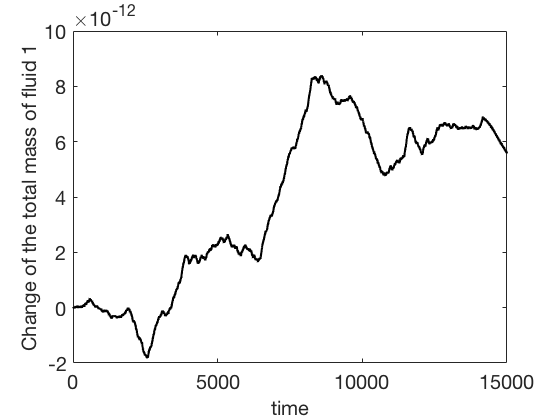}}
\subfigure[ Difference of total mass of the component 2 with its initial total mass ]{\includegraphics[width=0.3\textwidth]{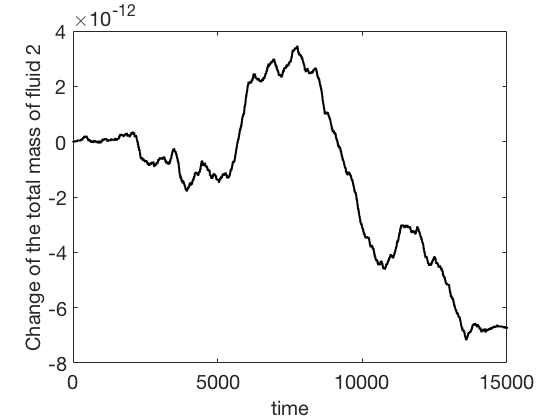}}
\caption{(a-d) Snapshots of $\rho_1$ at different times as solutions of system (\ref{eq:nondim}) with the Flory-Huggins mixing  energy (\ref{eq:fh}) without hydrodynamics. (e) The total free energy of system (\ref{eq:nondim}). Two major coarsening events bring the phase of the binary system into the final state shown in (d). $\rho_2$ is given by $1-\rho_1$. The total mass of both phases are conserved as shown in (f-g). }
\label{fig:rho1_without_velo}
\end{figure}

\rho_2(x, y, t=0) = 0.5  - 0.005 \times cos(10 \pi y).
\eea\een
Since $\rho_1 + \rho_2 = 1$ in the thermodynamic model without hydrodynamics,  we show the phase behavior of $\rho_1$ only.
The time evolution  of $\rho_1$ at a few selected times are depicted in Figure \ref{fig:rho1_without_velo}. Firstly, we observe that the growth rate of the numerical solutions $\rho_1$ near the equilibrium state is $\alpha = 0.2077$, which matches with the linear stability analysis result shown in Figure (\ref{fig:F_H_1D_unstable_mode}-b). In the long-time  behavior, we observe that $\rho_1$ develops small-scale structures and then coarsens to large-scale structures eventually. In Figure \ref{fig:rho1_without_velo}, we show numerical solutions at several time slots and the corresponding  total energy up to $t=15000$. The system goes through three coarsening events which are captured by the phase morphology at different times shown as well as the total energy evolution in Figure \ref{fig:rho1_without_velo}. The outcome at the end of the computation is a four-band structure.

\subsubsection{ Phase separation with hydrodynamics}

\noindent \indent When hydrodynamics is coupled with the thermodynamical phase evolution, its role must show up somewhere. Here, we investigate how hydrodynamic impact on phase separation dynamics.  Since the eigenvector corresponding to the unstable mode shown in Figure (\ref{fig:F_H_1D_unstable_mode}-b) is ($\rho_1^c, \rho_2^c$) = (1, -1, 0), we adopt the same initial conditions for $\rho_1$ and $\rho_2$ as before and a zero velocity condition:
\ben\bea{l}
\rho_1(x,y,0)= 0.5  + 0.005 \times cos(10 \pi y), \quad

\rho_2(x,y,0) = 0.5  - 0.005 \times cos(10 \pi y), \quad

\begin{figure}
\centering
\subfigure[ $\rho_1$ at t = 0  ]{\includegraphics[width=0.24\textwidth]{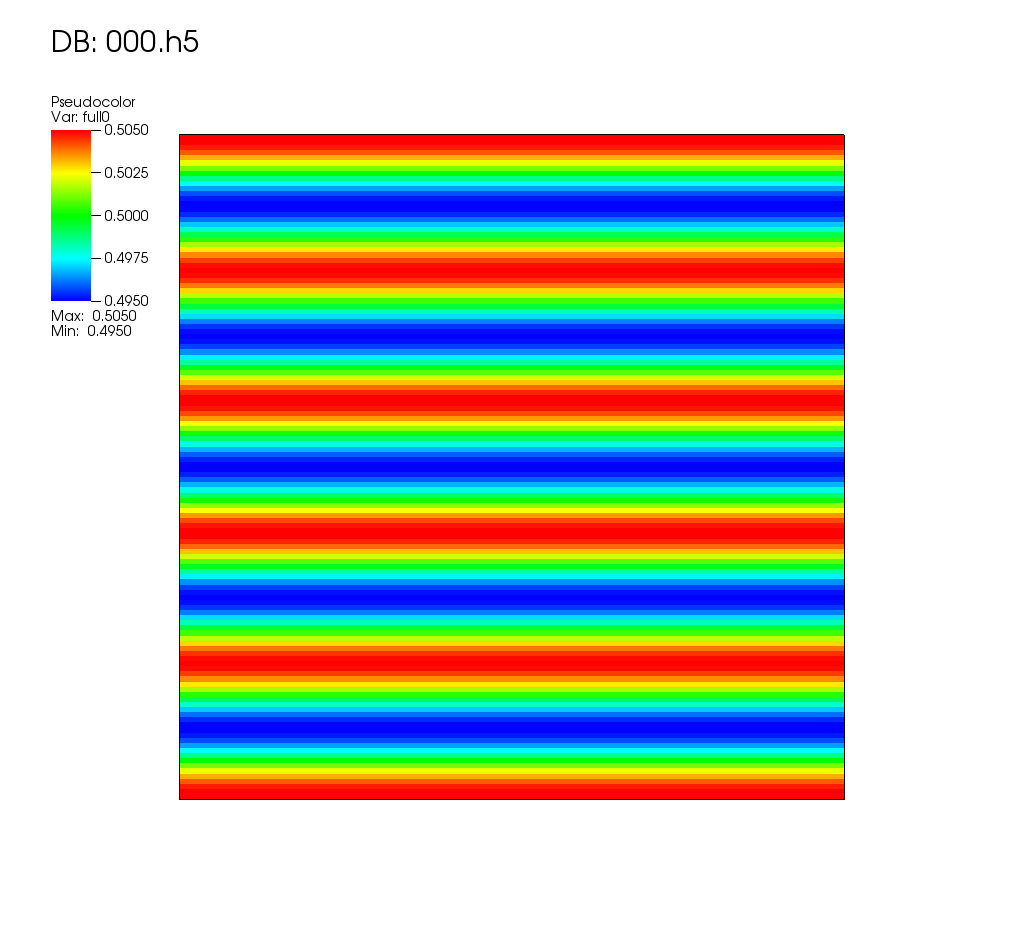}}
\subfigure[ $\rho_1$ at t = 50  ]{\includegraphics[width=0.24\textwidth]{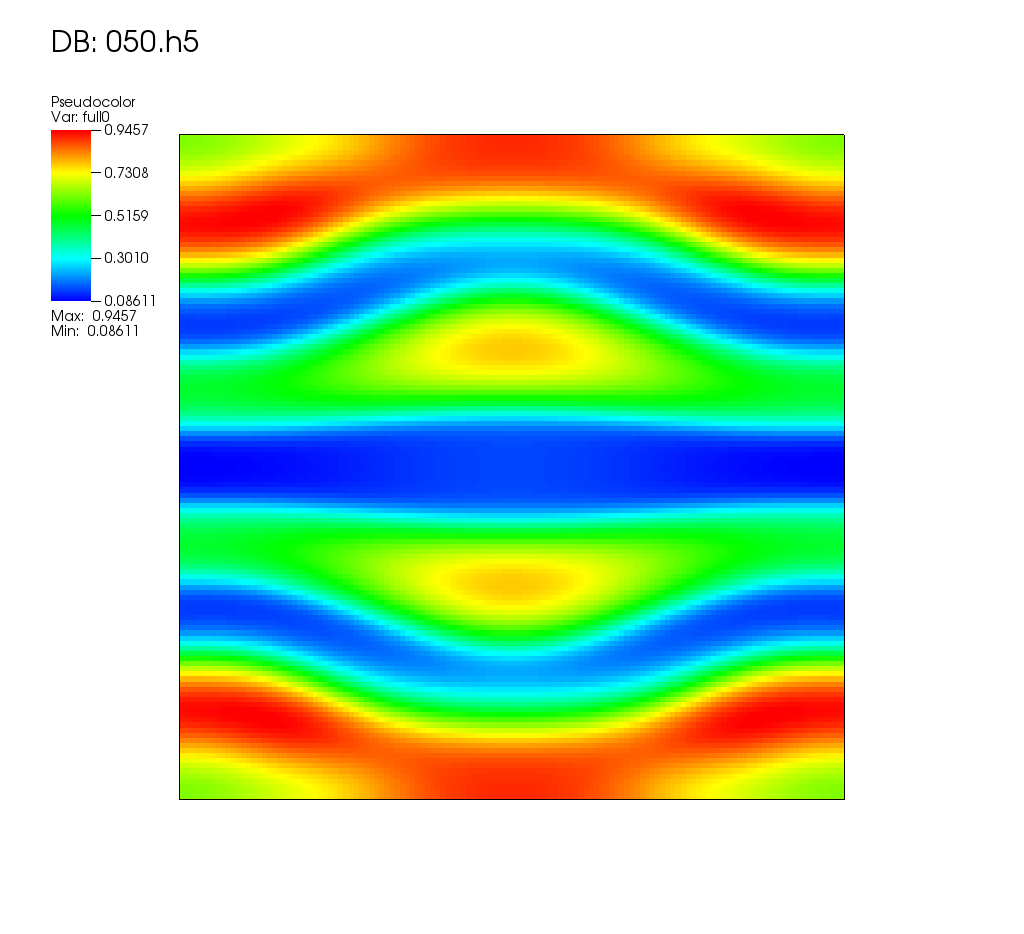}}
\subfigure[ $\rho_1$ at t = 100  ]{\includegraphics[width=0.24\textwidth]{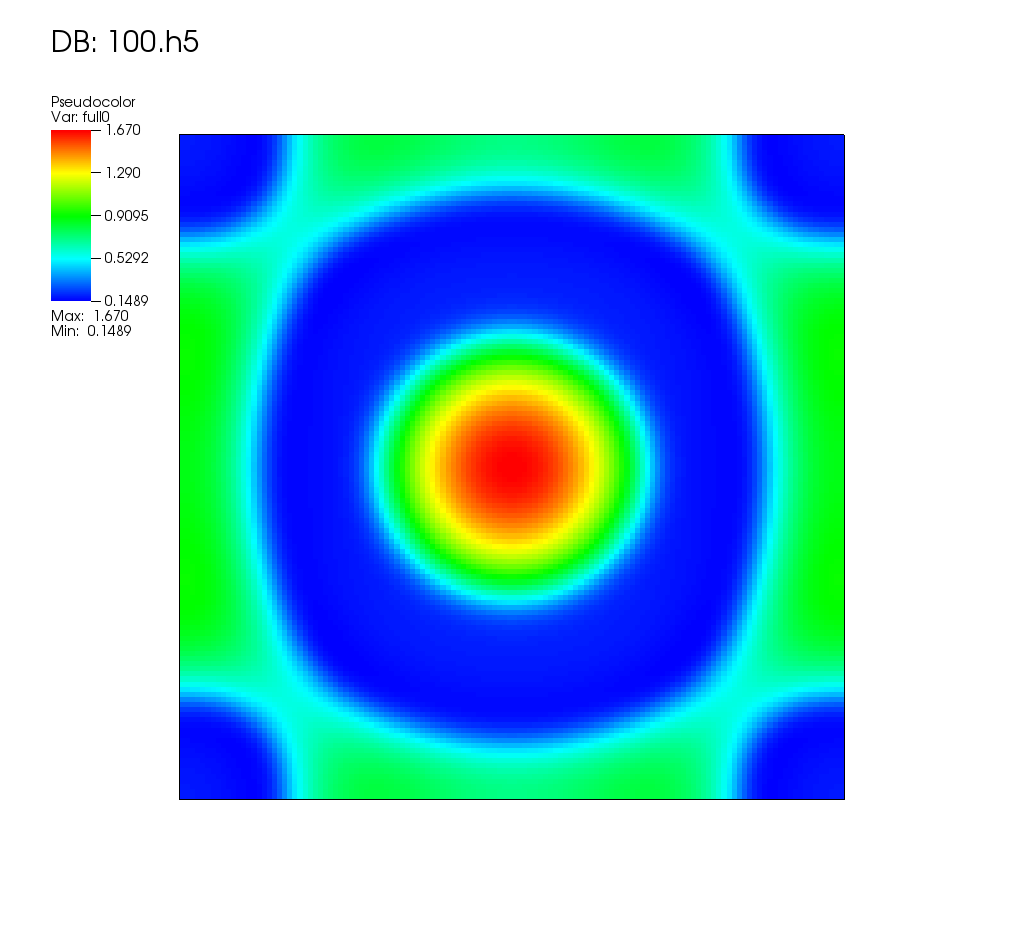}}
\subfigure[ $\rho_1$ at t = 150 ]{\includegraphics[width=0.24\textwidth]{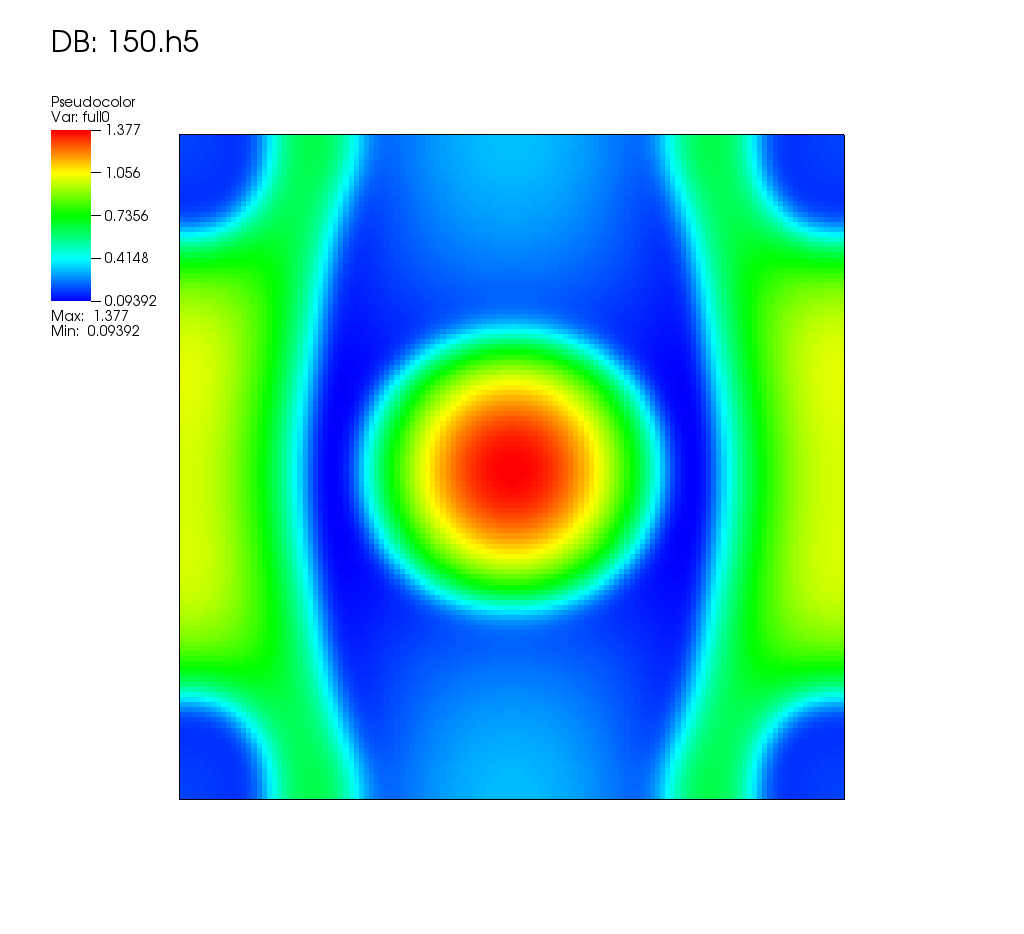}}
\subfigure[ $\rho_1$ at t = 200  ]{\includegraphics[width=0.24\textwidth]{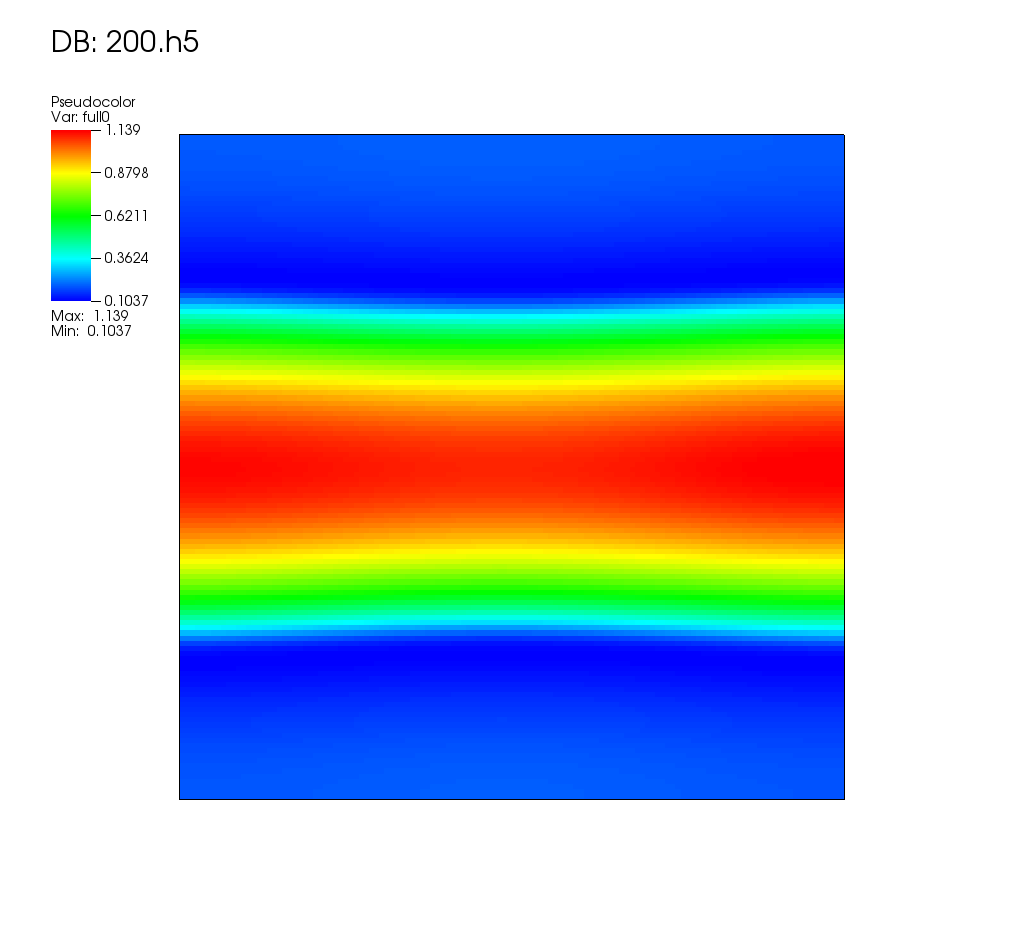}}
\subfigure[ $\rho_1$ at t = 400 ]{\includegraphics[width=0.24\textwidth]{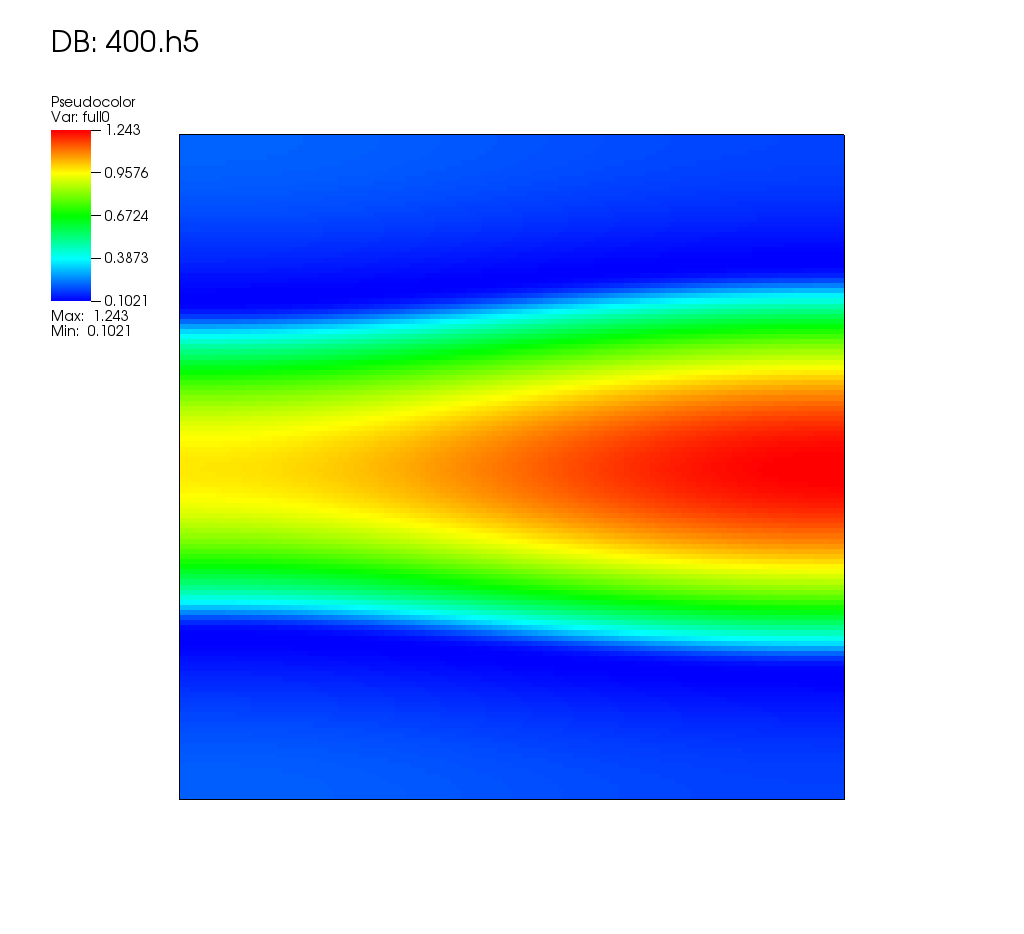}}
\subfigure[ $\rho_1$ at t = 600 ]{\includegraphics[width=0.24\textwidth]{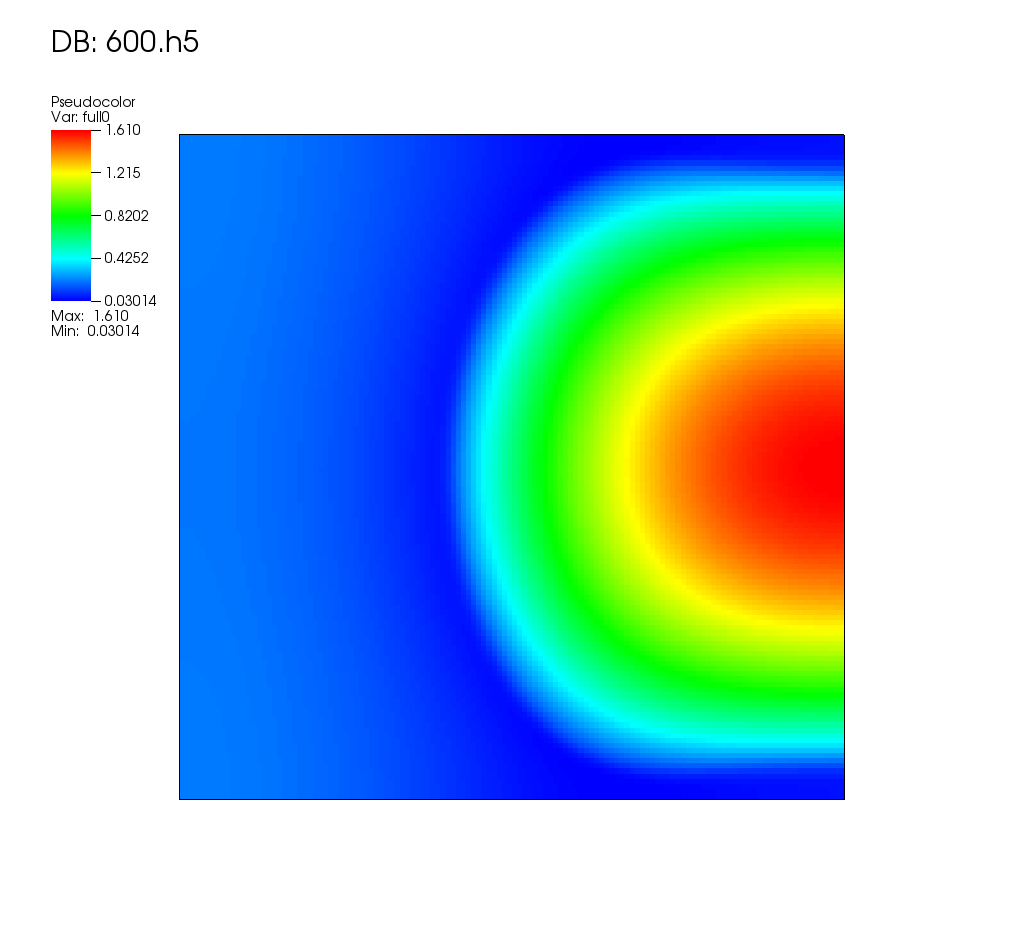}}
\subfigure[ $\rho_1$ at t = 1400  ]{\includegraphics[width=0.24\textwidth]{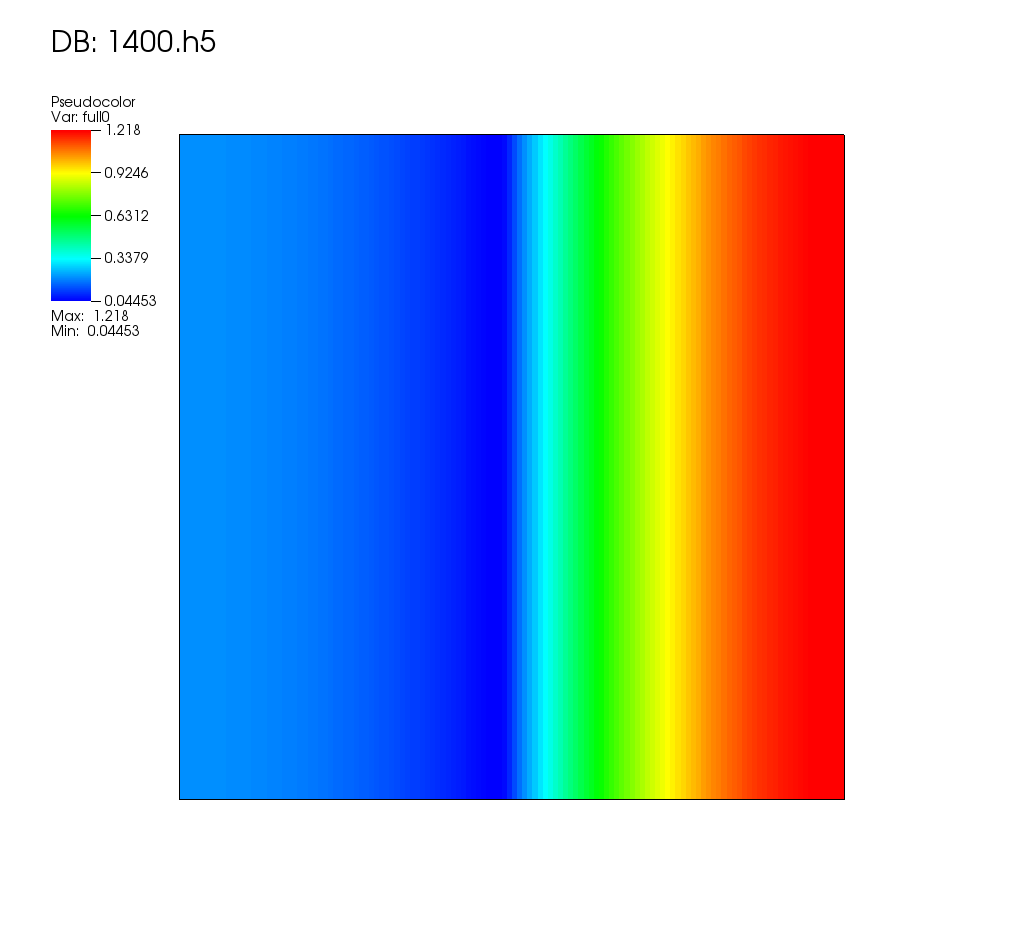}}

\subfigure[ Total Energy ]{\includegraphics[width=0.32\textwidth]{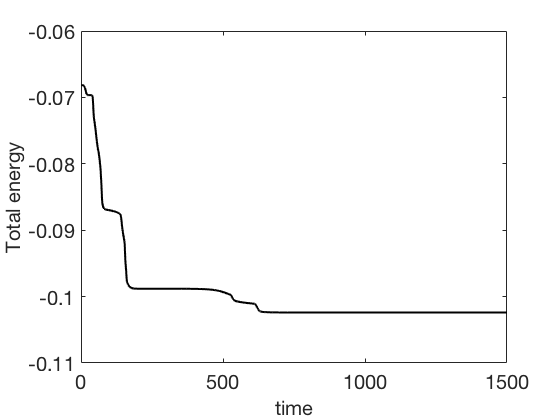}}
\subfigure[ Difference of total mass of the component 1 with its initial total mass ]{\includegraphics[width=0.32\textwidth]{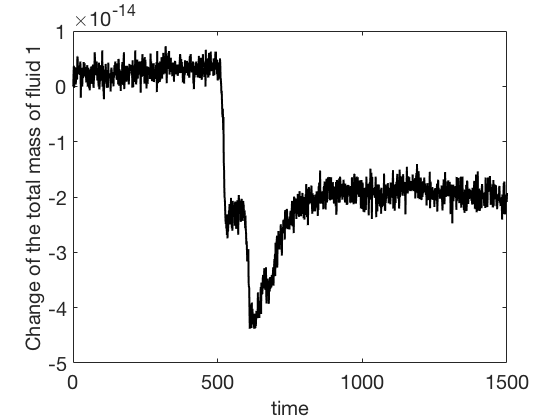}}
\subfigure[ Difference of total mass of the component 2 with its initial total mass ]{\includegraphics[width=0.32\textwidth]{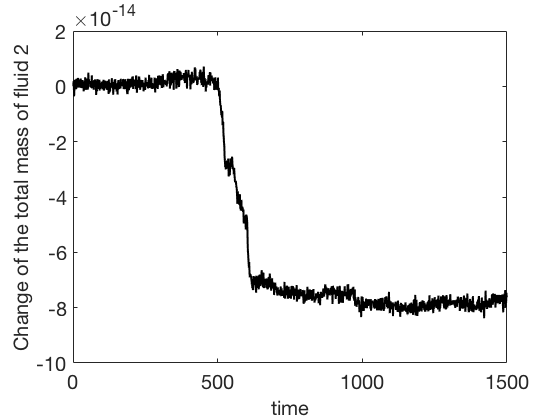}}
\caption{(a-h) Snapshots of $\rho_1$ at different times as a solution of system (\ref{eq:nondim}) with the Flory-Huggins mixing energy (\ref{eq:fh}) and hydrodynamic interaction. (i) Total energy of the system (\ref{eq:nondim}) with the Flory-Huggins bulk free energy (\ref{eq:fh}); (j, k) Difference of the total mass of component 1 and 2 compared with the initial mass, indicating mass conservation of both phases in the simulation.}
\label{fig:rho1}
\end{figure}

\begin{figure}
\centering
\subfigure[ $\rho_2$ at t = 0  ]{\includegraphics[width=0.24\textwidth]{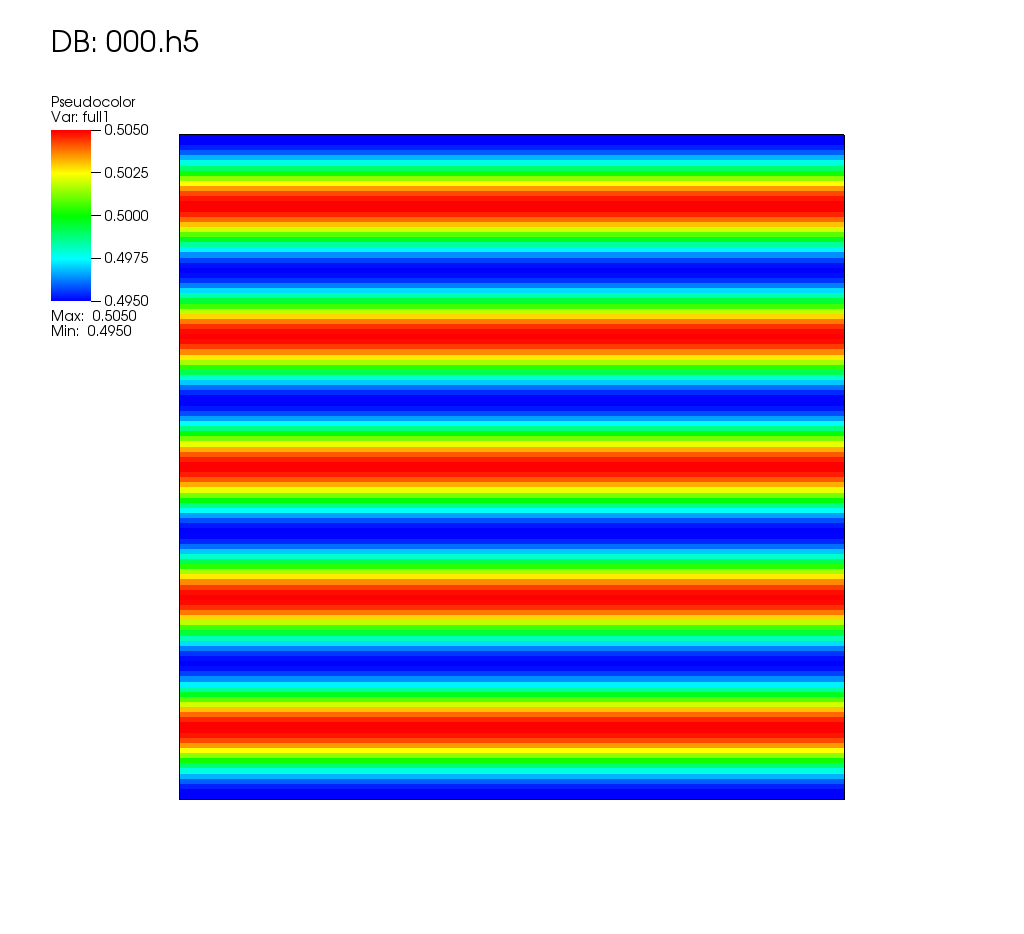}}
\subfigure[ $\rho_2$ at t = 50  ]{\includegraphics[width=0.24\textwidth]{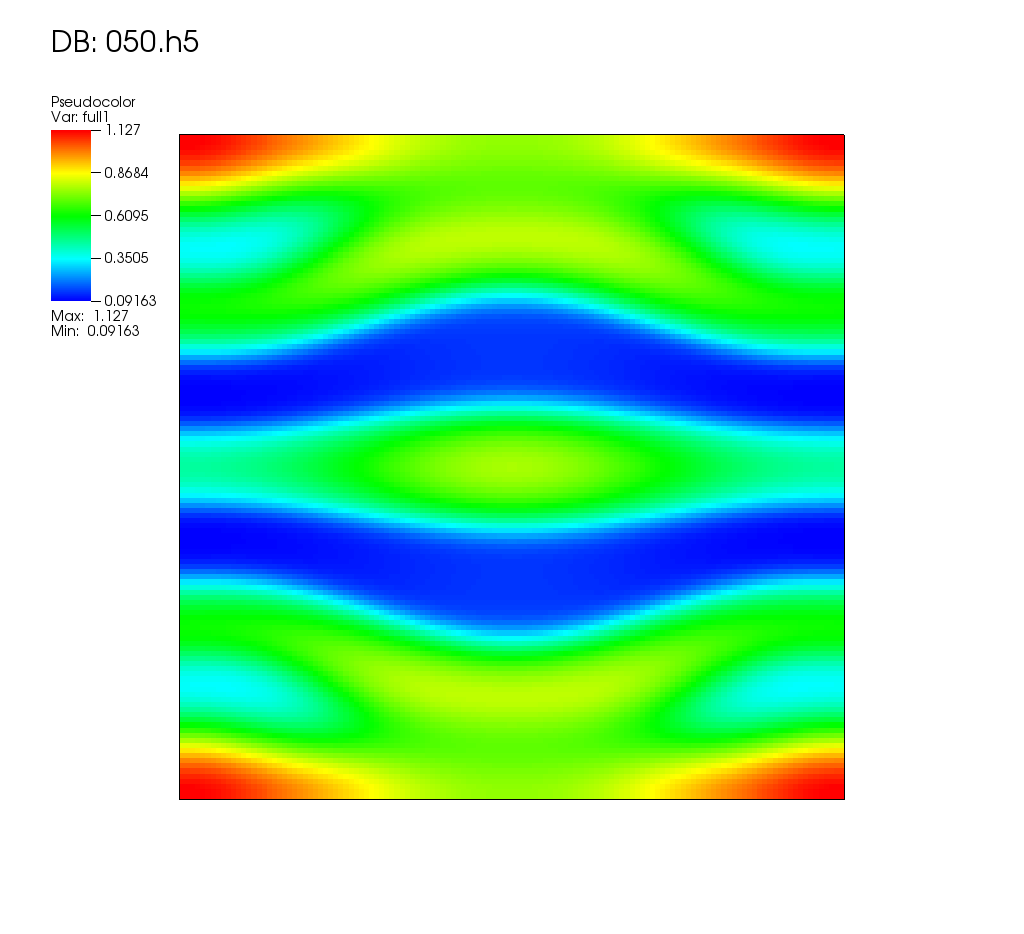}}
\subfigure[ $\rho_2$ at t = 100  ]{\includegraphics[width=0.24\textwidth]{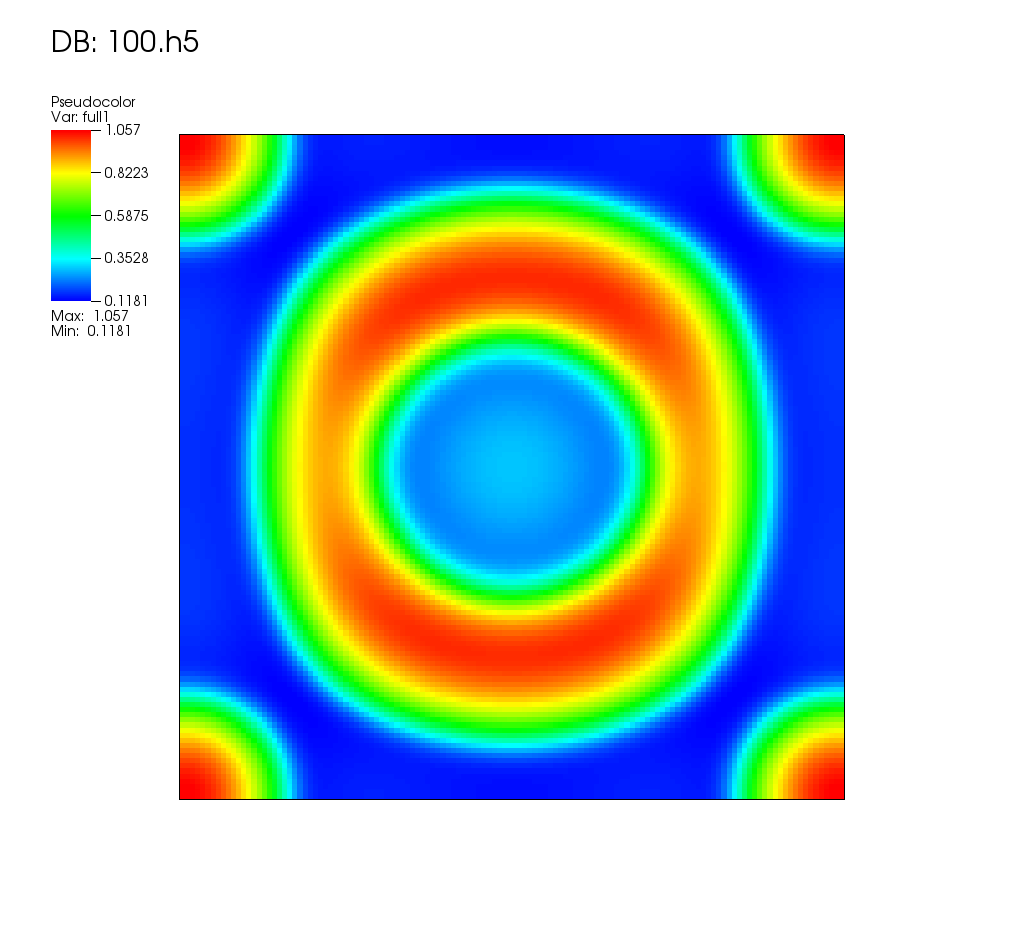}}
\subfigure[ $\rho_2$ at t = 150 ]{\includegraphics[width=0.24\textwidth]{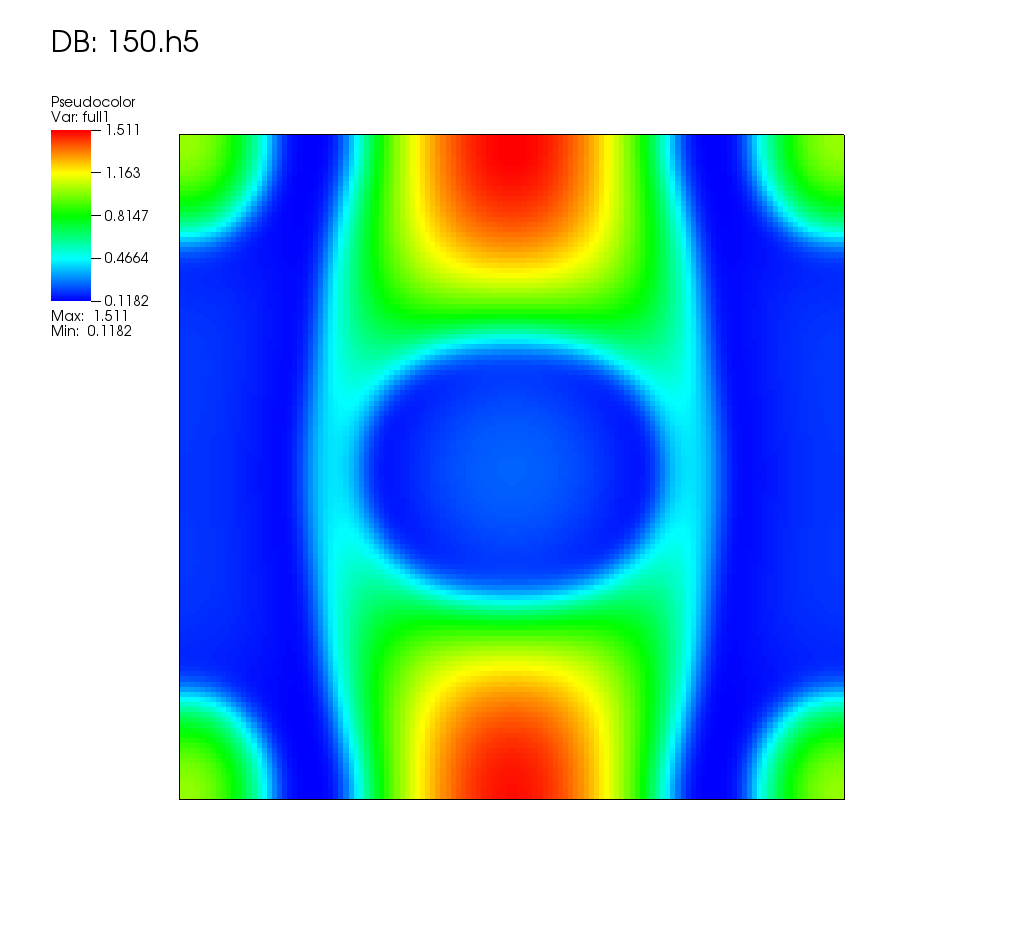}}
\subfigure[ $\rho_2$ at t = 200  ]{\includegraphics[width=0.24\textwidth]{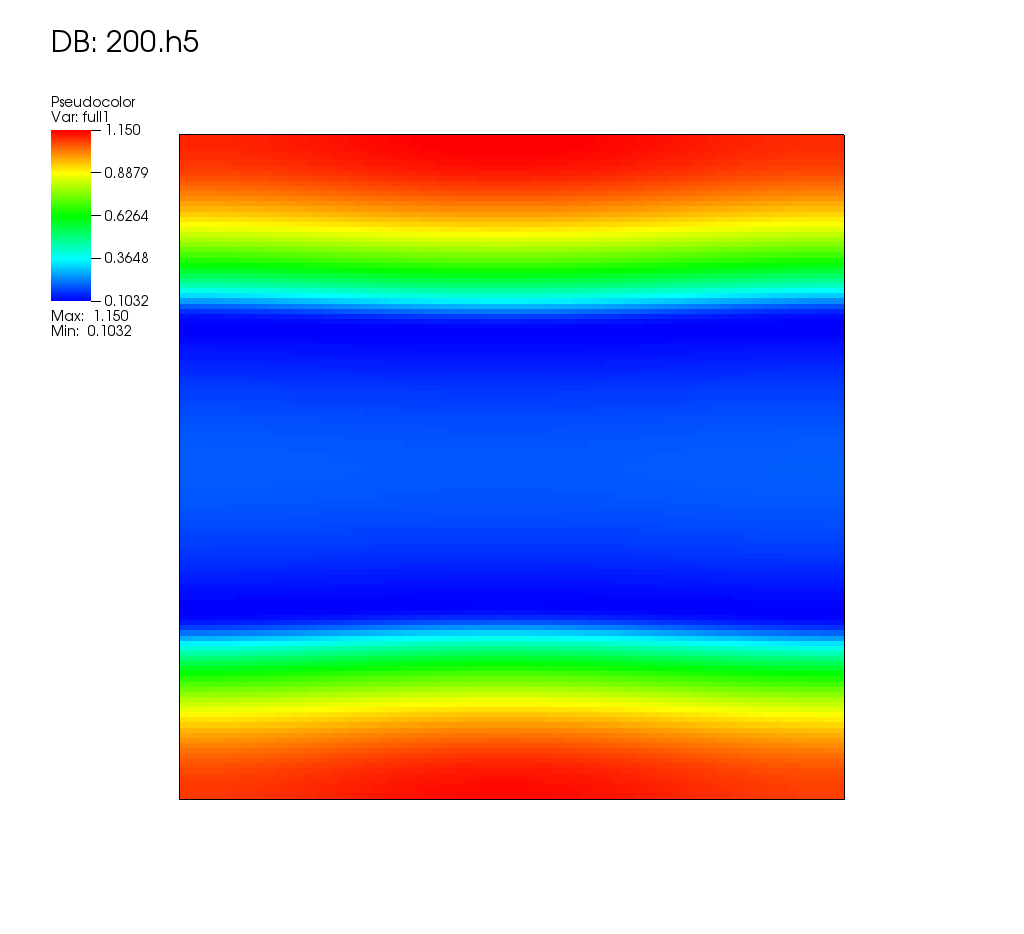}}
\subfigure[ $\rho_2$ at t = 400 ]{\includegraphics[width=0.24\textwidth]{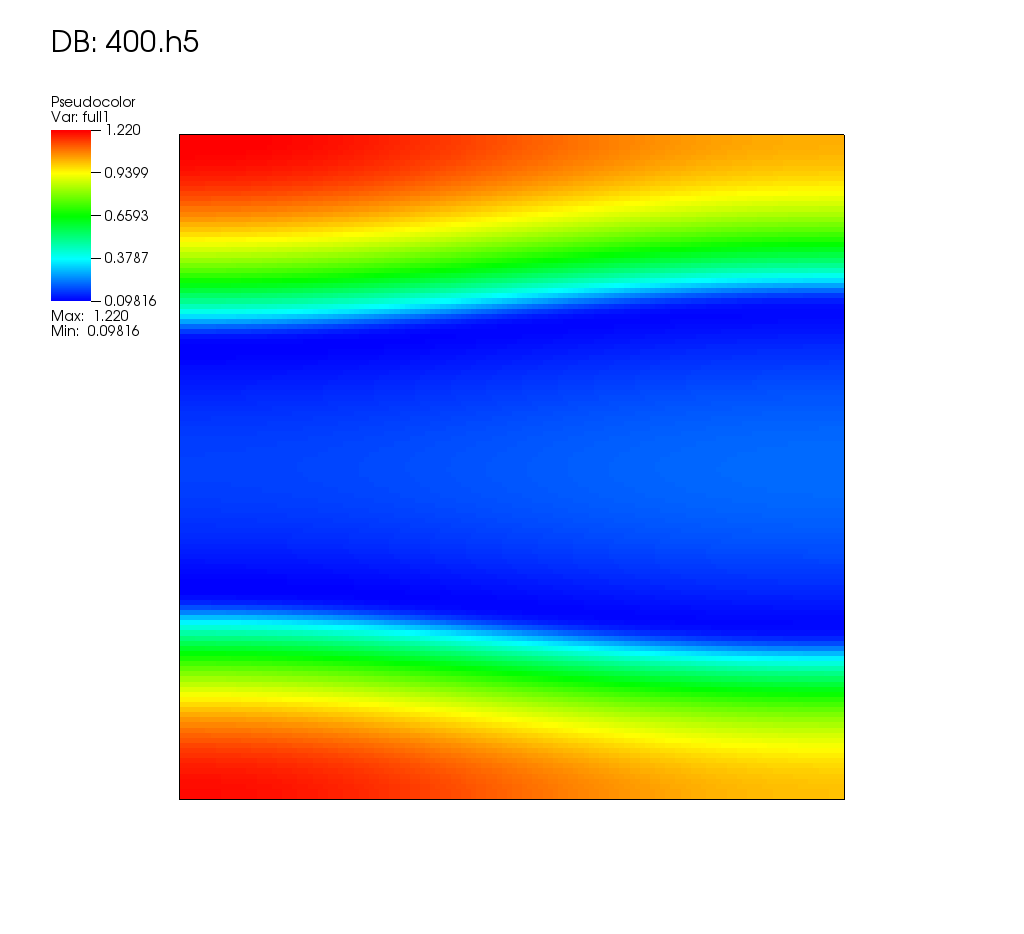}}
\subfigure[ $\rho_2$ at t = 600 ]{\includegraphics[width=0.24\textwidth]{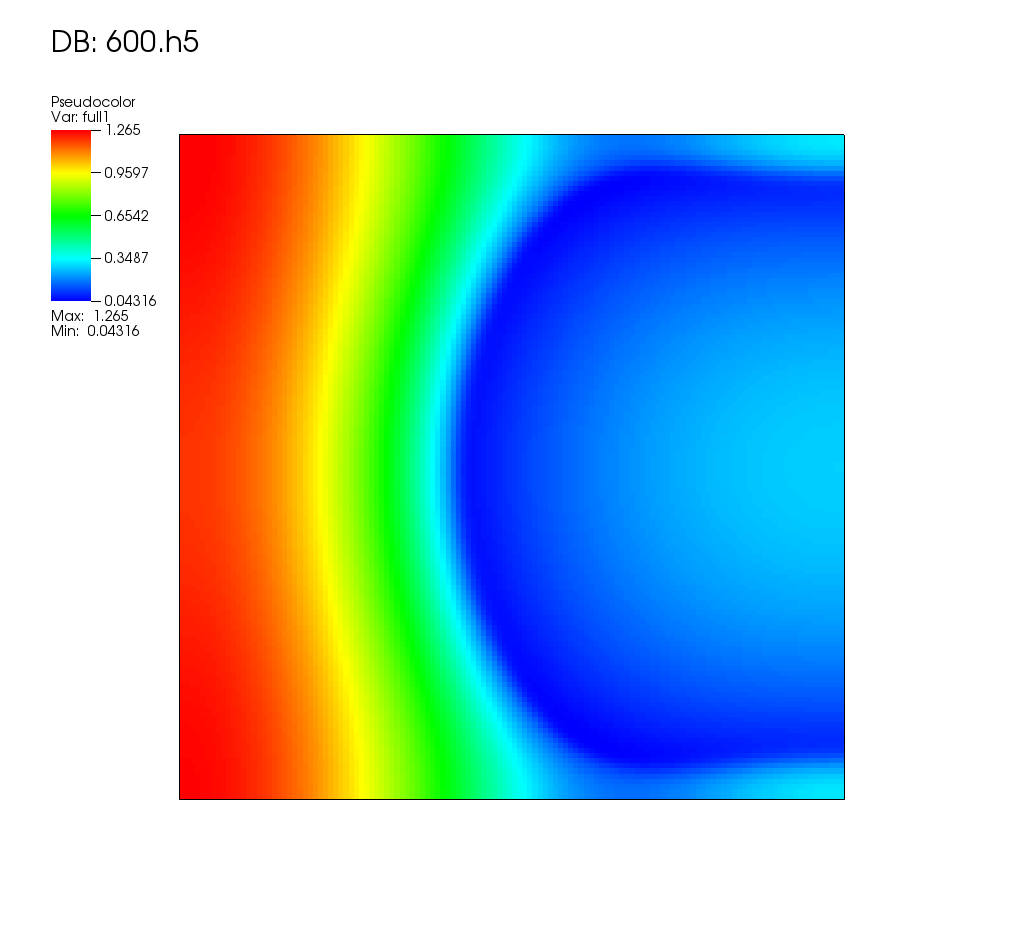}}
\subfigure[ $\rho_2$ at t = 1400  ]{\includegraphics[width=0.24\textwidth]{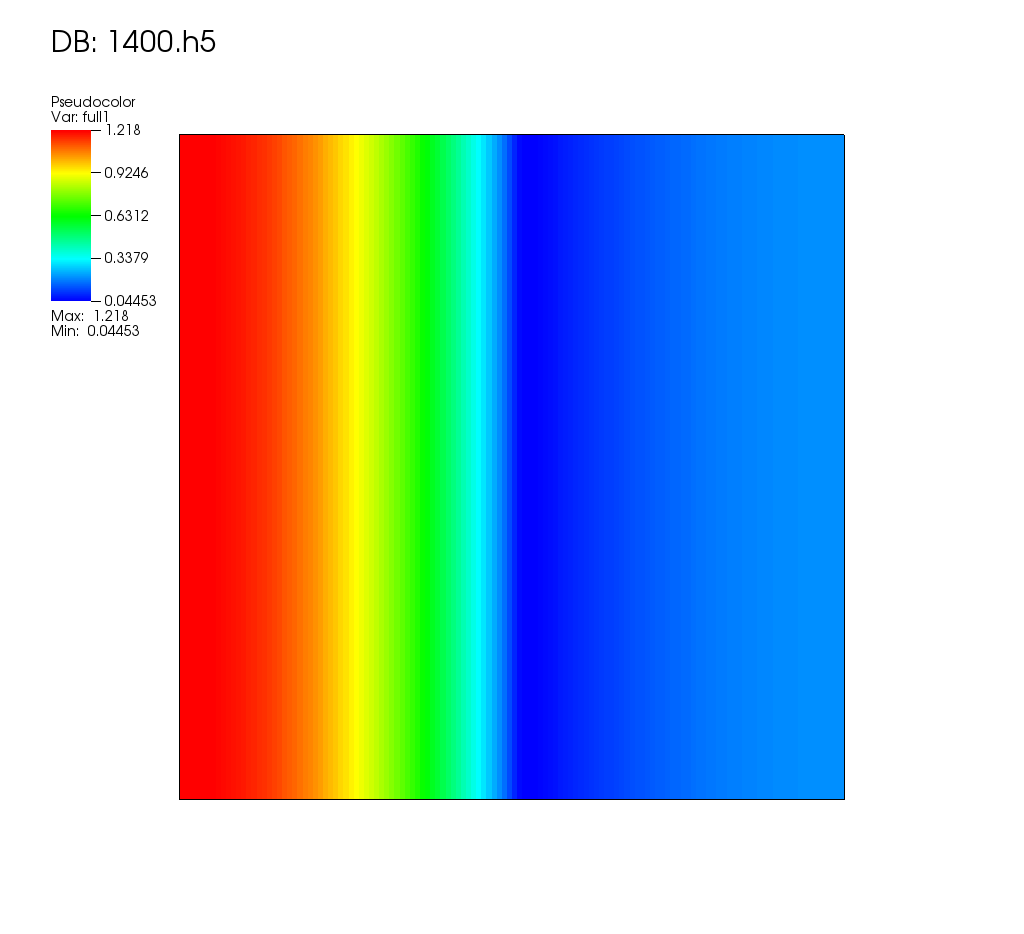}}
\subfigure[ ${\bf v}$ at t = 0  ]{\includegraphics[width=0.24\textwidth]{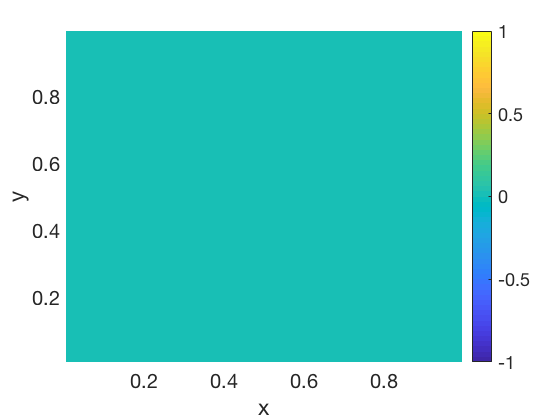}}
\subfigure[ ${\bf v}$ at t = 50  ]{\includegraphics[width=0.24\textwidth]{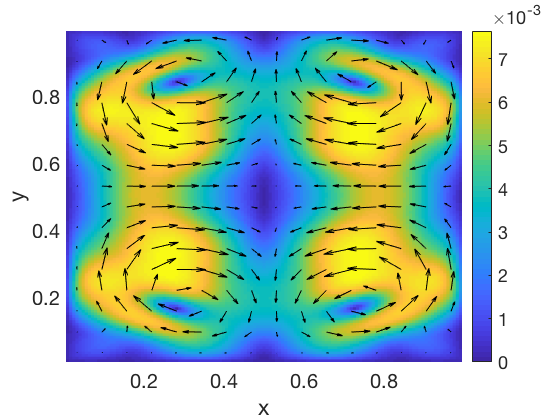}}
\subfigure[ ${\bf v}$ at t = 100  ]{\includegraphics[width=0.24\textwidth]{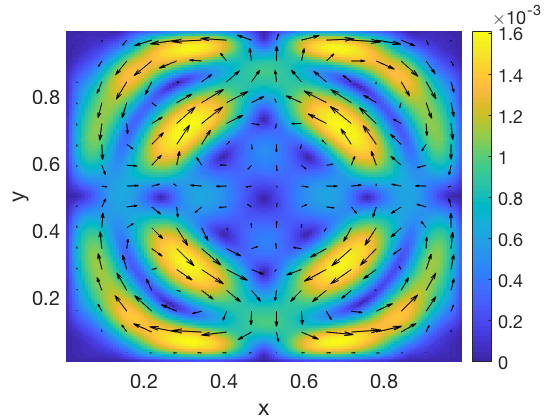}}
\subfigure[ ${\bf v}$ at t = 150 ]{\includegraphics[width=0.24\textwidth]{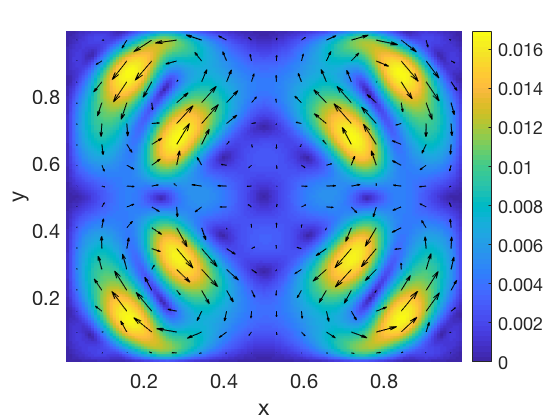}}
\subfigure[ ${\bf v}$ at t = 200  ]{\includegraphics[width=0.24\textwidth]{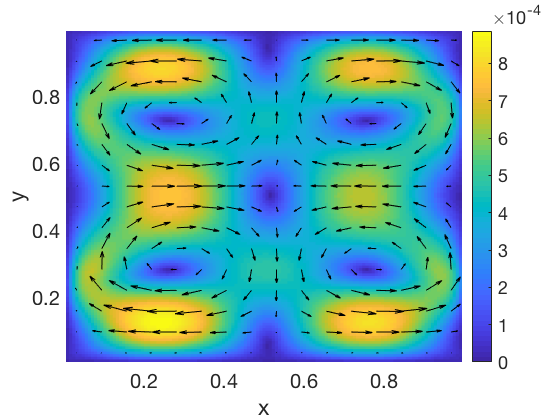}}
\subfigure[ ${\bf v}$ at t = 400 ]{\includegraphics[width=0.24\textwidth]{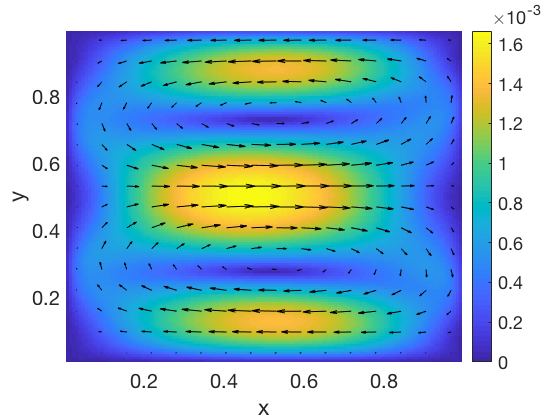}}
\subfigure[ ${\bf v}$ at t = 600 ]{\includegraphics[width=0.24\textwidth]{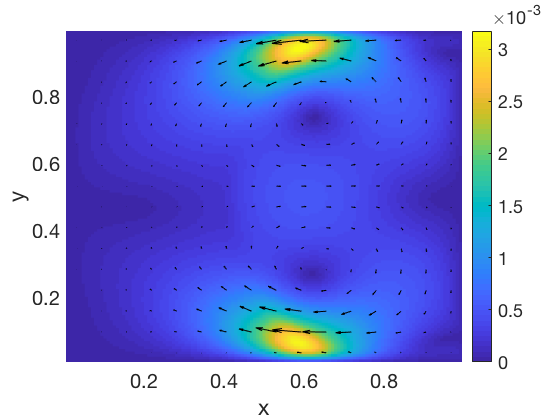}}
\subfigure[ ${\bf v}$ at t = 1400  ]{\includegraphics[width=0.24\textwidth]{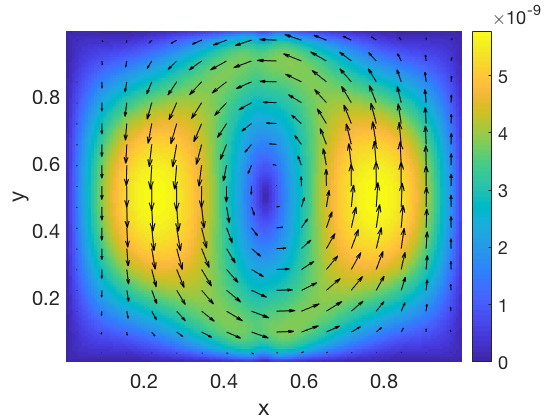}}

\caption{(a-h) Snapshots of $\rho_2$ at different times as a solution of system (\ref{eq:nondim}) with the Flory-Huggins mixing energy given in (\ref{eq:fh}) and hydrodynamic interaction. (i-p) Snapshots for velocity field ${\bf v} = (v_1, v_2)$ at different times. Weak flows are present due to hydrodynamic effect during the phase evolution. The nontrivial velocity leads to different phase morphology in the end compared to the case without hydrodynamic interaction at the end of our simulation. }
\label{fig:rho2}
\end{figure}

{\bf v}(x,y,0) = (0, 0).
\eea\een
When hydrodynamics is considered, the local total mass density $\rho$ is no longer spatially homogeneous anymore. However,  phase separation goes on as shown in Figure \ref{fig:rho1} and Figure \ref{fig:rho2}. In Figure \ref{fig:rho1}, we observe that the total energy of the system is dissipative and the total mass of component 1 and 2 are conserved in the domain globally. The velocity field in the domain is plotted at the selected times. Some vorticities form and disperse eventually as the phase morphology approaches a steady state. The induced nontrivial velocity field promotes the transport of materials and mixing across the domain leading to a two-band structure phase morphology eventually, which is a global energy stable state. In contrast, the final phase morphology developed in the phase separation without hydrodynamics may have only reached a local energy stable state, which can be explained by the comparison of the total energy evolutions shown in Figure (\ref{fig:rho1_without_velo}-e) and Figure (\ref{fig:rho1}-i), respectively. This tells us that hydrodynamics indeed changes local densities, the path of phase evolution and even the final energy steady states of fluid mixtures. This is alarming, indicating that hydrodynamic effects are instrumental in determining the correct spatial phase diagram for the binary fluid mixture. The total energy in the solution with hydrodynamics is smaller than that without it. So, hydrodynamics in a binary compressible fluid flow promotes fluid mixing and thereby speeds up phase separation.

\subsection{Dynamics of gas-liquid Mixtures}

\noindent \indent The compressible fluid model has many applications in the petroleum industry, where mixtures of non-hydrocarbons and hydrocarbons are abundant, such as in petroleum reservoirs or natural gas pipelines. Understanding their thermodynamic and hydrodynamic properties can help one to improve petroleum quality and yield significantly.

In the past, several equations of state had been developed to describe the relation among state variables (e.g. the volume, pressure and temperature) under a given set of physical conditions for compressible fluids. The Peng-Robinson equation of state (PR-EOS) \cite{Peng_Robinson} is one of the popular ones, which has been successfully applied to thermodynamic and volumetric calculations in both industries and academics. Specifically, PR-EOS provides a reasonable accuracy near the critical point, which makes it a good choice for gas-condensate systems in the petroleum industry. For this reason, we adopt it in a hydrocarbon mixture of methane and n-decane to show the performance of our model and numerical scheme in simulating hydrodynamics of the hydrocarbon mixtures. Many properties of the mixture can be studied by our mathematical model, such as mass adsorption of one component in the mixture on the interface between two phases near the equilibrium state, surface tension and even  verification of mixing rules in the mixture. In this example, we will focus on hydrodynamics of a hydrocarbon mixture with an unstable gas-liquid interface and study the mass adsorption phenomena at the interface from the point of view of the free energy near an equilibrium state.

The free energy density function derived from PR-EOS reads
\ben\bea{l}
f = f_b + h({\bf n}, T),
\eea \label{eq:Peng_Robinsoon_total_energy}
\een
where $f_b = \frac{1}{2} \sum_{i,j=1}^N c_{i,j} \nabla n_i  \cdot \nabla n_j$ is the conformational energy. The bulk free energy density function $h({\bf n}, T)$ is given in \eqref{eq:P_R}.

\rem{Since $f^{ideal}$ changes rapidly near the origin which may introduce singularity in numerical simulations, we regularize this term near the origin   as follows
\ben\bea{l}
f^{ideal} =
\begin{cases}
RT n_i(ln(\epsilon) - 1) + RT(\frac{1}{2 \epsilon} n_i^2 - \frac{\epsilon}{2 }) ,  \qquad   if \quad  n_i < \epsilon,  \\
RT n_i (ln(n_i)-1), \hbox{otherwise,}
\end{cases}
\eea\een
where $\epsilon>0$.
Corresponding to the modification, the chemical potential is changed to
\ben\bea{l}
\mu^{ideal} =
\begin{cases}
RT (ln(\epsilon) - 1) + RT(\frac{1}{ \epsilon} n_i) ,  \qquad   if \quad  n_i < \epsilon,  \\
RT ln (n_i), \hbox{otherwise.}
\end{cases}
\eea\een}

\begin{figure}
\centering
\subfigure[ Initial condition of n-decane ($C_{10}H_{22}$) ]{\includegraphics[width=0.3\textwidth]{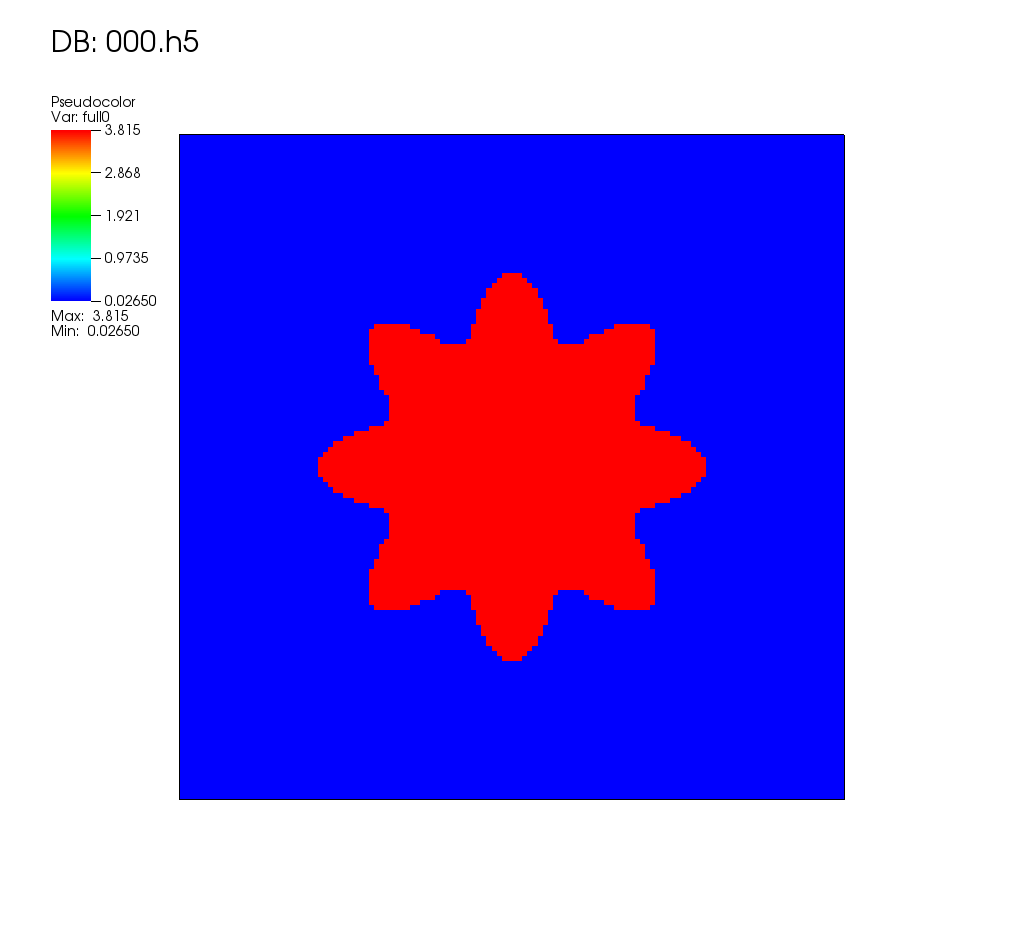}}
\subfigure[ Initial condition of methane ($CH_4$)]{\includegraphics[width=0.3\textwidth]{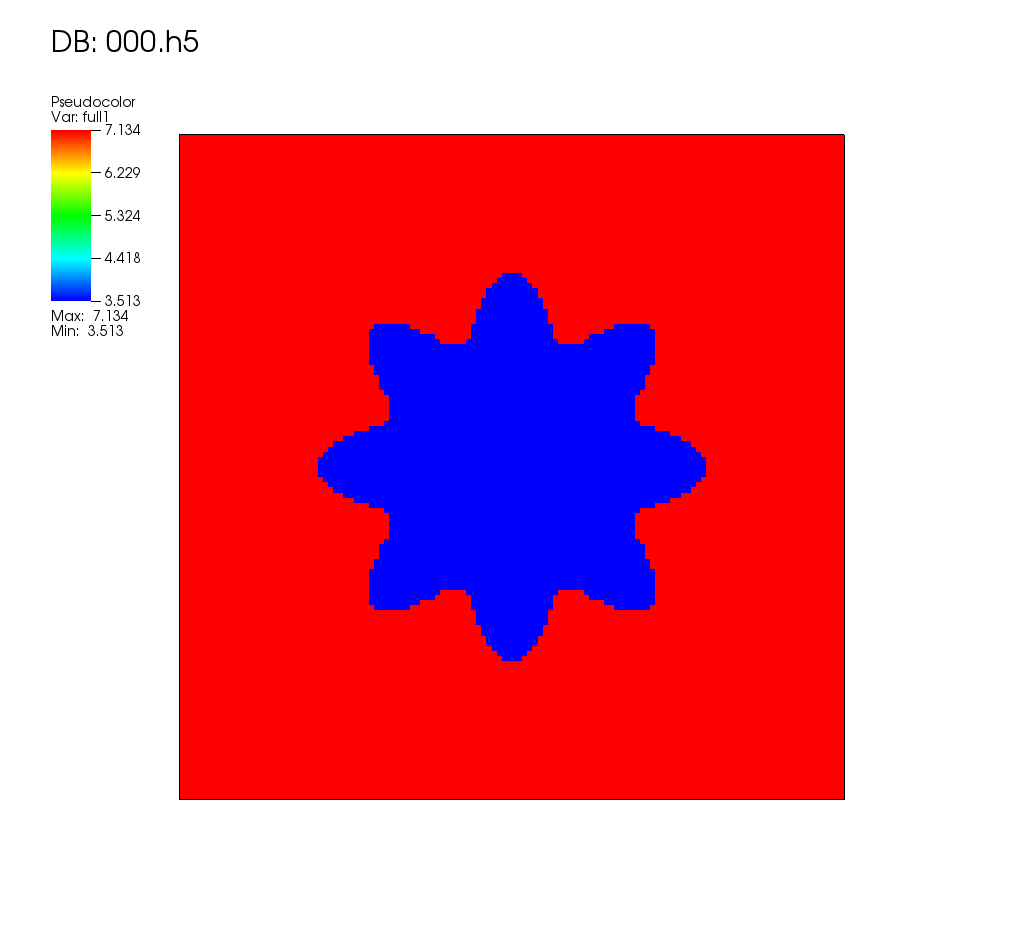}}
\caption{Initial conditions of two components in gas-liquid mixture}
\label{fig:Initial_cd}
\end{figure}

\begin{figure}
\centering
\subfigure[ $n_1$ at t = 1  ]{\includegraphics[width=0.225\textwidth]{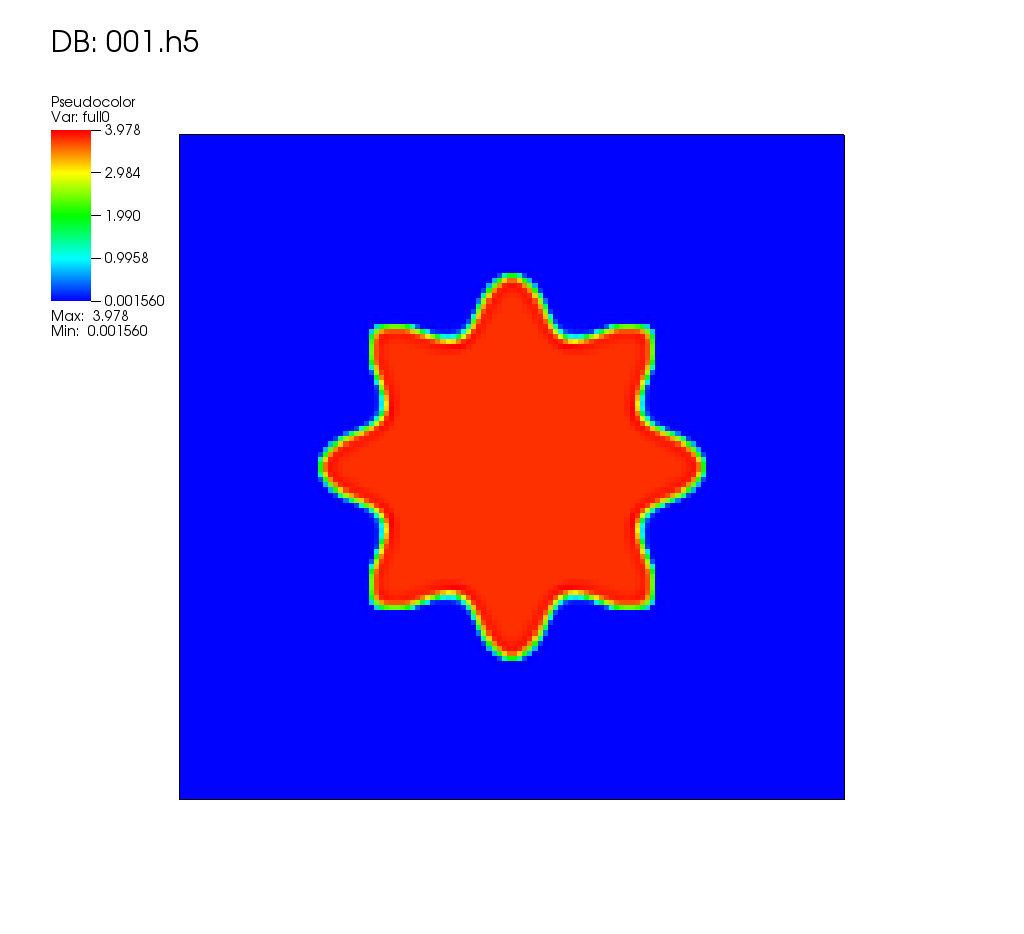}}
\subfigure[ $n_1$ at t = 3  ]{\includegraphics[width=0.225\textwidth]{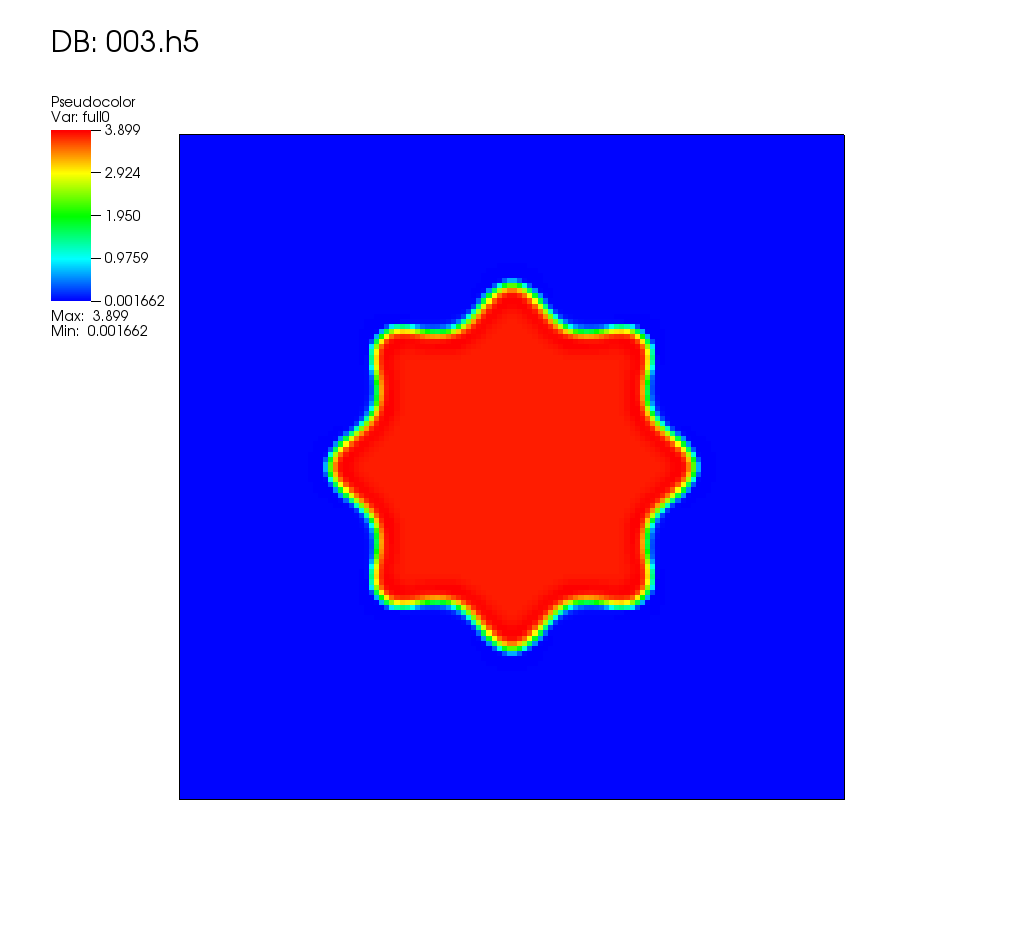}}
\subfigure[ $n_1$ at t = 5  ]{\includegraphics[width=0.225\textwidth]{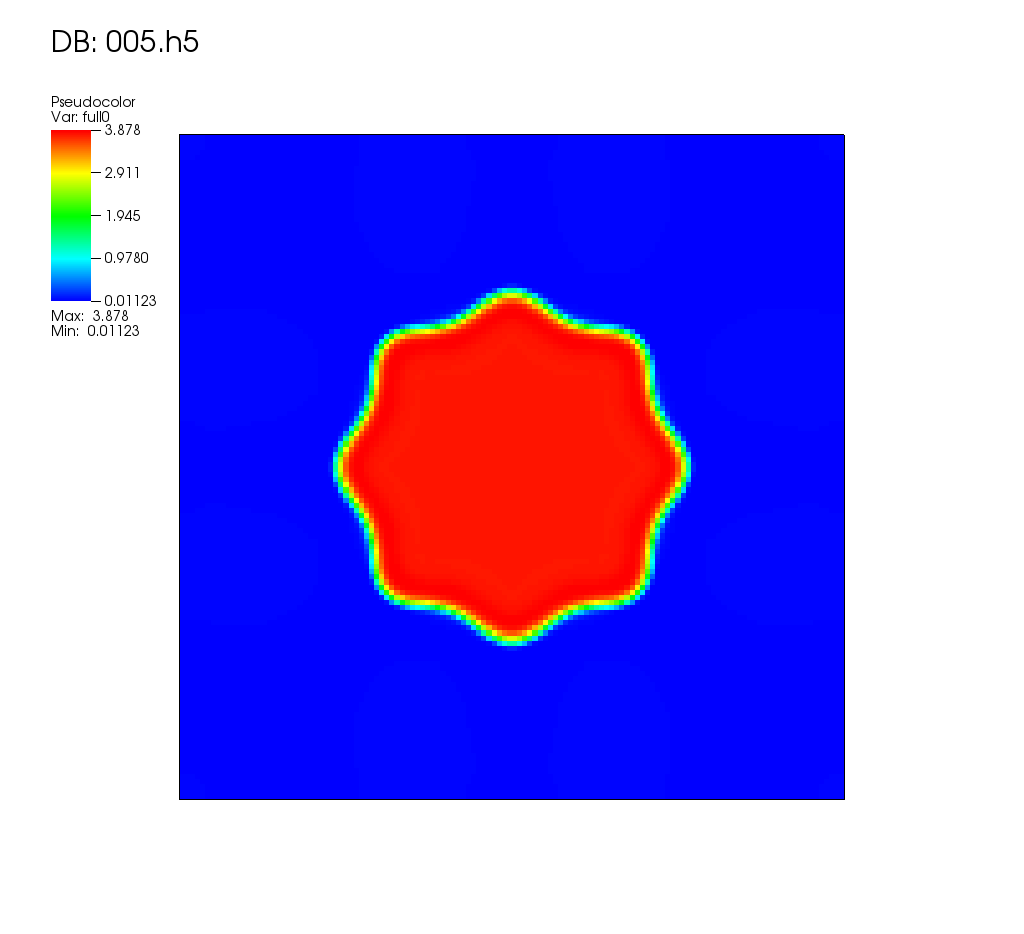}}
\subfigure[ $n_1$ at t = 6000  ]{\includegraphics[width=0.225\textwidth]{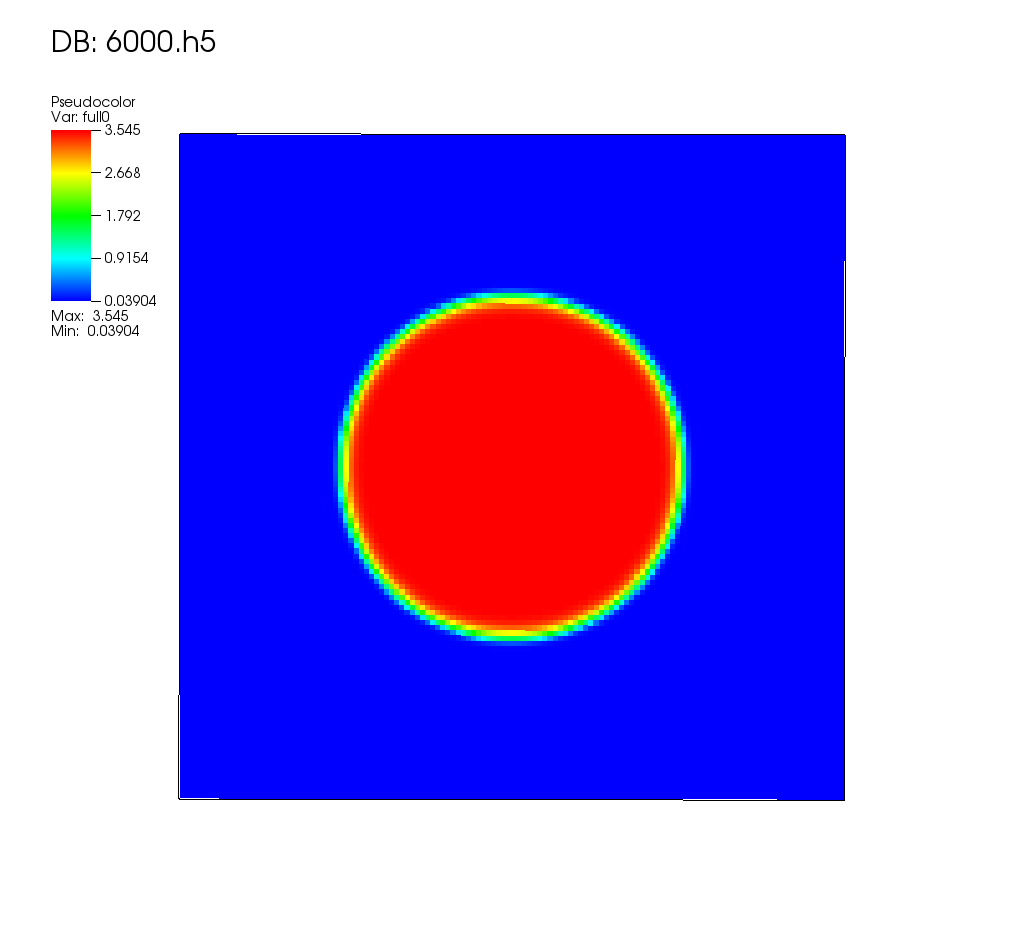}}\\
\subfigure[ $n_2$ at t = 1  ]{\includegraphics[width=0.225\textwidth]{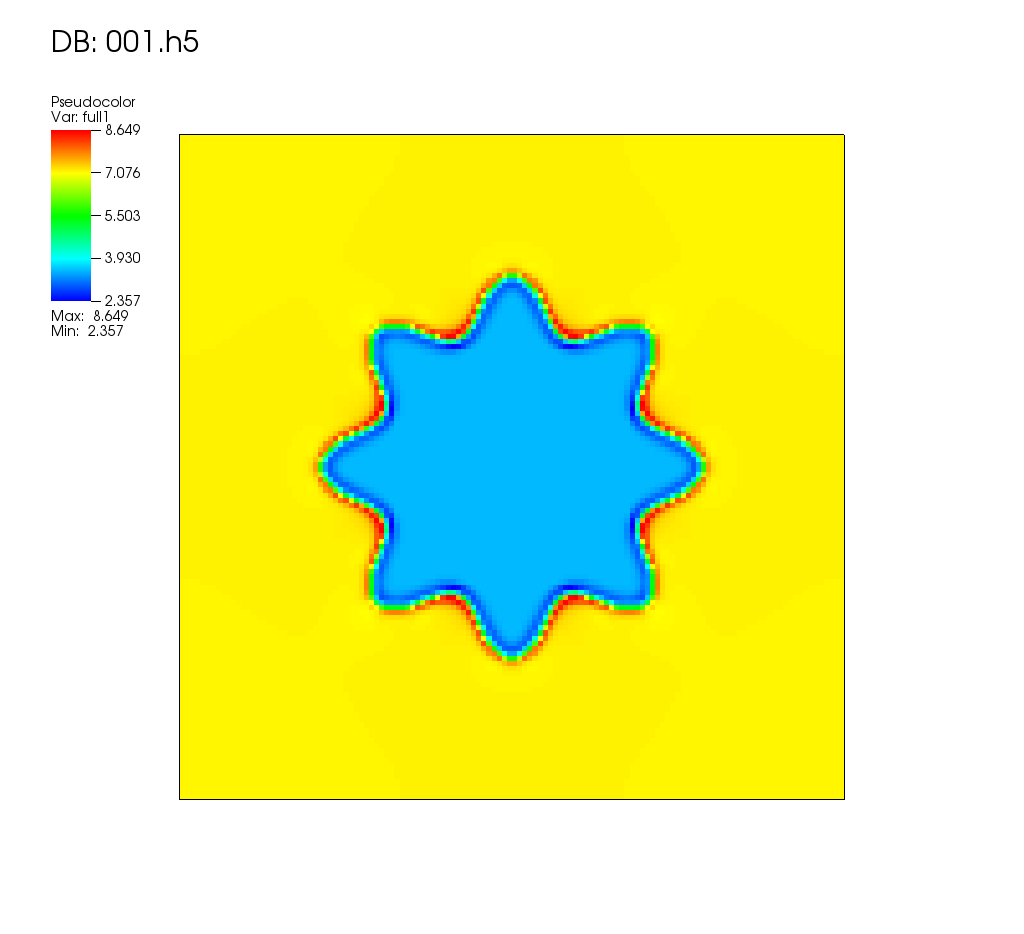}}
\subfigure[ $n_2$ at t = 3 ]{\includegraphics[width=0.225\textwidth]{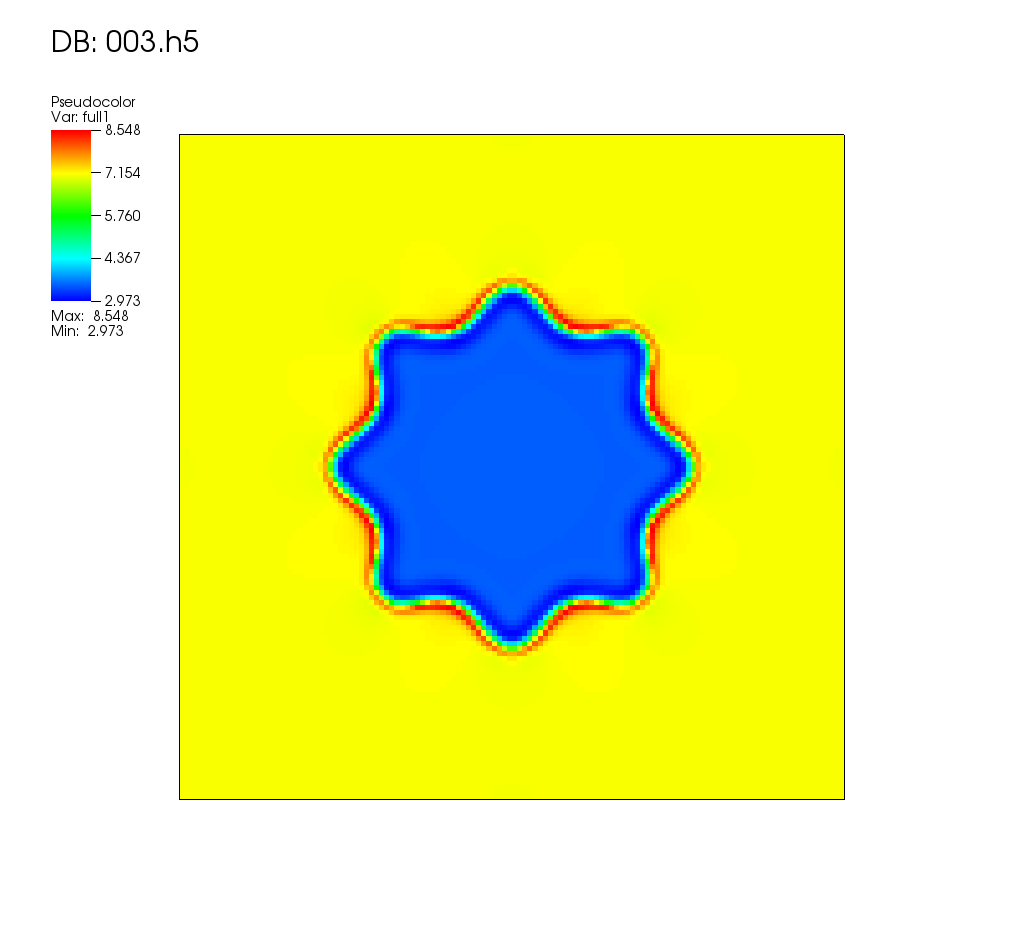}}
\subfigure[ $n_2$ at t = 5  ]{\includegraphics[width=0.225\textwidth]{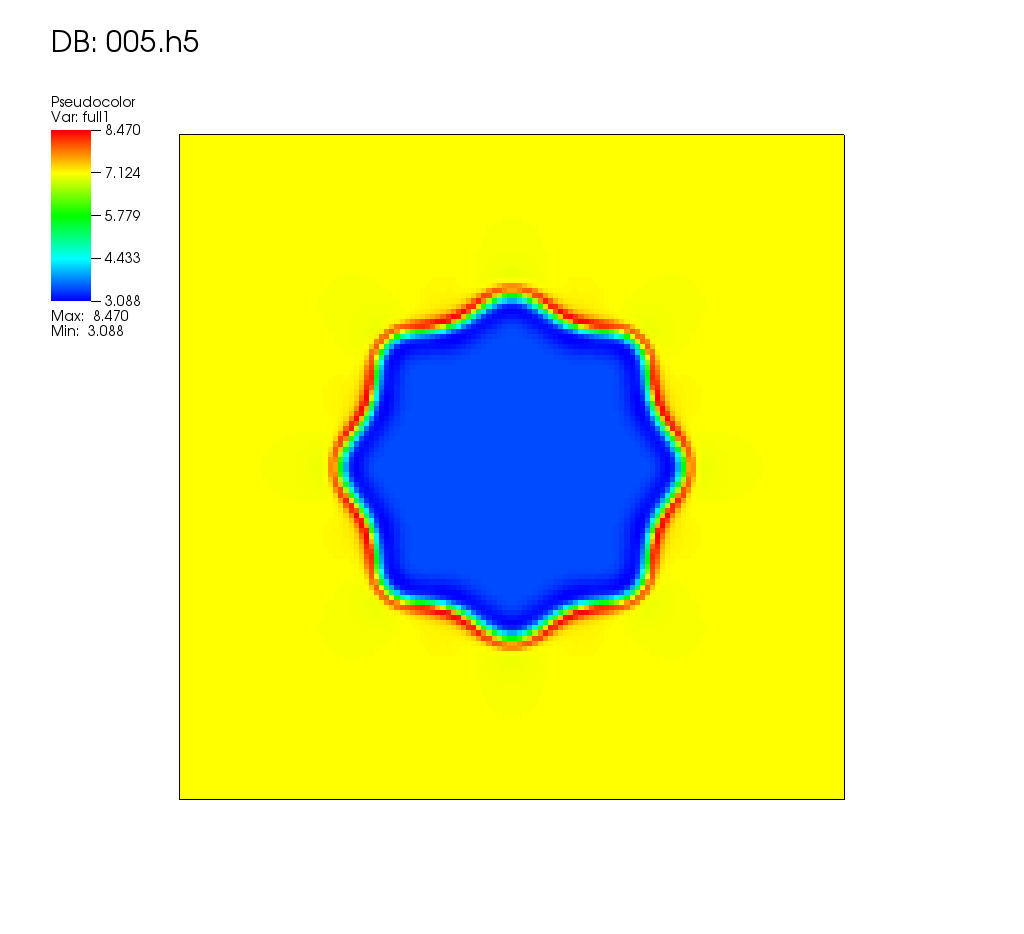}}
\subfigure[ $n_2$ at t = 6000 ]{\includegraphics[width=0.225\textwidth]{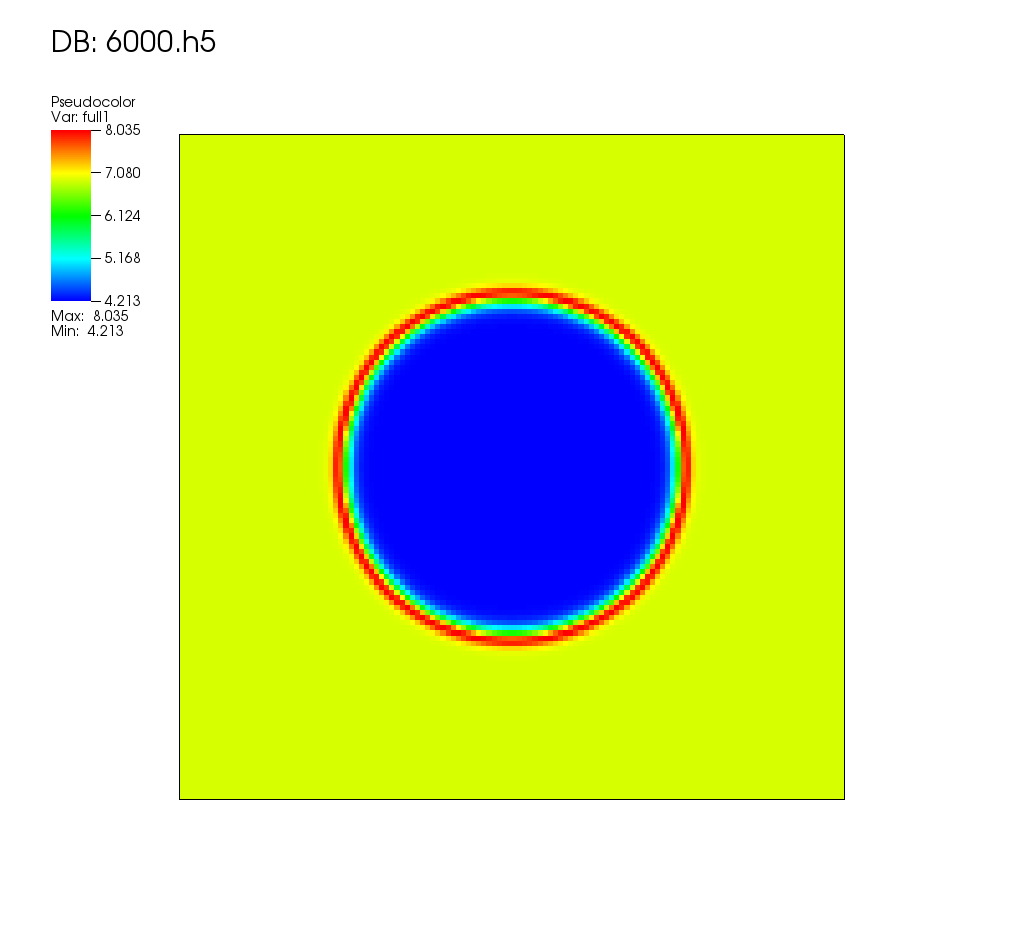}}\\
\subfigure[ ${\bf v}$ at t = 1  ]{\includegraphics[width=0.225\textwidth]{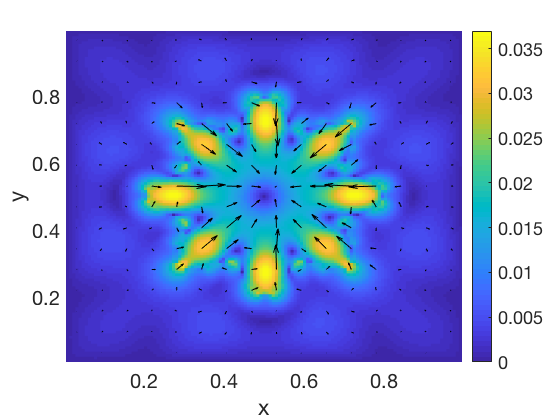}}
\subfigure[ ${\bf v}$ at t = 3  ]{\includegraphics[width=0.225\textwidth]{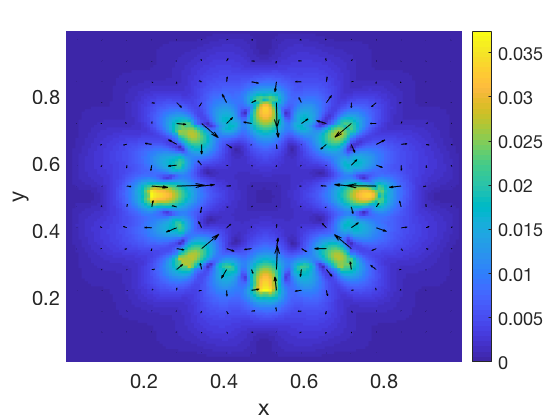}}
\subfigure[ ${\bf v}$ at t = 5  ]{\includegraphics[width=0.225\textwidth]{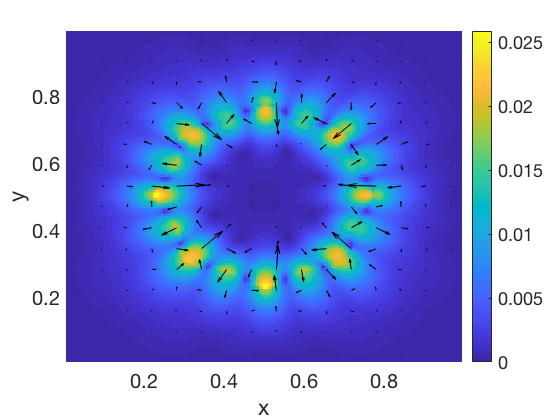}}
\subfigure[ ${\bf v}$ at t = 6000  ]{\includegraphics[width=0.225\textwidth]{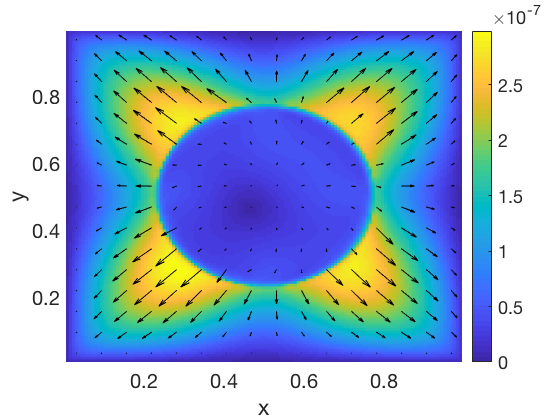}}
\caption{(a-d) Snapshots of $n_1$ at t = 1, 3, 5, 6000. (e-h) Snapshots of $n_2$. (i-l) The corresponding velocity fields.}
\label{fig:molar_density1}
\end{figure}

We consider a mixture of methane and n-decane in a square domain with the length of 80 nm on each side. We denote the molar density of n-decane as $n_1$ and that of methane as $n_2$, respectively. In Table \ref{table:parameters}, we list the dimensional parameters related to these two components. Other parameter values \cite{Cullick_Mathis1984} are chosen as follows
\ben\bea{l}
\eta_1 = \eta_2 = 1 \times 10^{-4} Pa \cdot s, \quad \overline{\eta}_1 = \overline{\eta}_2 = 0.33 \times 10^{-4} Pa \cdot s, \quad M_1 = 1 \times 10^{-12} m^2 \cdot s^{-1},\\
\kappa_{n_1 n_1} = 1.1246 \times 10^{-18}, \quad \kappa_{n_2 n_2} = 2.8649 \times 10^{-20}, \quad \kappa_{n_1 n_2} = 8.9748 \times 10^{-20}.
\eea\een
The gas constant is $R = 8.3144598 J \cdot mol^{-1} \cdot K^{-1}$, the temperature $T = 330 K$.
\begin{table}
\caption{Dimensional critical parameters}
\begin{center}
\begin{tabular}{c|c|c|c|c} \hline
Symbol &  $T_c(K)$ & $P_c(MPa)$ & w & m (kg $\cdot mol^{-1}$)\\ \hline
 n-decane ($C_{10}H_{22}$)&  617.7 & 2.103 & 0.4884 & 0.14228    \\ \hline
methane ($CH_4$)&  190.564 & 4.5992 & 0.01142 & 0.0160428      \\ \hline
 \end{tabular}
\end{center}
\label{table:parameters}
\end{table}\\
The initial conditions are given by
\ben\bea{l}
n_i = \begin{cases}
n_i^{liquid}, \qquad (x^2 + y^2) \leq (r_1 + r_2 \times cos(n \times arctan(\frac{x}{y})))^2 \quad in \quad [-4 \times 10^{-8} m, 4 \times 10^{-8} m]^2, \\
n_i^{gas}, \qquad otherwise \quad in \quad [-4 \times 10^{-8} m, 4 \times 10^{-8} m]^2,
\end{cases}
\eea\een
where $r_1 = 1$, $r_2 = 0.2$, $n = 8$ and
\ben\bea{l}
n_1^{liquid} = 3814.6 mol \cdot m^{-3}, \qquad n_1^{gas} = 26.5 mol \cdot m^{-3},\\
n_2^{liquid} = 3513.2 mol \cdot m^{-3}, \qquad n_2^{gas} = 7133.9 mol \cdot m^{-3}.
\eea\een
If we take characteristic molar density $n_0 = 10^3 mol \cdot m^{-3}$, characteristic density $\rho_0 = n_0 m_2 = 16.0428 kg \cdot m^{-3}$, characteristic length $h = 2 \times 10^{-8} m$, characteristic time $t_0 = 6.4171 \times 10^{-11}s$, and characteristic temperature $T_0 = 273 K$, we obtain dimensionless parameter values as follows
\ben\bea{l}
Re_{1s} = Re_{2s} =  1, \quad Re_{1v} = Re_{2v} = 3, \quad M_1 = 9.7136 \times 10^{-4},\\
\kappa_{n_1 n_1} = 0.0018, \quad \kappa_{n_2 n_2} = 4.5961 \times 10^{-5}, \quad \kappa_{n_1 n_2} = 1.4398 \times 10^{-4} .
\eea\een
Other dimensionless critical parameters of the methane and n-decane are given in table \ref{table:parameters2}. Through the non-dimensionalization, the  gas constant $R$ results in a constant $R_0 = 1.4566$, the dimensionless temperature $T = 1.2088$.
\begin{table}
\caption{Dimensionless critical parameters}
\begin{center}
\begin{tabular}{c|c|c|c|c} \hline
Symbol &  $T_c $ & $P_c $ & w & m  \\ \hline
 n-decane ($C_{10}H_{22}$)&  2.2626 & 1.3495 & 0.4884 & 8.8688    \\ \hline
methane ($CH_4$)&  0.6980 & 2.9513 & 0.01142 & 1      \\ \hline
 \end{tabular}
\end{center}
\label{table:parameters2}
\end{table}
The corresponding dimensionless initial conditions become
\ben\bea{l}
n_i = \begin{cases}
n_i^{liquid}, \qquad (x^2 + y^2) \leq (r_1 + r_2 \times cos(n \times arctan(\frac{x}{y})))^2 \quad in \quad [-2, 2] \times [-2, 2], \\
n_i^{gas}, \qquad otherwise \quad in \quad [-2, 2] \times [-2, 2],
\end{cases}
\eea\een
where $r_1 = 1$, $r_2 = 0.2$, $n = 8$ and
\ben\bea{l}
n_1^{liquid} = 3.8146, \quad n_1^{gas} = 0.0265 ,\quad
n_2^{liquid} = 3.5132 , \quad n_2^{gas} = 7.1339 .
\eea\een
Shown in Figure \ref{fig:Initial_cd}, we perturb the initial condition with certain roughness on the interface, which is unstable due to the surface tension. As time elapses, the roughness vanishes, leading to a surface with the minimal surface tension on it, shown in Figure (\ref{fig:peng-robinson_energy1}-b). The corresponding time evolution of velocities are depicted in Figure \ref{fig:molar_density1}, which show that hydrodynamics indeed speed up the evolution of the system to the steady states.

\subsubsection{Density profiles and mass absorption at the interface in equilibrium}

\begin{figure}
\centering
\subfigure[ Density profiles at the equilibrium state (t = 6000) ]{\includegraphics[width=0.3\textwidth]{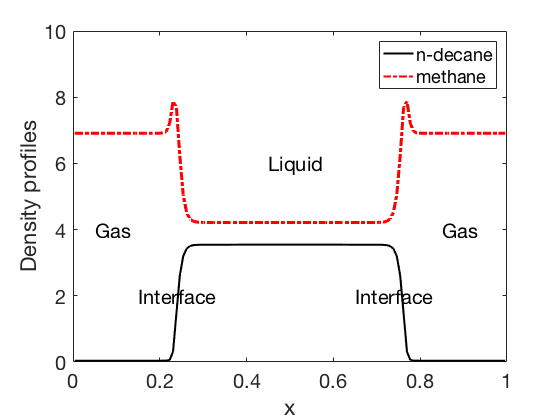}}
\subfigure[ Free energy contour at the equilibrium state (t = 6000) ]{\includegraphics[width=0.3\textwidth]{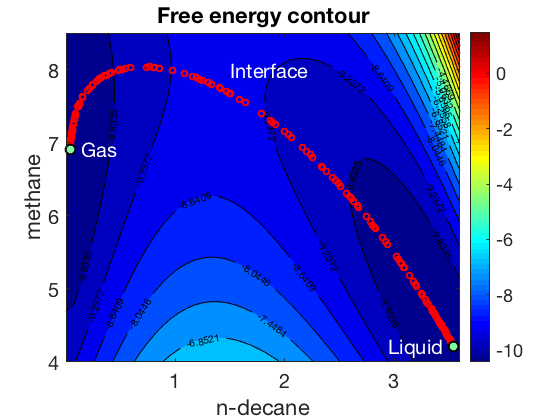}}
\caption{(a) Density profiles of n-decane and methane (y = 0) at the equilibrium state; (b) Free energy contour. Green points represent the densities of n-decane and methane at bulk area and red circles represent their densities on the interface at equilibrium state.}
\label{fig:density_profiles_surface}
\end{figure}

\noindent \indent Near equilibrium ($t=6000$), we show the density profiles of the two fluid components at $y = 0$  in Figure (\ref{fig:density_profiles_surface}-a) and observe  mass absorption of methane at the interface. At the equilibrium of co-existing phases, two (or more) bulk phases have equal chemical potentials, i.e. the corresponding bulk free energies lie on the same tangent line (or surface). For the  Peng-Robinson free energy, it is not straightforward to find the equilibrium states by observing the graph of the free energy function directly. Following the work reported in \cite{Row89, MU2017118}, we subtract the tangent line (or surface) from the Helmoholtz free energy density function to make the equilibrium states as the minimum points, which are then easily observed,
\ben\bea{l}
h_m ({\bf n}, T) = h({\bf n}, T) - \sum_{i=1}^2 \mu_i^0 n_i,
\eea\een
where $\mu_i^0, i = 1, 2$ represent the chemical potential of the ith component at the bulk equilibrium state. We show the modified free energy contour in Figure (\ref{fig:density_profiles_surface}-b). The circled curve represents the energy path of density profiles at the equilibrium state. To avoid high free energy, n-decane and methane change from one equilibrium state (Gas) to another equilibrium state (Liquid) through the saddle point of the free energy surface. Thus, the methane has a higher density on the interface than in the bulk states, leading to the mass absorption phenomena at the interface.

\begin{figure}
\centering
\subfigure[ Total Energy ]{\includegraphics[width=0.3\textwidth]{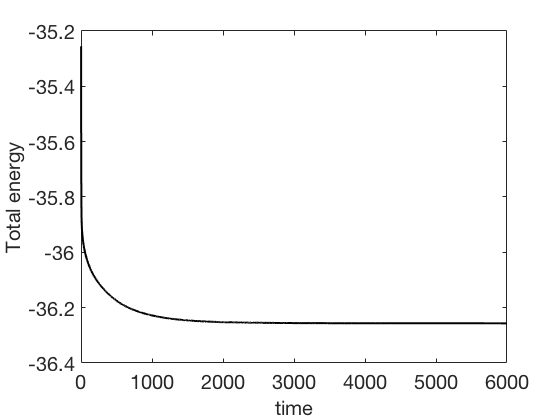}}
\subfigure[ Surface tension ]{\includegraphics[width=0.3\textwidth]{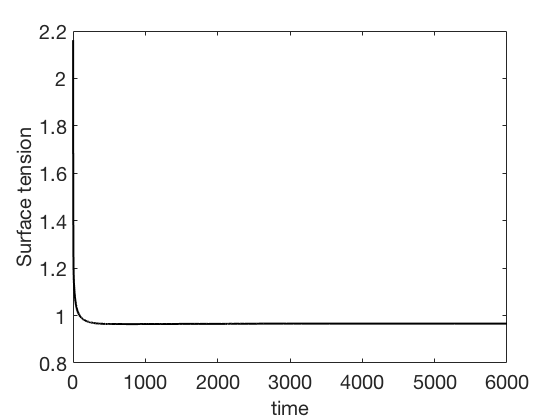}}\\
\subfigure[ The difference of total mass of the component 2 with its initial total mass ]{\includegraphics[width=0.3\textwidth]{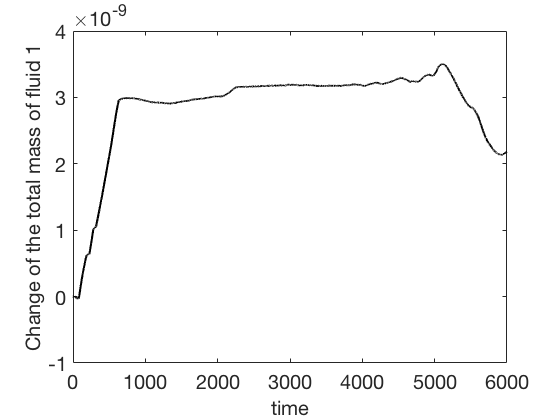}}
\subfigure[ The difference of total mass of the component 2 with its initial total mass ]{\includegraphics[width=0.3\textwidth]{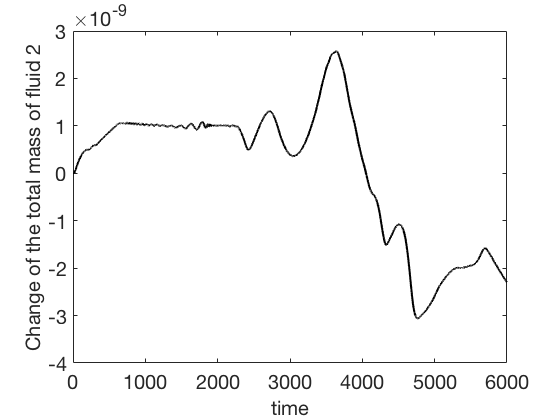}}
\caption{(a) Total energy of the system (\ref{eq:nondim_mole}) with the Peng-Robinson bulk free energy \eqref{eq:Peng_Robinsoon_total_energy}; (b) Surface tension of the mixture; (c, d) Total mass of the component 1 and 2 on the rectangular domain $\Omega = [-2, 2] \times [-2, 2]$, solved in the system (\ref{eq:nondim_mole}) with the Peng-Robinson bulk free energy \eqref{eq:Peng_Robinsoon_total_energy}. (e) Density profiles of n-decane and methane (y = 0) at the equilibrium state; (f) Free energy contour. Green points represent the densities of n-decane and methane at bulk area and red circles represent their densities on the interface at equilibrium state. }
\label{fig:peng-robinson_energy1}
\end{figure}

The total energy and total mass difference with the initial condition for each component are shown in Figure \ref{fig:peng-robinson_energy1}, which verifies  energy stability and mass conservation of our numerical scheme.

This numerical experiment not only demonstrates that our mathematical model can be applied to study  thermodynamic and hydrodynamic properties of the fluid mixture in an application relevant to the petroleum industry, but also showcases that our numerical scheme can handle the Navier-Stokes-Cahn-Hilliard equation system with a highly nonlinear free energy \eqref{eq:Peng_Robinsoon_total_energy}.



\section{Conclusion}
In this paper, we present a second order, fully-discrete, linear and unconditionally energy stable numerical scheme for the hydrodynamic phase field model of compressible fluid flow. Firstly, we reformulate the  model by introducing a couple of intermediate variables, based on the Energy Quadratization approach. Using the reformulated model equations, we develop a second order, energy stable, semi-discrete numerical scheme in time. Then,  we obtain a fully discrete numerical scheme applying the finite difference method on the staggered grid, which preserves a fully discrete energy dissipation law. In addition, the well-posedness of the linear system resulting from the linear numerical scheme is  proved rigorously. Several numerical experiments are presented to verify the accuracy, stability and efficiency of our numerical scheme. The comparison between the simulations with and without hydrodynamics is used to demonstrate the mixing role played by hydrodynamics in  phase separation phenomena in binary compressible fluid flows. The scheme can be readily extended to models  N-component compressible fluid flows with $N>2$.

\section{Appendix}
\subsection{Linear system resulting from the numerical scheme}

We summarize the linear system resulting from the numerical scheme as follows.
\ben\bea{l}
\begin{cases}
\big \{2 \frac{\rho_1}{\Delta t}  + d_x (A_x(   \overline{\rho}_1^{n + 1/2}    \overline{\frac{1}{\sqrt{\rho}}}^{n+1/2}    ) u) +  d_y (A_y(\overline{\rho}_1^{n + 1/2}    \overline{\frac{1}{\sqrt{\rho}}}^{n+1/2}    ) v) =   \\
M_1 \Delta_h \mu_1-   M_1 \Delta_h  \mu_2 + g_1 \big \}|_{i, j},
 i = 1, \cdots, N_x, j = 1, \cdots, N_y, \\
\\
\big \{2 \frac{\rho_2}{\Delta t} + d_x (A_x(\overline{\rho}_2^{n + 1/2}    \overline{\frac{1}{\sqrt{\rho}}}^{n+1/2}    ) u) +  d_y (A_y(\overline{\rho}_2^{n + 1/2}    \overline{\frac{1}{\sqrt{\rho}}}^{n+1/2}    ) v)=  \\
 - M_1 \Delta_h \mu_1 +  M_1 \Delta_h \mu_2 + g_2 \big \}|_{i, j},
 i = 1, \cdots, N_x, j = 1, \cdots, N_y, \\
\\
\big \{2 \frac{u}{\Delta t}    + \frac{1}{2} (\overline{u}^{n+1/2} D_x( \overline{\frac{1}{\sqrt{\rho}}}^{n+1/2} a_x u) + A_x( \overline{\frac{1}{\sqrt{\rho}}}^{n+1/2} d_x (\overline{ u}^{n+1/2} u))) \\

+ \frac{1}{2}(a_x (A_x \overline{v}^{n+1/2} D_y(A_x( \overline{\frac{1}{\sqrt{\rho}}}^{n+1/2} )u)) + A_x( \overline{\frac{1}{\sqrt{\rho}}}^{n+1/2} )d_y(A_yu A_x(\overline{v}^{n+1/2})) \\
= A_x( \overline{\frac{1}{\sqrt{\rho}}}^{n+1/2} )  (D_x( \frac{1}{Re_s} d_x  (A_x( \overline{\frac{1}{\sqrt{\rho}}}^{n+1/2} ) u)) + d_y( A_x(A_y \frac{1}{Re_s})  D_y  (A_x( \overline{\frac{1}{\sqrt{\rho}}}^{n+1/2} ) u))) \\

+ A_x( \overline{\frac{1}{\sqrt{\rho}}}^{n+1/2} ) D_x( \frac{1}{Re_s} d_x( A_x( \overline{\frac{1}{\sqrt{\rho}}}^{n+1/2} ) u))
+ A_x( \overline{\frac{1}{\sqrt{\rho}}}^{n+1/2} )  d_y( A_x(A_y \frac{1}{Re_s})  D_x( A_y( \overline{\frac{1}{\sqrt{\rho}}}^{n+1/2} ) v)) \\
 + A_x( \overline{\frac{1}{\sqrt{\rho}}}^{n+1/2} ) D_x( \frac{1}{Re_v} d_x( A_x( \overline{\frac{1}{\sqrt{\rho}}}^{n+1/2} ) u)) + A_x( \overline{\frac{1}{\sqrt{\rho}}}^{n+1/2} )  D_x( \frac{1}{Re_v} d_y( A_y( \overline{\frac{1}{\sqrt{\rho}}}^{n+1/2} ) v))  \\

- A_x(\overline{\rho_1}^{n+1/2}  \overline{\frac{1}{\sqrt{\rho}}}^{n+1/2} ) D_x(\mu_1) -  A_x( \overline{\rho_2}^{n+1/2}  \overline{\frac{1}{\sqrt{\rho}}}^{n+1/2} ) D_x(\mu_2)  + g_3  \big \}|_{i+\frac{1}{2}, j},\\
 i = 1, \cdots, N_x - 1, j = 1, \cdots, N_y ,\\
\\

\big \{ 2 \frac{v}{\Delta t}  + \frac{1}{2} (a_x(A_y \overline{u}^{n+1/2} D_x(A_y( \overline{\frac{1}{\sqrt{\rho}}}^{n+1/2} )v)) + A_y( \overline{\frac{1}{\sqrt{\rho}}}^{n+1/2} )d_x(A_y\overline{u}^{n+1/2} A_xv^{n+1/2})) \\
+ \frac{1}{2} (\overline{v}^{n+1/2} D_y( \overline{\frac{1}{\sqrt{\rho}}}^{n+1/2} a_yv^{n+1/2}) + A_y( \overline{\frac{1}{\sqrt{\rho}}}^{n+1/2} d_y(\overline{v}^{n+1/2}v^{n+1/2}))) \\
=  A_y( \overline{\frac{1}{\sqrt{\rho}}}^{n+1/2} )  (d_x( A_x(A_y \frac{1}{Re_s})  D_x(  A_y( \overline{\frac{1}{\sqrt{\rho}}}^{n+1/2} ) v)) + D_y( \frac{1}{Re_s} d_y ( A_y( \overline{\frac{1}{\sqrt{\rho}}}^{n+1/2} ) v))) \\

+  A_y( \overline{\frac{1}{\sqrt{\rho}}}^{n+1/2}  ) d_x( A_x(A_y \frac{1}{Re_s})  D_y( A_x( \overline{\frac{1}{\sqrt{\rho}}}^{n+1/2} ) u))

+  A_y( \overline{\frac{1}{\sqrt{\rho}}}^{n+1/2}  )  D_y( \frac{1}{Re_s} d_y( A_y( \overline{\frac{1}{\sqrt{\rho}}}^{n+1/2} ) v))  \\
+  A_y( \overline{\frac{1}{\sqrt{\rho}}}^{n+1/2}  )  D_y( \frac{1}{Re_v} d_x( A_x(\overline{\frac{1}{\sqrt{\rho}}}^{n+1/2}) u)) +  A_y( \overline{\frac{1}{\sqrt{\rho}}}^{n+1/2}   ) D_y( \frac{1}{Re_v} d_y( A_y( \overline{\frac{1}{\sqrt{\rho}}}^{n+1/2} ) v)) \\

 - A_y( \overline{\rho_1}^{n+1/2} \overline{\frac{1}{\sqrt{\rho}}}^{n+1/2}  ) D_y(\mu_1) -  A_y( \overline{\rho_2}^{n+1/2}  \overline{\frac{1}{\sqrt{\rho}}}^{n+1/2} ) D_y(\mu_2)  + g_4  \big \}|_{i, j +\frac{1}{2}},\\ i = 1, \cdots, N_x, j = 1, \cdots, N_y - 1.
\end{cases}
\eea \label{eq:full_discrete_0}
\een

\ben\bea{l}
\begin{cases}
\big \{ 4 \frac{q_1}{\Delta t}= 4 \overline{{\frac{\partial q_1}{\partial \rho_1}}}^{n+1/2} \frac{\rho_1}{\Delta t} + 4 \overline{\frac{\partial q_1}{\partial \rho_2}}^{n+1/2}  \frac{\rho_2}{\Delta t} + g_5 \big \} |_{i,j}, i = 1, \cdots, N_x , j=1, \cdots, N_y,\\
\\
\big \{ -  \frac{2}{\Delta t} \mu_1=  - 4q_1\frac{1}{\Delta t} \overline{\frac{\partial q_1}{\partial \rho_1}}^{n+1/2} +\frac{2}{\Delta t} \kappa_{\rho_1 \rho_1} \Delta_h \rho_1 + \frac{2}{\Delta t}\kappa_{\rho_1 \rho_2} \Delta_h \rho_2 + g_6 , \big \}|_{i, j},\\    i = 1, \cdots, N_x, j = 1, \cdots, N_y, \\
\\
\big \{ - \frac{2}{\Delta t} \mu_2 =  - 4q_1\frac{1}{\Delta t}  \overline{\frac{\partial q_1}{\partial \rho_2}}^{n+1/2} +\frac{2}{\Delta t} \kappa_{\rho_2 \rho_2} \Delta_h \rho_2 +\frac{2}{\Delta t} \kappa_{\rho_1 \rho_2} \Delta_h \rho_1 + g_7, \big \}|_{i, j}, \\    i = 1, \cdots, N_x, j = 1, \cdots, N_y,
\end{cases}
\eea \label{eq:full_discrete_00}
\een
where $\rho_i, \mu_i, i=1, 2$ and $q_1$ satisfy discrete homogeneous Neumann boundary conditions \eqref{eq:neumann_discrete}, $u, v$ the discrete homogeneous Dirichlet boundary conditions \eqref{eq:dirichlet_discrete}.
We define ${\bf D}_h$ as
\ben\bea{l}
 \left(
\bea{cc}
d_x  (A_x( \overline{\frac{1}{\sqrt{\rho}}}^{n+1/2} ) u))  & \frac{1}{2} (D_x( A_y( \overline{\frac{1}{\sqrt{\rho}}}^{n+1/2} ) v)  +  D_y  (A_x( \overline{\frac{1}{\sqrt{\rho}}}^{n+1/2} ) u))) \\
\frac{1}{2} (D_x( A_y( \overline{\frac{1}{\sqrt{\rho}}}^{n+1/2} ) v)  +  D_y  (A_x( \overline{\frac{1}{\sqrt{\rho}}}^{n+1/2} ) u)))  & d_y ( A_y( \overline{\frac{1}{\sqrt{\rho}}}^{n+1/2} ) v))
\eea
\right)
\eea \label{eq:D_h_definition}
\een

\bibliographystyle{plain}
\bibliography{xp_thesis}{}

\end{document}